\documentclass{amsart}

\usepackage{amssymb}
\usepackage{graphicx}

\newtheorem{theorem}{Theorem}[section]
\newtheorem{lemma}[theorem]{Lemma}
\newtheorem{proposition}[theorem]{Proposition}
\newtheorem{corollary}[theorem]{Corollary}

\newtheorem{_conjecture}[theorem]{Conjecture}

\newtheorem{_problem}[theorem]{Problem}
\newenvironment{problem}{\begin{_problem}\rm}{\end{_problem}}

\newtheorem{_algorithm}[theorem]{Algorithm}
\newenvironment{algorithm}{\begin{_algorithm}\rm}{\end{_algorithm}}

\newtheorem{_subroutime}[theorem]{Subroutine}
\newenvironment{subroutine}{\begin{_subroutime}\rm}{\end{_subroutime}}

\newtheorem{_claim}[theorem]{Claim}
\newenvironment{claim}{\begin{_claim}\rm}{\end{_claim}}

\newtheorem{_subclaim}[theorem]{Sub-claim}

\newtheorem{_definition}[theorem]{Definition}
\newenvironment{definition}{\begin{_definition}\rm}{\end{_definition}}

\newtheorem{_remark}[theorem]{\it Remark}
\newenvironment{remark}{\begin{_remark}\rm}{\end{_remark}}

\newtheorem{_example}[theorem]{Example}
\newenvironment{example}{\begin{_example}\rm}{\end{_example}}

\newtheorem*{maintheorem}{Main Theorem}

\numberwithin{equation}{section}
\numberwithin{table}{section}
\numberwithin{figure}{section}

\newcommand{\F}{\mathord{\mathbb F}}
\renewcommand{\P}{\mathord{\mathbb  P}}
\newcommand{\Q}{\mathord{\mathbb  Q}}
\newcommand{\R}{\mathord{\mathbb R}}
\newcommand{\Z}{\mathord{\mathbb Z}}

\newcommand{\AAA}{\mathord{\mathcal A}}

\newcommand{\CCC}{\mathord{\mathcal C}}

\newcommand{\EEE}{\mathord{\mathcal E}}

\newcommand{\III}{\mathord{\mathcal I}}

\newcommand{\LLL}{\mathord{\mathcal L}}
\newcommand{\MMM}{\mathord{\mathcal M}}

\newcommand{\OOO}{\mathord{\mathcal O}}

\newcommand{\UUU}{\mathord{\mathcal U}}
\newcommand{\VVV}{\mathord{\mathcal V}}

\newcommand{\ZZZ}{\mathord{\mathcal Z}}

\font\mathgot=eufm10

\newcommand{\SSSS}{\mathord{\hbox{\mathgot S}}}
\newcommand{\moduli}{\mathord{\hbox{\mathgot M}}}

\newcommand{\maprightsp}[1]{\; \smash{\mathop{\; \longrightarrow \; }\limits\sp{#1}}\; }

\newcommand{\mapdown}{\phantom{\Big\downarrow}\hskip -8pt \downarrow}

\newcommand{\mapdownsurj}{
\hbox{$\bigm\downarrow$}
\llap{\hbox{\raise 2pt\hbox{$\bigm\downarrow$}}}%
\vstrechmapdown
}

\newcommand{\inj}{\hookrightarrow}

\newcommand{\isom}{\smash{\mathop{\;\to\;}\limits\sp{\sim\,}}}

\newcommand{\set}[2]{\{\; {#1} \; \mid \; {#2} \;  \}}
\newcommand{\sset}[1]{\{#1\}}
\newcommand{\shortset}[2]{\{ {#1} \mid  {#2}   \}}

\newcommand{\map}[3]{ #1 \, : \, #2 \, \to \, #3}
\newcommand{\mapisom}[3]{ #1 \, : \, #2 \; \isom \; #3}

\newcommand{\shortmap}[3]{ #1 : #2 \to #3}

\newcommand{\sm}{\setminus}
\newcommand{\st}{\subset}

\newcommand{\sprime}{\sp\prime}

\newcommand{\spar}[1]{\sp{(#1)}}
\newcommand{\spprime}{\sp{\prime\prime}}

\newcommand{\sptimes}{\sp{\times}}
\newcommand{\sperp}{\sp{\perp}}
\newcommand{\dual}{\sp{\vee}}

\newcommand{\inv}{\sp{-1}}

\newcommand{\pione}{\pi\sb 1}

\renewcommand{\Im}{\operatorname{\rm Im}\nolimits}

\newcommand{\pr}{\operatorname{\rm pr}\nolimits}

\newcommand{\GL}{\operatorname{\it GL}\nolimits}
\newcommand{\PGL}{\operatorname{\it PGL}\nolimits}

\newcommand{\Aut}{\operatorname{\rm Aut}\nolimits}
\newcommand{\Hom}{\operatorname{\rm Hom}\nolimits}

\newcommand{\Sing}{\operatorname{\rm Sing}\nolimits}
\newcommand{\Spec}{\operatorname{\rm Spec}\nolimits}
\newcommand{\SSpec}{\operatorname{\bf Spec}\nolimits}
\newcommand{\Proj}{\operatorname{\rm Proj}\nolimits}
\newcommand{\rank}{\operatorname{\rm rank}\nolimits}
\newcommand{\disc}{\operatorname{\rm disc}\nolimits}
\newcommand{\Frob}{\operatorname{\rm Frob}\nolimits}

\newcommand{\Der}[2]{\frac{\partial #1}{\partial #2}}
\newcommand{\der}[2]{{\partial #1}/{\partial #2}}

\newcommand{\rmand}{\textrm{and}}

\newcommand{\quand}{\quad\rmand\quad}

\newcommand{\DK}{\mathord{\rm DK}}

\newcommand{\ZdG}{Z(dG)}

\newcommand{\Hz}[1]{H\sp 0 (\P\sp 2 , #1)}
\newcommand{\hz}[1]{h\sp 0 (\P\sp 2 , #1)}

\newcommand{\Uts}{\,\UUU\sb{2, 6}}
\newcommand{\Pt}{\P\sp 2}
\newcommand{\card}[1]{|#1|}
\newcommand{\gen}[1]{\langle #1\rangle}

\newcommand{\OPt}{\OOO\sb{\Pt}}
\newcommand{\Ob}{\Omega (b)}
\newcommand{\length}{\operatorname{\rm length}\nolimits}
\newcommand{\IZs}{\III\sb{Z(s)}}
\newcommand{\IZdG}{\III\sb{Z(dG)}}
\newcommand{\XG}{X\sb G}
\newcommand{\YG}{Y\sb G}
\newcommand{\KG}{K\sb G}
\newcommand{\piG}{\pi\sb G}
\newcommand{\phiG}{\phi\sb G}
\newcommand{\SG}{S\sb G}
\newcommand{\SzG}{S\sp0\sb G}

\newcommand{\Upb}{\,\UUU\sb{p, b}}
\newcommand{\Utb}{\,\UUU\sb{2, b}}
\newcommand{\tUtb}{\,\widetilde{\UUU}\sb{2, b}}
\newcommand{\tUts}{\,\widetilde{\UUU}\sb{2, 6}}
\newcommand{\tU}{\,\widetilde{\UUU}}
\newcommand{\Vpb}{\VVV\sb{p, b}}
\newcommand{\Vtb}{\VVV\sb{2, b}}
\newcommand{\Vts}{\VVV\sb{2, 6}}
\newcommand{\Ample}{H}
\newcommand{\AG}{\Ample\sb G}
\newcommand{\DG}{\operatorname{\rm DG}\nolimits}
\newcommand{\ctt}{\tt}
\newcommand{\tZ}{ {\mathord{\ctt Z}}}
\newcommand{\tP}{ {\mathord{\ctt P}}}
\newcommand{\tQ}{ {\mathord{\ctt Q}}}
\newcommand{\tS}{ {\mathord{\ctt S}}}
\newcommand{\sA}{ {\mathord{\mathbf A}}}
\newcommand{\te}{ {\mathord{\hbox{\hskip .7pt \ctt e}}}}
\newcommand{\tth}{ {\mathord{\hbox{\hskip .7pt  \ctt h}}}}
\newcommand{\tC}{ {\mathord{\ctt C}}}
\newcommand{\ttC}{ \widetilde{\tC}}
\newcommand{\lift}{\sp{\sim}}
\newcommand{\Pow}{\operatorname{\rm Pow}\nolimits}
\newcommand{\tCCG}{\widetilde{\CCC}\sb G}
\newcommand{\CCG}{\CCC\sb G}
\newcommand{\PT}[1]{D\sb{#1}} %proper transform
\newcommand{\HPT}[1]{F\sb{#1}} %half of proper transform

\newcommand{\wG}{w\sb G}
\newcommand{\Nodes}{\mathord{\it N}}
\newcommand{\res}{\mathord{\rm res}}
\newcommand{\wtC}{\widetilde C}
\newcommand{\wt}{\operatorname{\rm wt}\nolimits}
\newcommand{\wtenum}{\operatorname{\rm wtenum}\nolimits}
\newcommand{\tomg}{\tilde\omega}
\newcommand{\omg}{\omega}
\newcommand{\WT}{\mathord{\rm WT}}
\newcommand{\NS}{\mathord{{\it NS}}}
\newcommand{\Roots}{\mathord{\rm Roots}}
\newcommand{\vb}{\mathord{\bf b}}
\newcommand{\Sn}{\SSSS\sb n}
\newcommand{\St}{\SSSS\sb {21}}
\newcommand{\LL}{\mathord{\bf L}}
\newcommand{\RSts}{\Lambda\sb{2, \sigma}}
\newcommand{\YH}{Y\sb{|H|}}
\newcommand{\Bs}{\mathord{\rm Bs}}
\newcommand{\tGen}{\mathord{\tt Gen}}
\newcommand{\tCol}{\mathord{\tt Col}}
\newcommand{\CCL}[1]{\mathord{\bf C}\sb{#1}}

\begin{document}

\title[Supersingular $K3$ surfaces]{Supersingular  $K3$ surfaces in
characteristic $2$ \\
as double covers of a projective
plane}

\author{Ichiro Shimada}
\address{
Department of Mathematics,
Faculty of Science,
Hokkaido University,
Sapporo 060-0810,
JAPAN
}
\email{shimada@math.sci.hokudai.ac.jp
}


\begin{abstract}
For every supersingular $K3$ surface $X$ in characteristic $2$,
there exists a homogeneous polynomial $G$ of degree $6$
such that $X$ is birational to the purely inseparable double cover of $\Pt$
defined by $w^2=G$.
We present  an algorithm to calculate from $G$ a set of generators of the numerical N\'eron-Severi lattice of $X$.
As an application, we investigate
the stratification defined by the Artin invariant on a moduli space of supersingular $K3$ surfaces of degree $2$
in characteristic $2$.
\end{abstract}

\maketitle
\section{Introduction}\label{sec:intro}
We work over an algebraically closed field $k$ of characteristic $2$ in Introduction.
\par
In ~\cite{Shimada2003},
we have shown  that every supersingular $K3$ surface 
$X$ in characteristic $2$ is isomorphic to the minimal resolution
$X\sb G$ of a purely inseparable double cover $\YG$ of $\Pt$ defined by
$$
 w^2 =G(X\sb 0, X\sb 1, X\sb 2),
$$
where $G$ is a homogeneous polynomial of degree $6$ such that the singular locus 
$\Sing (\YG)$ of $\YG$ consists of $21$ ordinary nodes.
Conversely,
if $\YG$ has $21$ ordinary nodes
as its only singularities,
then $\XG$ is a supersingular $K3$ surface.
In characteristic $2$,
we can define the  differential $dG$ of 
a homogeneous polynomial $G$  of degree $6$ 
 as a global section  of the vector bundle $\Omega \sp 1 \sb{\P\sp 2} (6)$.
The condition that $\Sing (Y\sb G)$ consists of $21$ ordinary nodes
is equivalent to the condition that the subscheme $\ZdG$ of $\P\sp 2$ defined by $dG=0$
is reduced of dimension $0$.
The homogeneous polynomials  of degree $6$
satisfying this condition  form a Zariski open dense subset $\Uts$ of $\Hz{\OPt(6)}$.
The kernel of the linear homomorphism $G\mapsto dG$ is the linear subspace
$$
\Vts:=\set{H^2\in \Hz{\OPt(6)}}{H\in \Hz{\OPt(3)}} 
$$
of $\Hz{\OPt(6)}$.
If $G\in \Uts$, then $G+H^2\in \Uts$ holds for any $H^2\in \Vts$;
that is,
$\Vts$ acts on $\Uts$ by translation.
Let $G$ and $G\sprime$ be polynomials in $\Uts$.
The supersingular $K3$ surfaces $\XG$ and $X\sb{G\sprime}$ are isomorphic over $\Pt$
if and only if there exist $c\in k\sptimes$ and $H^2\in \Vts$ such that
$$
G\sprime = cG +H^2.
$$
Therefore we can construct a moduli space $\moduli$ of supersingular $K3$ surfaces of degree $2$
in characteristic $2$ by
\begin{equation*}\label{eq:moduli}
\moduli:= \P\sb * (\Uts/\Vts)/ \PGL(3, k).
\end{equation*}
\par
The purpose of this paper is to investigate 
the stratification
of $\Uts$ by the Artin invariant of the supersingular $K3$ surfaces.
Our investigation   yields 
 an algorithm to calculate
a set of generators of 
the numerical N\'eron-Severi
lattice
of  $\XG$ from the homogeneous polynomial $G\in\Uts$.
\par
\bigskip
Suppose that a polynomial  $G$ in  $\Uts$ is given.
The singular points of $\YG$ is mapped bijectively to the points of $\ZdG$
by the covering morphism.
We denote by 
$$
\phi\sb G: X\sb G \to \P\sp 2
$$
the composite of the minimal resolution $X\sb G\to Y\sb G$ 
and the covering morphism $Y\sb G \to \P\sp 2$.
The numerical N\'eron-Severi lattice of the supersingular $K3$ surface $X\sb G$ 
is denoted by $\SG$,
which is a hyperbolic lattice of rank $22$.
Let $\AG\subset X\sb G$ be the pull-back of a general line of $\P\sp 2$ by $\phiG$.
For a point $P\in \ZdG$,
we denote by  $\Gamma\sb {P}$ the $(-2)$-curve on $X\sb G$
that is contracted to $P$ by $\phi\sb G$.
It is obvious  that the sublattice $S\sp 0 \sb G$ of $\SG$
generated by the numerical equivalence classes $[\Gamma\sb P]$ $(P\in \ZdG)$
and $[\AG]$ is of rank $22$, and hence is of finite index in $\SG$.
\begin{definition}
Let $C\subset \P\sp 2$ be a reduced irreducible plane curve.
We say that $C$ is {\it splitting in $X\sb G$} if the proper transform $D\sb C$ of $C$ in $X\sb G$
is not reduced.
If $C$ is splitting in $X\sb G$,
then the divisor $D\sb C$ is written as $2\HPT{C}$,
where $\HPT{C}$ is a reduced irreducible curve on $X\sb G$.
\end{definition}
\begin{definition}
A pencil $\EEE$ 
of cubic curves on $\Pt$  is called a {\it regular pencil splitting in $\XG$}
if the following hold;
\begin{itemize}
\item the base locus of $\EEE$  consists of distinct $9$ points,
\item  every singular member of $\EEE$  is an irreducible nodal curve, and 
\item  every member of $\EEE$ is splitting in $\XG$.
\end{itemize}
\end{definition}
The correctness of our main  algorithm (Algorithm~\ref{algo:main})
is a consequence of the following:
\begin{maintheorem}
Suppose that $G\in \Uts$.
\par
{\rm (1)}
Let $\III\sb{\ZdG}\subset \OOO\sb{\P\sp 2}$ denote  the ideal sheaf of $\ZdG$.
Then the linear system $ | \III\sb{\ZdG}(5) |$
is of dimension $2$,
and the  general member of $ | \III\sb{\ZdG}(5) |$
is reduced, irreducible, and splitting in $X\sb G$.
\par
{\rm (2)}
A line $L\subset\Pt$ is splitting in $X\sb G$
if and only if $\card{L\cap \ZdG}=5$.
\par
{\rm (3)}
A smooth conic  $Q\subset\Pt$ is splitting in $X\sb G$
if and only if $\card{Q\cap \ZdG}=8$.
\par
{\rm (4)}
Let $\EEE$ be a regular pencil of cubic curves of $\Pt$
 splitting in $X\sb G$.
Then the base locus $\Bs(\EEE)$ of $\EEE$ is contained in  $\ZdG$.
\par
{\rm (5)}
The  lattice $\SG$
is generated by the sublattice
$\SG\sp 0$
and the classes $[\HPT{C}]$, 
where $C$ runs through the set of splitting curves of the following type: %{:}
\begin{itemize}
\item the general member of the linear system $ | \III\sb{\ZdG}(5) |$,
\item a  line splitting in $\XG$,
\item a smooth conic splitting in $\XG$, 
\item  a member of a regular pencil  of cubic curves 
 splitting in $X\sb G$.
\end{itemize}
\end{maintheorem}
\begin{example}\label{example:DK}
Consider the polynomial 
\begin{equation}\label{eq:DK}
G\sb{\DK}:= X\sb 0 X\sb 1 X\sb 2 (X\sb 0 ^3 + X\sb 1 ^3 + X\sb 2 ^3),
\end{equation}
which was discovered by 
Dolgachev and Kondo
in~\cite{DK}.
They  showed that every supersingular $K3$ surface
in characteristic $2$ with Artin invariant $1$ 
is isomorphic to $X\sb{G\sb{\DK}}$.
The subscheme $Z(d G\sb{\DK}) \st \Pt$ consists of the $\F\sb 4$-rational points of $\Pt$.
A line $L\st\Pt$ is splitting in $X\sb{G\sb{\DK}}$
if and only if $L$ is $\F\sb 4$-rational.
The numerical N\'eron-Severi lattice of $X\sb{G\sb{\DK}}$ is generated by the classes
of the $(-2)$-curves 
$$
\Gamma\sb P\quad (P\in \Pt (\F\sb 4))
\quand
\HPT{L}\quad (L\in (\Pt)\dual (\F\sb 4)).
$$
(The classes $[\Ample\sb{G\sb{\DK}}]$ and $[\HPT{C}]$,
where $C$ is the general member of $ | \III\sb{Z (dG\sb{\DK})}(5) |$, 
are written as linear combinations of $[\Gamma\sb P]$ and $[\HPT{L}]$.)
\end{example}
\begin{example}\label{example:1598}
Consider the polynomial
\begin{multline*}
G:={X_{{0}}}^{5}X_{{1}}+{X_{{0}}}^{5}X_{{2}}+{X_{{0}}}^{3}{X_{{1}}}^{3}+{
X_{{0}}}^{3}{X_{{1}}}^{2}X_{{2}}+{X_{{0}}}^{3}X_{{1}}{X_{{2}}}^{2}+\\+ {X_
{{0}}}^{3}{X_{{2}}}^{3}+{X_{{0}}}^{2}X_{{1}}{X_{{2}}}^{3}+X_{{0}}{X_{{
2}}}^{5}+{X_{{1}}}^{5}X_{{2}}.
\end{multline*}
We put
\begin{eqnarray*}
P\sb 0 &:=&[{\alpha}^{13}+{\alpha}^{11}+{\alpha}^{10}+{\alpha}^{9}+{\alpha}^{7}+{\alpha}^{4}+{\alpha}^{3}+{\alpha}^{2}, \\&&
\phantom{aaaaaa} {\alpha}^{12}+{\alpha}^{11}+{\alpha}^{9}+{\alpha}^{5}+{\alpha}^{3}+{\alpha}^{2}+\alpha,\;1], \quand\\ 
P\sb {7} & :=&[{\alpha}^{12}+{\alpha}^{11}+{\alpha}^{10}+{\alpha}^{7}+{\alpha}^{6}+{\alpha}^{5}+{\alpha}^{4}+\alpha, \\&&
\phantom{aaaaaa} {\alpha}^{13}+{\alpha}^{11}+{\alpha}^{9}+{\alpha}^{5}+{\alpha}^{4}+{\alpha}^{3}+{\alpha}^{2}+\alpha,\;1],
\end{eqnarray*}
where $\alpha$ is a root of the irreducible polynomial 
$$
{t}^{14}+{t}^{13}+{t}^{12}+{t}^{8}+{t}^{5}+{t}^{4}+{t}^{3}+{t}^{2}+1
 \;\in\; \F\sb 2 [t].
$$
The subscheme $\ZdG$ is reduced of dimension $0$
consisting of the points 
$$
P\sb{\nu}:=\Frob^{\nu} (P\sb 0) \quad (\nu=0, \dots, 6) \quand 
P\sb{7+\nu}:=\Frob^{\nu} (P\sb 7) \quad (\nu=0, \dots, 13),
$$
where $\Frob$ is the Frobenius morphism  $\alpha\mapsto \alpha^2$ over $\F\sb 2$.
(We have $\Frob^{7} (P\sb 0)=P\sb 0$ and $\Frob^{14} (P\sb{7})=P\sb{7}$.)
There exists a line $L$
that passes through the points $P\sb 0$, $P\sb 1$, $P\sb{3}$, $P\sb{7}$, $P\sb{14}$. 
There exists a smooth conic  $Q$
that passes through the points $P\sb 7$, $P\sb 8$, $P\sb 9$, $P\sb {11}$, $P\sb {14}$,
$P\sb {15}$, $P\sb {16}$, $P\sb{18}$.
The lattice $\SG$ is generated by the classes in  $\SzG$ and the classes $[\HPT{C}]$  
associated to the general member  of $|\III\sb{\ZdG} (5)|$,
the splitting  lines
 $\Frob\sp{\nu} (L)$ and  
the splitting  smooth conics 
$\Frob\sp{\nu} (Q)$ for  $\nu=0, \dots, 6$. (We have $\Frob^{7} (L)=L$ and $\Frob^{7} (Q)=Q$.)
The Artin invariant of $\XG$ is $4$.
\end{example}
\begin{example}\label{example:1330}
Consider the polynomial
$$
G:={X_{{0}}}^{5}X_{{2}}+{X_{{0}}}^{3}{X_{{1}}}^{3}+{X_{{0}}}^{3}{X_{{2}}}
^{3}+X_{{0}}X_{{1}}{X_{{2}}}^{4}+{X_{{1}}}^{5}X_{{2}}.
$$
The subscheme $\ZdG$ is reduced of dimension $0$
consisting of the point $[0,0,1]$ and the Frobenius orbit of the point 
\begin{multline*}
[{\alpha}^{19}+{\alpha}^{18}+{\alpha}^{16}+{\alpha}^{15}+{\alpha}^{8}+{\alpha}^{3}+{\alpha}^{2}+\alpha,
\\{\alpha}^{19}+{\alpha}^{17}+{\alpha}^{16}+{\alpha}^{15}+{\alpha}^{14}+{\alpha}^{9}+{\alpha}^{8}+{\alpha}^{7}+{\alpha}^{5}+{\alpha}^{3}+\alpha,
\;1],
\end{multline*}
where $\alpha$ is a root of the irreducible polynomial 
$$
{t}^{20}+{t}^{19}+{t}^{18}+{t}^{15}+{t}^{10}+{t}^{7}+{t}^{6}+{t}^{4}+1
 \;\in\; \F\sb 2 [t].
$$
There are no reduced irreducible  plane curves of degree $\le 3$ that
are splitting in $X\sb G$.
Hence  $\SG$ is generated by the classes in $\SzG$ and the class $[\HPT{C}]$ 
associated to the general member of $|\III\sb{\ZdG} (5)|$.
Therefore  the Artin invariant of $X\sb G$ is $10$.
Note that it is a non-trivial problem to find  explicit examples
of supersingular $K3$ surfaces with big Artin invariant.
See~\cite{Shioda87} and~\cite{Goto95, Goto96}.
\end{example}
\begin{example}\label{example:1509}
Consider the polynomial 
$$
G:={X_{{0}}}^{5}X_{{1}}+{X_{{0}}}^{3}{X_{{1}}}^{2}X_{{2}}+X_{{0}}{X_{{2}}
}^{5}+{X_{{1}}}^{5}X_{{2}}.
$$
We put
\begin{eqnarray*}
P\sb 0
&:=&
[{\alpha}^{13}+{\alpha}^{12}+{\alpha}^{10}+{\alpha}^{9}+{\alpha}^{8}+{\alpha}^{3}+{\alpha}^{2},
\; {\alpha}^{13}+{\alpha}^{8}+{\alpha}^{2},\; 1],\quand\\ 
 P\sb {14} &:=&
[{\alpha}^{13}+{\alpha}^{12}+{\alpha}^{11}+{\alpha}^{10}+{\alpha}^{9}+{\alpha}^{8}+{\alpha}^{7}+{\alpha}^{6}+{\alpha}^{2}, \\&&
\phantom{aaaaaa} {\alpha}^{10}+{\alpha}^{9}+{\alpha}^{7}+{\alpha}^{4},\;1],
\end{eqnarray*}
where $\alpha$ is a root of the irreducible polynomial 
$$
{t}^{14}+{t}^{13}+{t}^{12}+{t}^{8}+{t}^{5}+{t}^{4}+{t}^{3}+{t}^{2}+1
 \;\in\; \F\sb 2 [t].
$$
The subscheme $\ZdG$ is reduced of dimension $0$.
It consists  of the points $P\sb{\nu}:=\Frob^{\nu} (P\sb 0)$ $(\nu=0, \dots, 13)$ and 
$P\sb{14+\nu}:=\Frob^{\nu} (P\sb {14})$ $(\nu=0, \dots, 6)$.
(We have $\Frob^{14} (P\sb 0)=P\sb 0$ and $\Frob^{7} (P\sb{14})=P\sb{14}$.)
We put
$$
A:=\{ P\sb{0}, P\sb{1}, P\sb{3}, P\sb{7}, P\sb{8}, P\sb{10}, P\sb{14}, P\sb{18}, P\sb{19}\}.
$$
We have $\Frob\sp 7 (A)=A$.
For each $\nu=0, \dots, 6$,
there exists a regular pencil $\EEE\sb\nu$ of cubic curves splitting in $\XG$
such that  the  base locus  $\Bs (\EEE\sb\nu)$ is equal to  $\Frob\sp{\nu}(A)$.
The lattice $\SG$ is generated by the classes in $\SzG$ and the classes $[\HPT{C}]$  
associated to the general member of $|\III\sb{\ZdG} (5)|$ and  
the members  of $\EEE\sb \nu$ for $\nu=0, \dots, 6$.
The Artin invariant of $\XG$ is $7$.
\end{example}
The configuration of  irreducible curves of degree $\le 3$  splitting in $\XG$
is encoded by  the $2$-elementary group
$$
\CCG\lift:= \SG/ S\sp 0\sb G,
$$
which we will regard as a  linear code in the  $\F\sb 2$-vector space $(\SzG)\dual /\SzG$ of dimension $22$,
where $(\SzG)\dual$ is the dual lattice of $\SzG$.
Using the  basis
$$
[\Gamma\sb P]/2\quad (P\in \ZdG)\quand [\AG]/2
$$
of $(\SzG)\dual$,
we can identify the $\F\sb 2$-vector space $(\SzG)\dual /\SzG$
with
$$
\Pow (\ZdG) \oplus \F\sb 2,
$$
where $\Pow (\ZdG)$ is the power set of $\ZdG$ equipped with a structure of
the $\F\sb 2$-vector space by
$$
A+B=(A\cup B) \sm (A\cap B)\qquad (A, B\subset \ZdG).
$$
We define the code  $\CCG\st\Pow (\ZdG)$
to be the image of  $\CCG\lift$ by the projection $(\SzG)\dual /\SzG\to \Pow (\ZdG)$.
It turns out that we can recover  from  $\CCG$
the numerical N\'eron-Severi lattice $\SG$, 
and obtain
the configuration of   curves of degree $\le 3$  splitting in $\XG$.
In particular, we have
$$
\hbox{the Artin invariant of $\XG$}
\;=\;11-\dim\sb{\F\sb 2}\CCG.
$$
\begin{theorem}\label{thm:characterization1}
Let $\tZ$ be a finite set with $\card{\tZ}=21$,
and let $\tC\subset \Pow(\tZ)$
be a  code.
There exists a polynomial $G\in \Uts$
such that $\tC$ is mapped to $\CCG\subset \Pow(\ZdG)$ by a certain  bijection $\tZ\isom \ZdG$
if and only if $\tC$ satisfies the following conditions;
\begin{itemize}
\item[(a)] $\dim\sb{\F\sb 2} \tC\le 10$,
\item[(b)] the word $\tZ\in \Pow(\tZ)$ is contained in $\tC$, and
\item[(c)] $\card{A}\in \{ 0,5,8,9,12,13,16,21\}$ for every word $A\in \tC$.
\end{itemize}
\end{theorem}
We say that two codes $\tC$ and $\tC\sprime$ in $\Pow (\tZ)$ are said to be {\it $\St$-equivalent}
if there exists a permutation $\tau$ of $\tZ$ such that $\tau (\tC)=\tC\sprime$ holds.
By computer-aided calculation,
we have classified all the $\St$-equivalence classes of codes
satisfying the conditions (a), (b) and (c) in Theorem~\ref{thm:characterization1}.
The list is given in \S\ref{sec:codelist}.
\begin{theorem}\label{thm:nopsclass}
The number $r(\sigma)$ of the  $\St$-equivalence  classes of codes
with dimension $11-\sigma$
satisfying the conditions {\rm (b)} and {\rm (c)} in Theorem~\ref{thm:characterization1}
is given as follows:
\begin{equation}\label{eq:rsigmatable}
{
\renewcommand{\arraystretch}{1.2}
\begin{array}{|c||c|c|c|c|c|c|c|c|c|c|}
\hline
\sigma & 1 & 2 & 3 & 4 & 5 & 6 & 7 & 8 & 9 & 10 \\
\hline
 r(\sigma) & 1 & 3 & 13 &41 &58 &43 & 21 & 8 & 3 & 1 \\
\hline
\end{array}
} \hskip 3pt \lower 4mm \hbox{.}
\end{equation}
\vskip 5pt
\end{theorem}
From the list,
we obtain the following facts  about the stratification of $\Uts$ 
by the Artin invariant.
For $\sigma =1, \dots, 10$, we put
$$
\UUU\sb{\sigma} := \set{G\in\Uts}{\hbox{the Artin invariant of $X\sb G$ is $\sigma$}} \quand 
\UUU\sb{\le \sigma} := \bigcup\sb{\sigma\sprime\le \sigma} {\,\UUU}\sb{\sigma\sprime}.
$$
Note that each ${\,\UUU}\sb{\le \sigma}$ is Zariski closed in  $\Uts$.
\begin{corollary}\label{cor:nopsirred}
The number of the  irreducible components of ${\,\UUU}\sb{\sigma}$
is at least $r(\sigma)$,
where $r(\sigma)$ is  given in~\eqref{eq:rsigmatable}.
\end{corollary}
\begin{corollary}\label{cor:le9}
The Zariski closed subset ${\,\UUU}\sb{\le 9}$ of $\Uts$
consists of three irreducible hypersurfaces ${\,\UUU} [33]$,
${\,\UUU}[42]$ and ${\,\UUU}[51]$,
where ${\,\UUU}[ab]$ is the locus of all $G\in \Uts$ 
that can be written as $G=G\sb a G\sb b +H^2$,
where  $G\sb a$, $G\sb b$ and  $H$ are homogeneous polynomials of degree $a$, $b$ and $3$, 
respectively.
\end{corollary}
\begin{corollary}\label{cor:equalto1}
If the Artin invariant of $X\sb G$ is $1$,
then, via a linear automorphism of $\Pt$,
the covering morphism  $Y\sb G\to \Pt$ is isomorphic to the Dolgachev-Kondo surface $Y\sb{G\sb{\DK}}\to\Pt$
in Example~\ref{example:DK}.
In particular,
the locus ${\,\UUU}\sb{1}$ is irreducible,
and,
in the moduli space 
$\moduli= \P\sb * (\Uts/\Vts)/ \PGL(3, k)$,
the  locus of supersingular $K3$ surfaces
with Artin invariant $1$ consists of a single point.
\end{corollary}
\par
\medskip
Purely inseparable covers of the projective plane are called {\it Zariski surfaces},
and their properties have been studied by P.~Blass and J.~Lang~\cite{BL}.
In particular,
an algorithm to calculate the Artin invariant
has been established~\cite[Chapter~2, Proposition~6]{BL}.
Our algorithm gives us not only the Artin invariant
but also a geometric description
of generators of the numerical N\'eron-Severi group.
\par
\medskip
This paper is organized as follows.
\par
\medskip
As is suggested above,
the global section $dG$ of $\Omega\sp 1 \sb{\Pt} (6)$
plays an important role in the study of $\XG$.
In \S\ref{sec:global},
we study  global sections of $\Omega\sp 1 \sb{\Pt} (b)$ in general,
where $b$ is an  integer $\ge 4$.
The problem that is considered in this section  is
to characterize the subschemes defined by $s=0$, where $s$ is a global section  of $\Omega\sp 1 \sb{\Pt} (b)$,
among  reduced $0$-dimensional subschemes $Z$ of $\Pt$.
A characterization is given in terms of the linear system $|\III\sb{Z} (b-1)|$.
The results in this section hold in any characteristics.

In \S\ref{sec:geometric},
we assume that the ground field is of characteristic $p>0$,
and define a global section $dG$ of $\Omega\sp 1 \sb{\Pt} (b)$,
where $G$ is a homogeneous polynomial of degree $b$ divisible by $p$.
We then investigate   geometric properties of the purely inseparable cover  $\YG\to\Pt$ defined by $w^p=G$,
and the minimal resolution $\XG$ of $\YG$.
Many results of this section have been already presented in~\cite{BL}.

From \S\ref{sec:globalchar2},
we assume that the ground field is of characteristic $2$.
Let $b$ be an even integer $\ge 4$.
In \S\ref{sec:globalchar2},
we consider the problem to determine
whether a given global section of $\Omega\sp 1 \sb{\Pt} (b)$
is written as $dG$
by some homogeneous polynomial $G$.
In \S\ref{sec:codes},
we associate  to a homogeneous polynomial
$G$ 
a binary linear code $\CCG$ that describes the numerical N\'eron-Severi lattice
of $\XG$.
A notion of {\it geometrically realizable $\Sn$-equivalence classes of codes}
is introduced.
In \S\ref{sec:splitting},
we define  a  word $\wG (C)$  of  $\CCG$
for  each  curve $C$ splitting in $\XG$,
and  study the geometry of splitting curves.

From \S\ref{sec:ssk3},
we put $b=6$,
and study the supersingular $K3$ surfaces $\XG$ in characteristic $2$.
In \S\ref{sec:ssk3},
we review some known facts about $K3$ surfaces.
In \S\ref{sec:codelist},
the relation between the code $\CCG$ and the configuration of curves splitting in $\XG$
is explained.
We present the complete list of 
geometrically realizable $\St$-equivalence classes of codes.
Theorems and Corollaries stated above are proved in this section.
 In \S\ref{sec:algorithm},
we present an algorithm that calculates the code $\CCG$ from a given homogeneous polynomial $G\in\Uts$,
and give  concrete examples.
Some irreducible components of ${\,\UUU}\sb{\sigma}$ are described in detail.
\section{Global sections of $\Omega\sp 1 \sb{\Pt} (b)$ in arbitrary characteristic}\label{sec:global}
In this section, we work over an algebraically closed field $k$
of {\it arbitrary} characteristic.
\par
\medskip
Let $b$ be an integer $\ge 4$.
We consider the locally free sheaf
$$
\Ob:=\Omega\sp 1 \sb{\Pt} \otimes \OPt (b)
$$
of rank $2$ on the projective plane $\Pt$.
From the exact sequence
\begin{equation}\label{eq:exact}
0 \;\to\;  \Ob \;\to\; \OPt (b-1)\sp{\oplus 3} \;\to\; \OPt (b) \;\to\; 0,
\end{equation}
we obtain 
$$
n:= c\sb 2 (\Ob) = b^2 -3b+3.
$$
For a global section $s\in \Hz{\Ob}$,
we denote by $Z(s)$ the subscheme of $\Pt$ defined by $s=0$,
and by $\III\sb{Z(s)}\subset \OPt$ the ideal sheaf of $Z(s)$.
If $Z(s)$ is a reduced $0$-dimensional scheme,
then $Z(s)$ consists of $n$ reduced points.
\par
\medskip
The main result of this section is the following:
\begin{theorem}\label{thm:mainsecglobal}
Let $Z$ be a $0$-dimensional reduced subscheme of $\Pt$
with the ideal sheaf $\III\sb Z\subset \OOO\sb{\Pt}$.
Suppose that $\length \OOO\sb Z =n$.
Then the following two conditions are equivalent: 
\par
\smallskip
{\rm (i)}
There exists a global section $s$ of $\Ob$ such that $Z=Z(s)$.
\par
{\rm (ii)}
There exists a pair $(C\sb 0$, $C\sb 1)$ of members of the linear system
$| \III\sb{Z} (b-1) |$
such that the scheme-theoretic intersection $C\sb 0\cap C\sb 1$ is the union 
of $Z$ and a $0$-dimensional subscheme $\Gamma\st\Pt$ of $\length \OOO\sb{\Gamma} =b-2$
that is contained in a line disjoint from $Z$.
\par
\smallskip
If these conditions are satisfied,
then the global section $s$ with  $Z=Z(s)$ is 
unique up to multiplicative constants.
\end{theorem}
Let $[X\sb 0,X\sb 1, X\sb 2]$ be homogeneous coordinates of $\Pt$.
We put
$$
l\sb{\infty} :=\{ X\sb 2=0\},
\quad
U:=\Pt\sm l \sb{\infty},
$$
and let $(x\sb 0, x\sb 1)$ be the affine coordinates on $U$
given by 
$$
x\sb 0:=X\sb 0/ X\sb 2 \quand x\sb 1:=X\sb 1/X\sb 2.
$$
We also regard $[x\sb 0, x\sb 1]$ as homogeneous coordinates of $ l \sb{\infty}$.
Let $e\sb b$ be the global section of $\OPt (b)$
that corresponds to $X\sb 2^b\in \Hz{\OPt (b)}$.
A section 
\begin{equation}\label{eq:sU}
\sigma\sb 0 (x\sb 0, x\sb 1) dx\sb 0 \otimes e\sb b  + \sigma\sb 1 (x\sb 0, x\sb 1) dx\sb 1\otimes e\sb b
\end{equation}
of $\Ob$ on $U$ extends to a global section of $\Ob$ over $\Pt$
if and only if the following holds;
\begin{equation}\label{eq:extend1}
\parbox{10cm}{
the polynomials $\sigma\sb 0$, $\sigma\sb 1$, and $\sigma\sb 2 := x\sb 0 \sigma\sb 0 + x\sb 1 \sigma\sb 1$
are of degree $\le b-1$.
}
\end{equation}
For $i=0,1$ and $2$,
let $\sigma \sb i \sp{(b-1)} (x\sb 0, x\sb 1)$ be the homogeneous part of degree $b-1$ of $\sigma\sb i$.
Then the condition~\eqref{eq:extend1} is rephrased as follows;
\begin{equation}\label{eq:extend2}
\parbox{10.5cm}{
$\deg \sigma\sb 0 <b$, $\deg\sigma\sb 1<b$,
and there exists a homogeneous polynomial $\gamma (x\sb 0, x\sb 1)$
of degree $b-2$ 
such that $\sigma \sb 0 \sp{(b-1)}=x\sb 1 \gamma$ and   $\sigma \sb 1 \sp{(b-1)}=-x\sb 0 \gamma$.
}
\end{equation}
In particular, we have
$$
\hz{\Ob}=b^2-1.
$$
This equality also  follows from the exact sequence~\eqref{eq:exact}.
\begin{remark}\label{rem:Zsatinfty}
Suppose that a global section $s$ of $\Ob$  is given by~\eqref{eq:sU} on $U$.
The subscheme $Z(s)$ of $\Pt$  is defined on $U$ by $\sigma\sb 0=\sigma\sb 1=0$.
The intersection $Z(s)\cap  l \sb{\infty}$
is set-theoretically equal to the common zeros of the homogeneous polynomials
$\sigma\sb 0\sp{(b-1)}$,
$\sigma\sb 1\sp{(b-1)}$ and 
$\sigma\sb 2\sp{(b-1)}$ on $l \sb{\infty}$.
In particular,
if $s\in \Hz{\Ob}$ is chosen generally,
then $Z(s)$ is reduced of dimension $0$.
\end{remark}
Let $\Theta$
be the sheaf of germs of regular vector fields
on $\Pt$, that is,
$\Theta$ is the dual of $\Omega\sp 1\sb{\Pt}$.
Let $e\sb{-1}$ be the rational section of $\OPt (-1)$ that corresponds to $1/ X\sb 2$.
The vector space $\Hz{\Theta (-1))}$
is of dimension $3$,
and is generated by $\theta\sb 0, \theta\sb 1, \theta\sb 2$,
where
$$
\theta\sb 0 |U =\Der{}{x\sb 0}\otimes e\sb{-1}, \quad
\theta\sb 1| U=\Der{}{x\sb 1}\otimes e\sb{-1}, \quad
\theta\sb 2 |U= \left( x\sb 0 \Der{}{x\sb 0} + x\sb 1 \Der{}{x\sb 1}\right ) \otimes e\sb{-1}. 
$$
Since $c\sb 2 (\Theta (-1))=1$,
every non-zero global section $\theta$ of $\Theta (-1)$ has a single reduced zero,
which we will denote by $\zeta ([\theta])$,
where $[\theta]\in \P\sb * (\Hz{\Theta(-1)})$
is the one-dimensional linear subspace of $\Hz{\Theta(-1)}$  generated by $\theta$.
When $\theta$ is given by
\begin{equation*}\label{eq:theta1}
\theta |U = A \theta\sb 0 + B\theta \sb 1+ C\theta\sb 2\qquad (A, B, C\in k),
\end{equation*}
then $\zeta ({[\theta]})$ is equal to
$[A, B, -C]$
in terms of the homogeneous coordinates $[X\sb 0, X\sb 1, X\sb 2]$.
Thus we obtain  an isomorphism
$$
\mapisom{\zeta}{\P\sb * (\Hz{\Theta(-1)}) }{ \Pt}.
$$
For a hyperplane $V\subset \Hz {\Theta(-1)}$,
we denote by $l\sb V \subset \Pt$ the line corresponding to $V$ by $\zeta$.
For a line $l\subset \Pt$,
we denote by $V\sb l\subset \Hz {\Theta(-1)}$ the hyperplane corresponding to $l$ by $\zeta$.
\begin{remark}\label{rem:tau}
Suppose that a hyperplane $V$ of $\Hz{\Theta(-1)}$ is generated by $\tau\sb 0$ and $\tau\sb 1$.
Then there exist affine coordinates $(y\sb 0, y\sb 1)$
on $U\sb V :=\Pt\sm l\sb V$ and a rational section $e\sprime \sb{-1}$  of $\OPt (-1)$
having the  pole along $l\sb V$ such that
$$
\tau\sb 0 | U\sb V = \Der{}{y\sb 0}\otimes e\sprime \sb{-1},
\quad
\tau\sb 1 | U\sb V = \Der{}{y\sb 1} \otimes e\sprime \sb{-1}.
$$
\end{remark}
A global section $s$ of $\Ob$
defines a linear homomorphism
$$
\map{\varphi\sb s}{\Hz{\Theta(-1)}}{\Hz {\IZs (b-1)}}
$$
via the natural coupling $\Omega\sp 1\sb{\Pt} \otimes \Theta \to \OPt$.
Suppose that $s$ is given by~\eqref{eq:sU}.
For $i=0,1$ and $2$,
we put
$$
\tilde \sigma\sb i (X\sb 0, X\sb 1, X\sb 2):= X\sb 2 ^{b-1} \sigma\sb i(X\sb 0 / X\sb 2, X\sb 1 / X\sb 2).
$$
Then $\varphi\sb s$ is given by 
\begin{equation}\label{eq:varphiexplicit}
\varphi\sb s (\theta\sb i)=\tilde\sigma\sb i\qquad (i=0, 1, 2).
\end{equation}
\begin{proposition}\label{prop:1secglobal}
Let $s$ be a global section of $\Ob$ such that $Z(s)$ is reduced of dimension $0$.
Then the following hold: %{:}
\par
{\rm (1)}
The linear homomorphism $\varphi\sb s$ is an isomorphism.
\par
{\rm (2)}
Let $l\subset \Pt$ be a line such that $l\cap Z(s)=\emptyset$,
and let $P\sb{s, l} \subset | \IZs (b-1) |$ be the pencil corresponding to the hyperplane $V\sb l \subset \Hz{\Theta(-1)}$
via the isomorphism 
$\varphi\sb s$.
Then the base locus of $P\sb{s, l}$ is of the form
$$
Z(s) +\Gamma(s, l),
$$
where $\Gamma (s, l)$ is a $0$-dimensional scheme of $\length \OOO\sb{\Gamma(s, l)} =b-2$.
Moreover  the ideal sheaf $\III\sb{\Gamma(s, l)}\st \OPt$ of $\Gamma(s, l)$
contains the ideal sheaf $\III\sb l$ of the line $l$.
\end{proposition}
\begin{proof}
First we show that $\varphi\sb s$ is injective.
Suppose that there exists  a non-zero global section $\theta$
of $\Theta(-1)$ such that $\varphi\sb s (\theta)=0$.
We have affine coordinates $(y\sb 0, y\sb 1)$
on some affine part $U\sprime$ of $\Pt$ such that
$$
\theta|U\sprime = \Der{}{y\sb 0}\otimes e\sprime\sb{-1},
$$
where $e\sprime\sb{-1}$ is a rational section  of $\OPt (-1)$ that is regular on $U\sprime$.
We express $s$ by 
$$
s|U\sprime = (\sigma\sprime\sb 0 d y\sb 0 + \sigma\sprime\sb 1 d y\sb 1)\otimes e\sprime\sb {b}, 
$$
where $e\sprime\sb{b}:=1/ (e\sprime \sb{-1})\sp{\otimes b}$.
Since $\varphi\sb s (\theta)=0$,
we have $\sigma\sprime \sb 0=0$.
Because $Z(s)$ is of dimension $0$,
$Z(s)\cap U\sprime$ must be empty.
Hence  $\sigma\sprime \sb 1$ is a non-zero constant.
Because $b\ge 4$,
the line $\Pt\sm U\sprime$ at infinity
is contained in $Z(s)$ by Remark~\ref{rem:Zsatinfty},
which contradicts the assumption.
Therefore  $\varphi\sb s$ is injective.
\par
Next we prove (2).
We choose the homogeneous coordinates $[X\sb 0, X\sb 1, X\sb 2]$
in such a way that $l$ is defined by $X\sb 2=0$.
The hyperplane $V\sb l$ of $\Hz{\Theta(-1)}$
is generated by $\theta \sb 0$ and $\theta\sb 1$.
Since their images by $\varphi\sb s$ are
$\tilde\sigma\sb 0$ and $\tilde\sigma\sb 1$, 
the pencil $P\sb{s, l}\st |\III\sb{Z(s)} (b-1)|$
is spanned by the curves $C\sb 0$ and $C\sb 1$ of degree $b-1$
defined by $\tilde\sigma\sb 0=0$ and $\tilde\sigma\sb 1=0$.
Since $Z(s)\cap l=\emptyset$ by the assumption,
we see from Remark~\ref{rem:Zsatinfty}
that the scheme-theoretic intersection $C\sb 0\cap C\sb 1\cap U$
coincides with $Z(s)$,
and at least one of  $C\sb 0$ or $C\sb 1$ does not contain $l$ as an irreducible component.
Hence the base locus of $P\sb{s, l}$ is $Z(s) + \Gamma (s, l)$,
where $\Gamma(s, l)$ is a $0$-dimensional scheme whose support is contained in $l$.
We have
$$
\length\OOO\sb{\Gamma(s, l)}=(b-1)^2-n =b-2.
$$
Note that the support of $\Gamma (s, l)$
is the zeros on $l$ of the homogeneous polynomial $\gamma$ of degree $b-2$ that has appeared in~\eqref{eq:extend2}.
Suppose that $s$ is general.
Then $\gamma$ is a reduced polynomial, and hence $\Gamma (s, l)$ is equal to
the reduced scheme defined by $X\sb 2 =\gamma (X\sb 0, X\sb 1)=0$,
because their supports and lengths coincide.
In particular, the ideal sheaf  $\III\sb{\Gamma (s, l)}$ of $\Gamma (s, l)$
contains the ideal sheaf $\III\sb l$  of $l$.
By the specialization argument,
we see that $\III\sb{\Gamma (s, l)}$ contains $\III\sb l$
for any $s$ such  that $Z(s)$ is reduced,  of dimension $0$ and disjoint from $l$.
\par
It remains to show that $\varphi\sb s$ is surjective.
It is enough to show that 
$$
\hz{\IZs (b-1)}= 3.
$$
We follow the argument of~\cite[pp.\,712-714]{GH}.
Let $\pi : S\to\Pt$ be the blow-up of $\Pt$ at the points of $Z(s)$,
and let $E$ be the union of $(-1)$-curves on $S$ that are contracted by $\pi$.
We have 
$$
E^2=-n, \quad
K\sb S\cong  \pi\sp * \OPt (-3) \otimes  \OOO\sb S (E),
\quand
h\sp 0 (S, K\sb S)=h\sp 1 (S, K\sb S)=0.
$$
Let $L\to S$ be the line bundle corresponding to the invertible sheaf
$$
\pi\sp * \OPt (b-1) \otimes  \OOO\sb S (-E).
$$
There exists a natural isomorphism
\begin{equation}\label{eq:eq2}
H\sp 0 (S, L) \cong \Hz {\IZs (b-1)}.
\end{equation}
From $h^2(S, L)= h\sp 0 (S, K\sb S-L)=0$ and $\chi (\OOO\sb S)=1$,
we obtain from the  Riemann-Roch theorem that
\begin{equation}\label{eq:eq3}
h\sp 0 (S, L) = h\sp 1 (S, L) - (b^2-7b+6)/2.
\end{equation}
Let $\xi\sb 0$ and $\xi\sb 1$ be 
the global sections of the line bundle $L$
corresponding to the homogeneous polynomials $\varphi\sb s (\theta\sb 0)= \tilde \sigma \sb 0$ 
and $ \varphi\sb s (\theta\sb 1)= \tilde\sigma \sb 1 $ in $ \Hz{\IZs (b-1)}$
by the natural isomorphism~\eqref{eq:eq2}.
Since $Z(s)$ is reduced,
the curves $C\sb 0=\{\tilde\sigma\sb 0=0\}$
and $C\sb 1=\{ \tilde\sigma\sb 1=0\}$
are smooth at each point of $Z(s)$, and 
they intersect transversely at each point of $Z(s)$.
Hence the divisors on $S$
defined by $\xi\sb 0=0$ and $\xi\sb 1=0$
have no common points on $E$.
Therefore we can construct 
the Koszul complex
$$
0\;\to\; \OOO\sb S (K\sb S-L) \;\to\; 
\OOO\sb S (K\sb S) \oplus \OOO\sb S (K\sb S) \;\to\; 
\III\sb{\pi\inv (\Gamma(s, l))} (K\sb S +L)\;\to\; 0
$$
from $\xi\sb 0$ and $\xi\sb 1$,
where $\III\sb{\pi\inv (\Gamma(s, l))}\st\OOO\sb S$
is the ideal sheaf of $\pi\inv (\Gamma(s, l))$.
From this complex, 
we obtain
\begin{equation}\label{eq:eq4}
h\sp 1 (S, L)= h\sp 0 (S, \III\sb{\pi\inv (\Gamma(s, l))} (K\sb S +L)) =\hz {\III\sb{\Gamma(s, l)} (b-4)}.
\end{equation}
\par
Suppose that  $b=4$.
Then we have $\hz {\III\sb{\Gamma(s, l)} (b-4)}=0$,
and hence,
from~\eqref{eq:eq2}-\eqref{eq:eq4},
we obtain $\hz{\IZs (b-1)}= 3$.
\par
Suppose that $b\ge 5$.
Assume that the general member $D$ of 
$|\III\sb{\Gamma(s, l)} (b-4)|$
satisfies $l\not\subset D$.
Then the length of the scheme-theoretic intersection of $l$ and $D$ is $b-4$.
Since $\III\sb D \subset \III\sb {\Gamma (s, l)}$ and $\III\sb l \subset \III\sb {\Gamma (s, l)}$,
this contradicts $\length \OOO\sb{\Gamma(s, l)}=b-2$.
Therefore the linear system $|\III\sb{\Gamma(s, l)} (b-4)|$
possesses $l$ as a fixed component.
Since $\III\sb{\Gamma (s, l)} \supset \III\sb l$, we have
\begin{equation}\label{eq:eq5}
\hz {\III\sb{\Gamma(s, l)} (b-4)} = \hz {\OPt (b-5)} =3+(b^2-7b+6)/2.
\end{equation}
Combining~\eqref{eq:eq2}-\eqref{eq:eq5},
we obtain $\hz{\IZs (b-1)}= 3$.
\end{proof}
\begin{remark}\label{rem:Phi}
Let $s\in \Hz{\Ob}$ be as in Proposition~\ref{prop:1secglobal}.
The $2$-dimensional linear system $|\IZs (b-1)|$
defines a morphism
$$
\Phi\sb s \;:\; \Pt\sm Z(s)  \;\to\; \P\sp * (\Hz{\IZs(b-1)}) \cong (\Pt)\dual, 
$$
where the second isomorphism is obtained from the isomorphism $\varphi\sb s$ and the dual of $\zeta$.
Let $l\in (\Pt)\dual$ be a general line of $\Pt$.
The inverse  image  of $l$ by $\Phi\sb s$ coincides with $\Gamma (s, l)$.
Therefore $\Phi\sb s$ is generically finite of degree $b-2$.
\end{remark}
\begin{remark}\label{rem:psi}
Let $s$, $l$, $V\sb l$ and $P\sb{s, l}$  be as in Proposition~\ref{prop:1secglobal}.
We have isomorphisms $P\sb{s, l} \cong \P\sb * (V\sb l)$ by $\varphi\sb s$,
and $\P\sb * (V\sb l) \cong l$ by $\zeta$.
By composition, we obtain an isomorphism 
$$
\mapisom{\psi\sb {s, l}}{P\sb{s, l}}{l}.
$$
The restriction of the pencil $P\sb{s, l}$ to $l$
consists of the fixed part $\Gamma (s, l)$
and one moving point.
The isomorphism $\psi\sb {s, l}$ maps  $C\in P\sb{s, l}$
to the moving point of the divisor $C\cap l$ of $l$.
Indeed, let us fix  affine coordinates $(x\sb 0, x\sb 1)$
on $U=\Pt\sm l$ as in the proof of Proposition~\ref{prop:1secglobal}
so that $V\sb l$ is generated by $\theta\sb 0$ and $\theta\sb 1$.
The isomorphism $\P\sb * (V\sb l) \cong l$ is written explicitly as 
$$
\zeta ([\theta\sb 0  + t \theta \sb 1])=[1,t, 0]\in l.
$$
On the other hand,
the projective plane curve of degree $b-1$ defined by the homogeneous  polynomial 
$$
\varphi\sb s (\theta\sb 0  + t \theta \sb 1) = \tilde \sigma \sb 0 + t \tilde \sigma\sb 1
$$
passes through the point $[1,t, 0]$ by~\eqref{eq:extend2}.
\end{remark}
\begin{corollary}\label{cor:1secglobal}
Let $s$ be a global section of $\Ob$ such that
$Z(s)$ is reduced of dimension $0$.
Then the linear system $|\IZs (b-1)|$
is of dimension $2$,
and its base locus coincides with $Z(s)$.
The general member of $|\IZs (b-1)|$
is reduced and irreducible.
\end{corollary}
\begin{proof}
The last statement follows from the assumption  that  $Z(s)$ is reduced and  from Bertini's theorem 
applied to the morphism  $\Phi\sb s$ in Remark~\ref{rem:Phi}.
\end{proof}
\begin{proof}[Proof of Theorem~\ref{thm:mainsecglobal}]
The implication from (i) to (ii) has been already proved in Proposition~\ref{prop:1secglobal}.
Suppose that $|\III\sb Z  (b-1)|$ has the property (ii).
We will construct a global section $s$ of $\Ob$ such that $Z=Z(s)$.
Let $l$ be the  line of $\Pt$  containing the subscheme $\Gamma$.
We choose homogeneous coordinates $[X\sb 0, X\sb 1, X\sb 2]$
such that $l$ is defined by $X\sb 2=0$.
Let $\tilde\sigma\sb 0(X\sb 0, X\sb 1, X\sb 2)=0$
and $\tilde\sigma\sb 1(X\sb 0, X\sb 1, X\sb 2)=0$
be the defining equations of $C\sb 0$ and $C\sb 1$,
respectively.
We put
\begin{eqnarray*}
&\sigma\sb 0 (x\sb 0, x\sb 1):= \tilde\sigma\sb 0 (x\sb 0, x\sb 1, 1),
\quad&
\sigma\sb 1 (x\sb 0, x\sb 1):= \tilde\sigma\sb 1 (x\sb 0, x\sb 1, 1),\\
&\sigma\sp{(b-1)} \sb 0 (x\sb 0, x\sb 1):= \tilde\sigma\sb 0 (x\sb 0, x\sb 1, 0),
\quad&
\sigma\sp{(b-1)} \sb 1 (x\sb 0, x\sb 1):= \tilde\sigma\sb 1 (x\sb 0, x\sb 1, 0).
\end{eqnarray*}
Let $\gamma (x\sb 0, x\sb 1)$
be the homogeneous polynomial of degree $b-2$ such that $\gamma=0$ defines the subscheme $\Gamma$ on the line $l$.
Since $C\sb 0\cap C\sb 1$ is scheme-theoretically equal to $Z+\Gamma$,
and $l$ is disjoint from $Z$,
the scheme-theoretic intersection $C\sb 0\cap C\sb 1 \cap l$ coincides with $\Gamma$.
Hence there 
exist linearly independent homogeneous linear forms $\lambda\sb 0(x\sb 0, x\sb 1)$
and $\lambda\sb 1(x\sb 0, x\sb 1)$ such that
$$
\sigma\sp{(b-1)} \sb 0=\lambda\sb 0 \gamma,\quad
\sigma\sp{(b-1)} \sb 1=\lambda\sb 1 \gamma.
$$
By linear change of coordinates $(x\sb 0, x\sb 1)$,
we can assume that $\lambda\sb 0=x\sb 1$
and $\lambda\sb 1=-x\sb 0$.
Then the section 
$(\sigma\sb 0 d x\sb 0 + \sigma\sb 1 d x\sb 1)\otimes e\sb b $ of $\Ob$
on $\Pt\sm l$ extends to a global section $s$ of $\Ob$.
We have 
$Z(s)\cap (\Pt\sm l)=C\sb 0\cap C\sb 1 \cap (\Pt\sm l)=Z$.
Because $l\not\st Z(s)$,
the subscheme $Z(s)$ is of dimension $0$.
Since the length $n=c\sb 2 (\Ob)$ of $\OOO\sb{Z(s)}$ is equal to that of $\OOO\sb Z$,
we have $Z=Z(s)$.
\par
Next we prove the uniqueness (up to multiplicative constants)  of $s$ satisfying $Z=Z(s)$.
Let $s\sprime$ be another global section of $\Ob$
such that $Z(s\sprime)=Z$.
The morphism
$$
\widetilde{\Phi} \sb{Z} \;:\; \Pt\sm Z  \;\to\; \P\sp * (\Hz{\III\sb Z (b-1)})
$$
defined by the  linear system $|\III\sb{Z} (b-1)|$
does not depend on the choice of $s$.
Let $\widetilde{P}\in  \P\sp * (\Hz{\III\sb Z (b-1)})$
be a general point.
By Remark~\ref{rem:Phi},
there exist lines $l$ and $l\sprime$ of $\Pt$ such that
$\widetilde{\Phi}\sb Z\inv (\widetilde {P})$
is equal to $\Gamma(s, l)=\Gamma(s\sprime, l\sprime)$.
On the other hand,
since the length $b-2$ of $\OOO\sb{\Gamma(s, l)}$  is $\ge 2$ by the assumption $b\ge 4$,
the subscheme $\Gamma(s, l)$ determines the line $l$ containing $\Gamma(s, l)$ uniquely.
Hence we have $l=l\sprime$,
which implies that $\Phi\sb s=\Phi\sb{s\sprime}$.
Therefore the linear isomorphisms $\varphi\sb s$ and $\varphi\sb{s\sprime}$
are equal up to a multiplicative constant,
and hence so are $s$ and $s\sprime$ by~\eqref{eq:varphiexplicit}.
\end{proof}
\begin{remark}
If there exists a pair $(C\sb 0, C\sb 1)$ of members of $|\III\sb Z (b-1)|$
satisfying the condition in Theorem~\ref{thm:mainsecglobal}\,(ii),
then the  {\it general } pair  of  members of $|\III\sb Z (b-1)|$
also satisfies it.
\end{remark}
\section{Geometric properties of purely inseparable covers of $\Pt$}\label{sec:geometric}
In this section,
we assume that the ground field $k$ is of positive characteristic $p$.
We fix a  multiple $b$ of $p$ greater than or equal to $4$.
\subsection{Definition of $\Upb$.}
Let $\MMM$ and $\LLL$ be line bundles on $\Pt$
corresponding to the invertible sheaves $\OPt (b/p)$ and $\OPt (b)$,
respectively.
We have a canonical  isomorphism
\begin{equation}\label{eq:isomMpL}
\MMM\sp{\otimes p} \;\;\isom \;\; \LLL.
\end{equation}
Using  this isomorphism,
we have  local trivializations of the line bundle $\LLL$ such that 
the transition functions 
are $p$-th powers,
and hence the usual differentiation of  functions 
defines a linear homomorphism
\begin{equation*}\label{eq:d}
\Hz {\LLL} \;\; \to\;\; \Hz{\Omega\sp 1\sb{\Pt} \otimes \LLL}=\Hz{\Ob},
\end{equation*}
which we  denote by $G\mapsto dG$.
We put 
$$
\Vpb:=\set{H\sp p\in \Hz{\LLL}}{ H\in \Hz{\MMM}}.
$$
Note that $\Vpb$ is a {\it linear} subspace of $\Hz{\LLL}$, 
because we are in characteristic $p$.
In fact, 
the kernel of the linear homomorphism $G\mapsto dG$ is  equal to $\Vpb$.
\par
Let $[X\sb 0, X\sb 1, X\sb 2]$ be homogeneous coordinates of $\Pt$,
and let $U$ be the affine part $\{X\sb 2\ne 0\}$ of $\Pt$,
on which affine coordinates $x\sb0 :=X\sb 0/X\sb 2$ and  $x\sb 1:=X\sb 1/X\sb 2$
are defined.
Suppose that a global section $G$ of $\LLL$ is given by a homogeneous polynomial 
$G(X\sb 0, X\sb 1, X\sb 2)$ of degree $b$.
Then $dG$ is given by
$$
dG|U =\left( \Der{g}{x\sb 0} dx\sb 0 +\Der{g}{x\sb 1} dx\sb 1\right) \otimes e\sb b,
$$
where $g(x\sb0, x\sb 1):=G(x\sb 0, x\sb 1, 1)$,
and $e\sb b$ is the section of $\LLL$ corresponding to $X\sb 2^b$
\begin{definition}
Let $G$ and $G\sprime$ be  global sections of $\LLL$.
We write $G\sim G\sprime$
if there exist a non-zero constant $c$ and a global section $H$ of $\MMM$
such that $G=c\,G\sprime +H\sp p$.
\end{definition}
\begin{remark}\label{rem:equivalgo}
For a homogeneous polynomial 
$G:=\sum\sb{i+j+k=b} a\sb{ijk} X\sb 0 ^i X\sb 1^j X\sb 2^k$
of degree $b$, we put
$$
\bar G:=\sum\sb{(i, j, k) \not\equiv (0, 0, 0)\,\bmod p} a\sb{ijk} X\sb 0 ^i X\sb 1^j X\sb 2^k.
$$
Let $G$ and $G\sprime$ be two global sections of $\LLL$.
Then $G\sim G\sprime$ holds
if and only if there exists a non-zero constant $c$ such that $\bar G=c\,\bar G\sprime$.
\end{remark}
Let $G$ be a global section of $\LLL$.
Using the isomorphism~\eqref{eq:isomMpL},
we can define a subscheme $\YG$ of the total space of the line bundle $\MMM$
by the equation
$$
w\sp p =G,
$$
where $w$ is a fiber coordinate of $\MMM$.
We denote by 
$$
\map{\piG}{\YG}{\Pt}
$$
the canonical projection,
which is a purely inseparable finite morphism of degree $p$.
It is easy to see that,
set-theoretically, we have
$$
\pi\sb G\inv (\ZdG)=\Sing (\YG).
$$
\begin{remark}
If $G\sim G\sprime$, then we have $Z(dG)=Z(dG\sprime)$,
and 
the schemes $\YG$ and $Y\sb{G\sprime}$ are isomorphic over $\Pt$.
\end{remark}
\begin{proposition}\label{prop:Upb}
For a global section $G$ of $\LLL$,
the following conditions are equivalent to each other: %{:}
\begin{itemize}
\item[{\rm(i)}]
The subscheme $\ZdG$ of $\Pt$ is reduced of dimension $0$.
\item[{\rm(ii)}]
For any $G\sprime$ with $G\sprime \sim G$,
the curve defined by $G\sprime =0$ 
has only ordinary nodes as its singularities.
\item[{\rm(iii)}]
The surface $\YG$ has only rational double points
of type $A\sb{p-1}$ as its singularities.
\end{itemize}
If $G$ is chosen generally from $\Hz{\LLL}$,
then $G$ satisfies these conditions.
\end{proposition}
\begin{proof}
Let $P$ be an arbitrary point of $\Pt$,
and  $Q$  the unique point of $\YG$
such that $\piG (Q)=P$.
We fix affine coordinates $(x\sb 0, x\sb 1)$ 
with the origin $P$ on an affine part $U\st \Pt$.
Let $G$ be expressed on $U$ by an inhomogeneous polynomial of $x\sb 0$ and $x\sb 1$;
$$
G|U= c\sb{00} + c\sb{10}x\sb 0 + c\sb{01} x\sb 1 + c\sb{20} x\sb 0^2 + c\sb{11} x\sb 0 x\sb 1 + c\sb{02}  x\sb 1^2
+ \textrm{(terms of higher degrees)}.
$$
Let $G\sprime$ be another global section of $\LLL$ that is expressed on $U$ by 
$$
G\sprime|U= c\sprime\sb {00} + c\sprime\sb{10}x\sb 0 + c\sprime\sb{01} x\sb 1 + c\sprime\sb{20} x\sb 0^2 +
c\sprime\sb{11} x\sb 0 x\sb 1 + c\sprime\sb{02}  x\sb 1^2 + \textrm{(terms of higher degrees)}.
$$
If $G\sim G\sprime$, there exists a non-zero constant $c$ such that
$$
c\sprime\sb{10}=c\,c\sb{10},
\quad
c\sprime\sb{01}= c\,c\sb{01},
\quand
c\sprime\sb{11}= c\,c\sb{11}.
$$
If $p>2$,  we also have
$$
c\sprime \sb{20}= c\,c\sb{20},
\quand 
c\sprime\sb{02}= c\,c\sb{02}.
$$
Since $\ZdG$ is defined by 
$$
\Der{(G|U)   }{x\sb 0}=\Der{(G|U)  }{x\sb 1}  =0
$$
locally around $P$, 
we have the following equivalences,
from which the equivalence of the conditions (i), (ii) and (iii) follows:
\begin{eqnarray*}
&& P\notin \ZdG\\
&\Longleftrightarrow& c\sb{10}\ne 0 \;\;\textrm{or}\;\; c\sb{01}\ne 0 \\
&\Longleftrightarrow&
\hbox{if $G\sprime \sim G$ and $G\sprime (P)=0$, then the curve defined by $G\sprime=0$ is 
smooth at $P$}\\ 
&\Longleftrightarrow&
\hbox{$\YG$ is smooth at $Q$};\\
%\end{eqnarray*}
%
%\vskip -30pt
%
%\begin{eqnarray*}
&& \hbox{$P$ is a reduced isolated point of $\ZdG$}\phantom{\vrule height 14pt}\\
&\Longleftrightarrow& c\sb{10}= c\sb{01}= 0 \quand 4 c\sb{20} c\sb{02}-c\sb{11}^2 \ne 0 \\
&\Longleftrightarrow&
\hbox{if $G\sprime \sim G$ and $G\sprime (P)=0$, then the curve defined by $G\sprime=0$ is reduced at $P$}\\
& & \hbox{and has an ordinary node  at $P$}\\
&\Longleftrightarrow&
\hbox{$\YG$ has a rational double point of type $A\sb{p-1}$ at  $Q$}.
\end{eqnarray*}
As was shown above, the locus
$$
N\sb P:=
\Biggl\{ 
\,
G\in \Hz{\LLL}
\;
\Bigg|
\;
\vcenter{
\hbox{
\vbox{
\hbox{$P\in \ZdG$, and}
\hbox{$P$ is {\it not} a reduced isolated point of $\ZdG$}
}
}
}
\Biggr\}
$$
is of codimension $3$ in $\Hz{\LLL}$ for any $P\in\Pt$.
Therefore,
if $G\in \Hz{\LLL}$ is general, 
$G$ is not contained in $N\sb P$ for any $P\in \Pt$,
and hence $\ZdG$ is reduced of dimension $0$.
\end{proof}
\begin{definition}\label{def:Upb}
We denote by $\Upb$ the Zariski open dense subset of $\Hz{\LLL}$
consisting of all $G$ satisfying the conditions in Proposition~\ref{prop:Upb}.
Note that, if $G\in \Upb$ and $G\sprime \sim G$,
then
$G\sprime \in \Upb$.
For $G\in \Upb$,
we put 
$$
k\sptimes G +\Vpb := \set{c\,G+H^p}{c\in k\sptimes, \; H\in\Hz{\MMM}}=\set{G\sprime \in \Upb}{G\sim G\sprime}.
$$
\end{definition}
\begin{remark}\label{rem:linearsystemIZdG}
By the linear homomorphism 
$$
\map{\varphi\sb{dG}}{\Hz{\Theta(-1)}}{\Hz{\III\sb{\ZdG} (b-1)}}
$$
that is an isomorphism for $G\in \Upb$ in virtue of  Proposition~\ref{prop:1secglobal}, 
we see that the $2$-dimensional linear system $|\III\sb{\ZdG} (b-1)|$ is spanned by the three curves
defined by $ \der{G}{X\sb 0}=0$, $ \der{G}{X\sb 1}=0$ and $ \der{G}{X\sb 2}=0$.
\end{remark}
\subsection{Geometric properties of  $\XG$ for $G\in \Upb$.}
From now on,
we fix a polynomial $G\in \Upb$.
Then $\Sing(\YG)$ 
consists of $n=b^2 -3b+3$ rational double points of type $A\sb{p-1}$.
Let
$$
\map{\phiG}{\XG}{\Pt}
$$
denote the composite of the minimal resolution $\XG\to\YG$ of $\YG$ and the purely inseparable finite morphism $\piG$.
We denote by $\AG \subset \XG$ the pull-back of a general line of $\Pt$ via $\phiG$.
\begin{proposition}\label{prop:KX}
The canonical divisor $K\sb G$ of the nonsingular surface $\XG$ is linearly equivalent
to $(b-b/p-3) \AG $.
\end{proposition}
\begin{proof}
Let $(x\sb 0, x\sb 1)$ be affine coordinates on an affine part $U$ of $\Pt$ that contains
$\ZdG$,
and let $g(x\sb 0, x\sb 1)$ be the  inhomogeneous  polynomial that corresponds to  $G$ on $U$.
On the surface $\YG$, we have
$$
0= d(w^p) =\Der{g}{x\sb 0} dx\sb 0 + \Der{g}{x\sb 1} dx\sb 1.
$$
The rational $2$-form
$$
\frac{dw\wedge d x\sb 0}{\partial g/ \partial x\sb 1} =-\frac{dw\wedge d x\sb 1}{\partial g/ \partial x\sb 0} 
$$
is therefore regular and nowhere vanishing on the Zariski open dense subset
$$
\piG\inv (U\sm \ZdG) =\piG\inv (U)\sm \Sing(\YG)
$$
of $\YG$.
By direct calculation,
we can show that this rational $2$-form has a zero of order $b-b/p-3$
along the pull-back $\piG\inv ( l \sb{\infty})$
of the line $ l \sb{\infty}:=\Pt\sm U$ at infinity.
Since  $\Sing (\YG)$ consists of only rational double points,
the canonical divisor of $\XG$ is $(b-b/p-3)$ times $\phiG\inv ( l \sb{\infty})$.
\end{proof}
\begin{definition}
We denote by $\SG$ the numerical N\'eron-Severi lattice of $\XG$,
and by $\SzG$ the sublattice of  $\SG$
that is generated by the class $[\AG]$,
and the classes $[\Gamma\sb i]$ $(i=1, \dots, n(p-1))$
of smooth rational curves $\Gamma\sb i$ on $\XG$ that are contracted to  the singular points of $\YG$.
\end{definition}
\begin{proposition}\label{prop:pelem}
The quotient group $\SG/\SzG$ is a finite  elementary $p$-group.
\end{proposition}
\begin{proof}
Let $C$ be a reduced irreducible curve on $\XG$.
If $\phiG(C)$ is a point, then $C$ is one of the  curves $\Gamma\sb i$,
and hence $[C]\in \SzG$.
Suppose that $\phiG (C)$ is of dimension $1$.
Let $D\st \Pt$ denote the curve $\phiG (C)$
with the reduced structure,
and let $\widetilde D \subset \XG$
be the proper transform of $D$ by $\phiG$.
Obviously 
we have 
$[\widetilde D]\in \SzG$.
If the morphism $\phiG|\sb{C} : C\to D$ is birational,
then $\widetilde D= pC$ holds,
because $\phiG$ is purely inseparable of degree $p$ 
over the generic point of $D$.
Hence we have $p[C]\in \SzG$.
If $\phiG |\sb{C} : C\to D$ is of degree $>1$,
then it must be of degree $p$ and $C=\widetilde D$ holds,
and hence  $[C]$ is contained in $\SzG$.
\end{proof}
Since $[\AG]$ and $[\Gamma\sb i]$ $(i=1, \dots, n(p-1))$ are linearly independent in $\SzG\otimes \Q$,
we obtain the following:
\begin{corollary}\label{cor:rank}
The rank of $\SG$ is equal to $n(p-1)+1$.
\end{corollary}
\begin{definition}
A non-singular projective surface  $X$ is called {\it supersingular}
(in the sense of Shioda)
if the rank of the numerical N\'eron-Severi lattice of  $X$ is equal to the second Betti number $b\sb 2 (X)$.
\end{definition}
\begin{definition}
A reduced irreducible  surface  $X$ is called {\it unirational}
if there exists a dominant rational map from $\Pt$ to $X$.
\end{definition}
\begin{proposition}\label{prop:unirational}
The surface $\XG$ is unirational and supersingular.
\end{proposition}
\begin{proof}
Let $k (x\sb 0, x\sb 1)$
be  the  rational function field of $\Pt$.
Since $\phiG:\XG\to\Pt$ is purely inseparable of degree $p$,
the function field of $\XG$ is contained in the purely transcendental extension
$k (x\sb 0\sp{1/p}, x\sb 1\sp{1/p})$ of $k$.
Therefore $\XG$ is unirational.
The supersingularity of $\XG$ then follows from~\cite[Corollary 2]{Shioda74}.
\end{proof}
\begin{remark}
Note that the second Betti number $n(p-1)+1$ of $\XG$ is equal to that of 
a $p$-th cyclic cover of a {\it complex} projective plane branched along a nonsingular plane curve of 
degree $b$.
\end{remark}
\section{Global sections of $\Ob$ in characteristic $2$}\label{sec:globalchar2}
From  this section, we assume that $p=2$.
Let $b$ be an even integer $\ge 4$.
\par
\medskip
Let $s$ be a global section of $\Ob$
such that $Z(s)$ is reduced of dimension $0$.
Recall from Remark~\ref{rem:Phi} that 
the $2$-dimensional linear system $|\IZs (b-1)|$
defines a morphism
$$
\Phi\sb s \;:\; \Pt\sm Z(s) \;\to\; \P\sp * (\Hz{\IZs(b-1)}) \cong (\Pt)\dual.
$$
\begin{proposition}
There exists a polynomial $G\in \Utb$  such that $s=dG$ holds
if and only if the  morphism 
$\Phi\sb s$ is inseparable.
\end{proposition}
\begin{proof}
Recall that, for a general $l\in (\Pt)\dual$, 
the inverse image of $l$ by $\Phi\sb s$ is the divisor $\Gamma (s, l)$
of $l$ with degree $b-2$ defined in Proposition~\ref{prop:1secglobal}.
Therefore the following three conditions on $s$ are equivalent to each other:
\begin{itemize}
\item[(i)]
The  morphism $\Phi\sb s$ is inseparable.
\item[(ii)]
For a general line $l\subset\Pt$,
there exists a divisor $\Delta (s, l)$ of $l$ with degree $b/2-1$ such that $\Gamma(s, l)=2\Delta(s, l)$ holds.
\item[(iii)]
Let $(x\sb 0, x\sb 1)$ be general affine coordinates of $\Pt$,
and let $s$ be given on the affine part by $(\sigma\sb 0 dx\sb 0 + \sigma\sb 1 dx\sb 1)\otimes e\sb b$.
Then there exists a homogeneous polynomial $\delta (x\sb 0, x\sb 1)$ 
of degree $b/2-1$ 
such that $\sigma\sb 0\sp{(b-1)}=x\sb 1 \delta\sp 2 $
and $\sigma\sb 1\sp{(b-1)}=x\sb 0 \delta\sp 2 $ hold.
\end{itemize}
Suppose that there exists $G\in \Utb$ such that $s=dG$.
Let $(x\sb 0, x\sb 1)$ be  general affine coordinates
on an affine part $U$.
Then $G|U$ is written as follows;
$$
\gamma\sb{00} (x\sb 0, x\sb 1)\sp2 
+ x\sb 0 \gamma\sb{10} (x\sb 0, x\sb 1)\sp2 
+ x\sb 1 \gamma\sb{01} (x\sb 0, x\sb 1)\sp2 
+ x\sb 0 x\sb 1 \gamma\sb{11} (x\sb 0, x\sb 1)\sp2,
$$
where $\gamma\sb{00}$ is an inhomogeneous polynomial of degree $\le b/2$,
and $\gamma\sb{10}$, $\gamma\sb{01}$ and $\gamma\sb{11}$ are 
inhomogeneous polynomials of degree $\le b/2-1$.
Then $s=dG$ is written on $U$ as
$$
((\gamma\sb{10}\sp 2 + x\sb 1 \gamma\sb{11}\sp 2) dx\sb 0 
+
(\gamma\sb{01}\sp 2 + x\sb 0 \gamma\sb{11}\sp 2) dx\sb 1)\otimes e\sb b.
$$
Therefore  the homogeneous part of $\gamma\sb{11}$
of degree $b/2-1$  yields the polynomial $\delta$ required in the condition (iii). 
\par
Conversely,
suppose that the condition (ii) holds.
Again we choose affine coordinates 
$(x\sb 0, x\sb 1)$
of $\Pt$ defined on an affine part $U\st\Pt$ containing $Z(s)$,
and let $s$ be given by $(\sigma\sb 0 dx\sb 0 + \sigma\sb 1 dx\sb 1)\otimes e\sb b$ on $U$.
Let $l$ be a line defined by
$$
x\sb 0 + A x\sb 1 + B=0 \quad(A, B\in k).
$$
Then the hyperplane $V\sb l \subset \Hz{\Theta(-1)}$ corresponding to $l$ 
via $\zeta$ is 
generated by $\theta\sb{\infty}$ and $\theta\sb{0}$,
where
$$
\theta\sb{\infty}|U = \left ( A\Der{}{x\sb 0} + \Der{}{x\sb 1}\right) \otimes e\sb{-1},
\quand
\theta\sb{0}|U = \left ( B\Der{}{x\sb 0} + x\sb 0\Der{}{x\sb 0} + x\sb 1 \Der{}{x\sb 1} \right) \otimes e\sb{-1}.
$$
For $u\in k$, we put
$$
\theta\sb u := u\theta\sb{\infty} + \theta\sb{0}\;\;\in\;\; V\sb l.
$$
The zero point $\zeta([\theta\sb u])$ of $\theta\sb u$ is $(Au+B, u)\in l$.
The member $C\sb u$ of the pencil $P\sb{s, l} \subset |\IZs(b-1)|$
corresponding to $\theta\sb u$ via the isomorphism $\varphi\sb s$ is defined by
$$
\varphi\sb s (\theta\sb u)=(Au+B)\sigma\sb0 + u \sigma\sb 1 + (x\sb 0 \sigma\sb 0 + x\sb1 \sigma\sb 1)=0.
$$
We put $t:=x\sb 1|\sb{l}$,
which is an affine parameter of the line $l$.
The divisor of $l$ cut out by $C\sb u$ is defined by the polynomial
$$
\varphi\sb s (\theta\sb u) (At+B, t)=(u+t) (A \sigma\sb 0 (At+B, t) +\sigma\sb 1 (At+B, t))
$$
of $t$.
Therefore the pencil $\{l\cap C\sb u\}$ of divisors on $l$ cut out by $P\sb{s, l}$
has a unique  moving point  $(Au+B, u)$ corresponding to the factor $u+t$,
and the fixed part
$$
\Gamma(s, l)=\{ A \sigma\sb 0 (At+B, t) +\sigma\sb 1 (At+B, t)=0\}.
$$
By the assumption,
we see that
\begin{eqnarray*}
&& \frac{d}{dt} (A \sigma\sb 0 (At+B, t) +\sigma\sb 1 (At+B, t))\\
&=& A^2 \Der{\sigma\sb 0}{x\sb 0} (At+B, t) + 
A\Bigl( \Der{\sigma \sb 0}{x\sb 1} +\Der{\sigma \sb 1}{x\sb 0}\Bigr)(At+B, t)+\Der{\sigma
\sb 1}{x\sb 1} (At+B, t)
\end{eqnarray*}
is zero for generic (and hence all) $A$, $B$ and $t$.
Therefore we have
$$
\Der{\sigma\sb 0}{x\sb 0}\equiv 0, \quad 
\Der{\sigma \sb 0}{x\sb 1} \equiv \Der{\sigma \sb 1}{x\sb 0}, \quad
\Der{\sigma\sb 1}{x\sb 1}\equiv 0.
$$
This implies that there exist polynomials $\alpha$, $\beta$ and $\gamma$ such that
$$
\sigma\sb 0= \alpha\sp 2 + x\sb 1 \gamma\sp 2,
\quad
\sigma\sb 1 =\beta\sp 2 +x\sb 0 \gamma\sp 2.
$$
We put
$$
g:= x\sb 0 \alpha \sp 2 + x\sb 1 \beta\sp 2 + x\sb 0 x\sb 1 \gamma\sp 2,
$$
and let $G$ be the homogeneous polynomial of degree $b$ obtained from $g$
by homogenization.
Since $\der{g}{x\sb 0}=\sigma\sb 0$ and $\der{g}{x\sb 1}=\sigma\sb 1$,
we have  $dG=s$.
\end{proof}
\section{Codes arising from purely inseparable double covers of $\Pt$}\label{sec:codes}
We assume  that $p=2$
and that $b$ is an even integer $\ge 4$.
\par\medskip
{\bf Remark on  notation.}
From this section,
we use  typewriter fonts $\tZ$, $\tS\sp0\sb{\tZ}$, $\tC$, $\tS\sb{\tZ} (\tC)$, $\tth$, $\te\sb{\tP}$   
and  $\tP\in \tZ$
in the situation where we are dealing 
with abstract codes and lattices 
in order to distinguish them from the corresponding objects
$\ZdG$, $\SG\sp 0$, $\CCG$, $\SG$, $[\AG]$, $[\Gamma\sb P]$  and $P\in \ZdG$
of geometric origin.
\subsection{The discriminant group of a lattice.}\label{subsec: discrim}
In this subsection,
we review the theory of  discriminant groups of  lattices due to Nikulin~\cite{N1}.
\par
\medskip
A {\it lattice} is a free $\Z$-module of finite rank with a non-degenerate symmetric bilinear form
$$
\Lambda\times\Lambda\;\;\to\;\;\Z
$$
denoted by $(u, v)\mapsto uv$.
A lattice $\Lambda$ is said to be {\it even} if $u\sp 2 \in 2\Z$
holds for every $u\in \Lambda$.
For a lattice $\Lambda$,
let $\Lambda\dual$ denote the $\Z$-module $\Hom (\Lambda,\Z)$.
We have a natural injective homomorphism
$\Lambda\inj\Lambda\dual$,
whose cokernel
$$
\DG (\Lambda):=\Lambda\dual/\Lambda
$$
is called the {\it discriminant group} of $\Lambda$.
The order of $\DG(\Lambda)$ is equal, up to sign,  to the discriminant $\disc \Lambda$ of $\Lambda$.
We denote by
$$
\map{\pr\sb{\Lambda}}{\Lambda\dual}{\DG(\Lambda)}
$$
the natural projection.
We have a $\Q$-valued symmetric bilinear form on $\Lambda\dual$
that extends the symmetric bilinear form on $\Lambda$.
Hence a symmetric bilinear form
$$
\map{b\sb{\Lambda}}{\DG(\Lambda)\times \DG(\Lambda)}{\Q/\Z}
$$
is defined.
When $\Lambda$ is an even lattice,
the quadratic form $u\mapsto u^2$ on $\Lambda\dual $ 
induces  a quadratic form
$$
\map{q\sb{\Lambda}}{\DG(\Lambda)}{\Q/2\Z}
$$
on $\DG (\Lambda)$ that relates to $b\sb\Lambda$ by
$$
b\sb\Lambda (u, v)=\frac{1}{2}\bigl(\,q\sb\Lambda (u+v)-q\sb\Lambda(u)-q\sb\Lambda(v)\,\bigr).
$$
\begin{definition}
For a subgroup $H$ of $\DG (\Lambda)$,
we put
$$
H\sperp :=\set{u\in \DG(\Lambda)}{b\sb\Lambda (u, v)=0\quad\textrm{for all}\quad v\in H}.
$$
A subgroup $H$ of $\DG (\Lambda)$
is called {\it $b$-isotropic} if $H$ is contained in $H\sperp$.
When $\Lambda$ is even,
we say that $H$ is {\it $q$-isotropic}
if $q\sb\Lambda (u)=0$ holds for every $u\in H$.
\end{definition}
An {\it overlattice} of $\Lambda$ is a submodule $\Lambda\sprime$ of $\Lambda\dual$
such that $\Lambda\sprime$ contains $\Lambda$ and that the $\Q$-valued symmetric bilinear form
of $\Lambda\dual$ takes values in $\Z$  on $\Lambda\sprime$.
Let $\Lambda\spprime$ be a lattice,
and suppose that there exists an injective isometry $\Lambda\inj \Lambda\spprime$
such that $\Lambda\spprime/\Lambda$ is finite.
Then we have a canonical injection $\Lambda\spprime\inj \Lambda\dual$,
and $\Lambda\spprime$ can be regarded as an overlattice of $\Lambda$.
When $\Lambda\sprime$ is an overlattice of $\Lambda$,
we have a sequence 
$$
\Lambda \;\subset \; \Lambda\sprime \; \subset\; (\Lambda\sprime)\dual\;\subset\; \Lambda\dual
$$
of submodules of $\Lambda\dual$ such that $[\Lambda\sprime : \Lambda]=[\Lambda\dual : (\Lambda\sprime)\dual]$.
\begin{proposition}[Nikulin~\cite{N1}]\label{prop:nikulin}
Let $\Lambda$ be a lattice.
\par
{\rm (1)}
The correspondence
$$
\Lambda\sprime \mapsto H\sb{\Lambda\sprime} :=\pr\sb \Lambda (\Lambda\sprime),
\quad
H\mapsto \Lambda\sprime\sb H:= \pr\sb \Lambda\inv (H)
$$
gives rise to a bijection between the set of overlattices of $\Lambda$ and the set of
$b$-isotropic subgroups of $\DG(\Lambda)$.
We have $\Lambda\sprime\sb H/\Lambda=H$
and $(\Lambda\sprime\sb H)\dual/\Lambda=H\sperp$.
In particular, the discriminant group 
$\DG (\Lambda\sprime\sb H)$ is isomorphic to $H\sperp /H$.
\par
{\rm (2)}
Suppose that $\Lambda$ is even.
Then the above correspondence yields
a bijection between the set of even overlattices of $\Lambda$ and
the set of $q$-isotropic subgroups of $\DG (\Lambda)$.
\end{proposition}
\subsection{Certain hyperbolic $2$-elementary lattices and associated codes.}\label{subsec:hyperbolic}
\begin{definition}
A lattice $\Lambda$ is called {\it hyperbolic}
if the signature of the real quadratic form on $\Lambda\otimes\R$ is $(1, \rank \Lambda-1)$.
\end{definition}
\begin{definition}
A lattice $\Lambda$ is called {\it $2$-elementary}
if the finite abelian group $\DG(\Lambda)$ is $2$-elementary,
that is,
if $\DG (\Lambda)$ is an $\F\sb 2$-vector space of dimension $\log\sb 2 |\disc \Lambda|$.
\par
A $2$-elementary lattice $\Lambda$ is called {\it of type I}\,
if $u\sp 2 \in \Z$ holds for every $u\in \Lambda\dual$, 
that is, if $b\sb\Lambda (x, x)=0$ holds for every $x\in \DG(\Lambda)$.
\end{definition}
\par
\medskip
Let $\tZ$ be a finite set.
(See Remark on  notation.)
We identify the $\F\sb 2$-vector space $\F\sb 2\sp{\tZ}$ of functions from $\tZ$ to $\F\sb 2$
with the power set $\Pow (\tZ)$
of $\tZ$ by
$$
v\in \F\sb 2 \sp\tZ \;\;\mapsto\;\;v\inv (1)\subset\tZ.
$$
A structure of the $\F\sb 2$-vector space on $\Pow (\tZ)$
is therefore defined by
$$
A+B=(A\cup B)\setminus (A\cap B)\qquad (A, B\st \tZ).
$$
An element of $\Pow (\tZ)$ is called a {\it word}.
For a word $A\st \tZ$,
the cardinality $\card{A}$ is called the {\it weight} of $A$.
\par
We consider an even hyperbolic $2$-elementary lattice
$$
\tS\sp 0 \sb{\tZ} := \bigoplus\sb{\tP\in\tZ} \Z \te\sb{\tP} \oplus \Z \tth
$$
with the symmetric bilinear form  given by 
$$
\te\sb{\tP}  \te\sb{\tQ} =\begin{cases}
-2  & \textrm{if $\tP=\tQ$}\\
0  & \textrm{if $\tP\ne\tQ$}
\end{cases},
\quad
\te\sb{\tP} \tth=0,
\quad
\tth\sp 2 =2.
$$
Then we have
$$
(\tS\sp 0 \sb{\tZ})\dual = \bigoplus\sb{\tP\in\tZ} \Z \,(\te\sb{\tP}/2) \oplus \Z \,(\tth/2)
\quad
\subset \quad
\tS\sp0 \sb{\tZ}\otimes \Q.
$$
The discriminant group $\DG (\tS\sp 0 \sb{\tZ})$ is therefore naturally identified with
$$
\F\sb 2 \sp{\tZ} \oplus \F\sb 2 = \Pow(\tZ) \oplus \F\sb 2
$$
in such a way that
a vector
$$
\sum (a\sb{\tP}/2) \te\sb{\tP}+(b/2)\tth
\qquad(a\sb{\tP}, b\in \Z)
$$
of $(\tS\sp 0 \sb{\tZ})\dual$ corresponds to
$$
 (A, b\,\bmod 2) \;\in\; \Pow (\tZ) \oplus \F\sb 2 ,
\qquad\hbox{where $A=\shortset{\tP\in\tZ}{a\sb{\tP}\equiv 1\,\bmod 2}$.}
$$
Hence we can consider 
subgroups of $\DG(\tS\sp 0\sb{\tZ})$ 
as binary linear codes in $\Pow(\tZ)\oplus \F\sb 2$.
Under this identification,
the symmetric bilinear form  $b\sb{\tS\sp 0 \sb{\tZ}}$ on $\DG(\tS\sp0\sb{\tZ})$ 
is given by 
$$
((A, \alpha), (A\sprime, \alpha\sprime))\mapsto 
\begin{cases}
(-\card{A\cap A\sprime} + 1)/2 \,\bmod \Z  & \textrm{if $\alpha=\alpha\sprime=1$},\\
-\card{A\cap A\sprime}/2 \,\bmod \Z  & \textrm{otherwise},
\end{cases}
$$
and the quadratic form $q\sb{\tS\sp 0 \sb{\tZ}}$ on $\DG(\tS\sp0\sb{\tZ})$ is given by
$$
(A, \alpha)\mapsto 
\begin{cases}
(-\card{A}+1 )/2 \,\bmod 2\Z  & \textrm{if $\alpha=1$, }\\
-\card{A}/2 \,\bmod 2\Z  & \textrm{if $\alpha=0$}.
\end{cases}
$$
Therefore, from Proposition~\ref{prop:nikulin},
we obtain the following:
\begin{corollary}\label{cor:codes}
Let $\ttC$ be a code in $\Pow(\tZ)\oplus \F\sb2$, 
which is considered as a subgroup of  $\DG(\tS\sp0\sb{\tZ})$
by the identification above.
\par
{\rm (1)}
If the submodule $\pr\sb{\tS\sp0\sb{\tZ}}\inv (\ttC)$ of $(\tS\sp 0\sb{\tZ})\dual$
corresponding to $\ttC$ 
is an overlattice of $\tS\sp0\sb{\tZ}$, then
the following holds;
\begin{equation}\label{eq:overlattice}
\card{A} \,\bmod2 \equiv\alpha\qquad\hbox{for every $(A, \alpha)\in \ttC$}.
\end{equation}
\par
{\rm (2)}
The submodule $\pr\sb{\tS\sp0\sb{\tZ}}\inv (\ttC)$
is an even overlattice of $\tS\sp0\sb{\tZ}$
if and only if every $(A, \alpha)\in \ttC$ satisfies
$$
\card{A}\equiv\begin{cases}
0\, \bmod 4  & \textrm{if $\alpha=0$},\\
1\,  \bmod 4  & \textrm{if $\alpha=1$}.
\end{cases}
$$
\end{corollary}
We denote by
$$
\map{\rho\sb{\tZ}}{\Pow(\tZ)\oplus\F\sb2}{\Pow(\tZ)}
$$
the projection onto the first factor.
\begin{definition}
Let $\tC$ be an arbitrary code in $\Pow(\tZ)$.
We put
$$
\tC\lift:=\set{ (A, \alpha)\in \Pow(\tZ)\oplus \F\sb 2}{A\in\tC\quand \card{A}\,\bmod2 =\alpha},
$$
and call it the {\it lift }of $\tC$.
It is obvious  that $\tC\lift$ is a linear subspace of $\Pow(\tZ)\oplus \F\sb 2$,
that $\dim \tC\lift$ is equal to $\dim \tC$, 
and that $\tC\lift$ is the unique code satisfying~\eqref{eq:overlattice}
and 
$\rho\sb{\tZ} (\tC\lift)=\tC$.
\par
We denote by $\tS\sb{\tZ} (\tC)$ the submodule 
$\pr\sb{\tS\sp0\sb{\tZ}}\inv (\tC\lift)$ of $(\tS\sp0\sb{\tZ})\dual$.
\end{definition}
If the submodule $\tS\sb{\tZ} (\tC)$ of $(\tS\sp0\sb{\tZ})\dual$ is an overlattice of $\tS\sp0\sb{\tZ}$,
then we have 
\begin{equation}\label{eq:disc}
|\disc (\tS\sb{\tZ} (\tC) ) | = 2\sp{n+1}/ \card{\tC}^2.
\end{equation}
Moreover  the lattice  $\tS\sb{\tZ} (\tC)$ is   hyperbolic and $2$-elementary,
because so is $\tS\sp0\sb{\tZ}$.
From Proposition~\ref{prop:nikulin},
we obtain the following:
\begin{proposition}\label{prop:evencode}
The  submodule $\tS\sb{\tZ} (\tC)$ of $(\tS\sp0\sb{\tZ})\dual$
 is an even overlattice of $\tS\sp0\sb{\tZ}$ if and only if $\card{A}\equiv 0$ or $1 \,\bmod 4$
holds for every $A\in\tC$.
\end{proposition}
\begin{proposition}\label{prop:typeIcode}
Suppose that $n=\card{\tZ}$ is odd,
and that 
 $\tS\sb{\tZ} (\tC)$  is an overlattice of $\tS\sp0\sb{\tZ}$.
If $\tC$ contains the word $\tZ$, then the $2$-elementary lattice $\tS\sb{\tZ} (\tC)$ is of type I.
\end{proposition}
\begin{proof}
Suppose that  $\tC$ contains $\tZ$.
Then $\tC\lift$ contains $(\tZ, 1)$ because $\card{\tZ}$ is odd.
If $(A, \alpha)\in (\tC\lift)\sperp$, then
$$
b\sb{\tS\sp0\sb{\tZ}} ( (\tZ, 1), (A, \alpha))=(-\card{A}+\alpha)/2 =0 \quad\textrm{in}\quad \Q/\Z,
$$
and hence
$$
b\sb{\tS\sp0\sb{\tZ}} ( (A, \alpha), (A, \alpha))=(-\card{A}+\alpha)/2 =0.
$$
If $u\in (\tS\sb{\tZ}(\tC))\dual$, then  $u\,\bmod \tS\sp0\sb{\tZ}\in \DG(\tS\sp0\sb{\tZ})$ is contained in $(\tC\lift)\sperp$,
and therefore $u^2\in \Z$ holds.
Hence $\tS\sb{\tZ} (\tC)$ is of type I.
\end{proof}
\subsection{The  lattice $\SG$ and the associated code.}
We fix a polynomial $G\in \Utb$.
Then $\Sing (\YG)$ consists of  $n=b^2-3b+3$ ordinary nodes
that are mapped bijectively to the points of $\ZdG$.
\begin{definition}
For a point $P\in \ZdG$,
we denote by $\Gamma\sb P$ the $(-2)$-curve on $\XG$ that is contracted to $P$
by $\phiG:\XG\to\Pt$.
\end{definition}
In the numerical N\'eron-Severi lattice $\SG$ of $\XG$, we have
$$
[\Gamma\sb P]   [\Gamma\sb Q]  =\begin{cases}
-2  & \textrm{if $P=Q$}\\
0  & \textrm{if $P\ne Q$}
\end{cases},
\quad
[\Gamma\sb P] [\AG]=0,
\quad
[\AG]\sp 2 =2.
$$
By sending $\te\sb P$ to $[\Gamma\sb P]$ and $\tth$ to $[\AG]$,
we obtain  an isomorphism
\begin{equation}\label{eq:isomS}
\tS\sp0\sb{\ZdG}\;\cong \;\SG\sp 0.
\end{equation}
Hence $\DG (\SG\sp 0)$ is identified with $\Pow(\ZdG) \oplus \F\sb 2$.
Since $\SG/\SzG$ is finite by Proposition~\ref{prop:pelem},
we can regard $\SG$ as an overlattice of $\SG\sp 0$.
\begin{definition}
We put
\begin{eqnarray*}
\tCCG &:= & \SG/\SG\sp0\;\subset\; \DG(\SG\sp0)=\Pow(\ZdG)\oplus\F\sb2, \quand\\
\CCG &:= &\rho\sb{\ZdG} (\tCCG)\; \subset\; \Pow(\ZdG).
\end{eqnarray*}
Note that $\tCCG$ is the lift $\CCG\lift$ of $\CCG$,
and that the overlattice  $\tS\sb{\ZdG} (\CCG)=\pr\sb{\tS\sp 0\sb{\ZdG}}\inv (\tCCG)$ of $\tS\sp 0 \sb{\ZdG}$
corresponding to $\CCG$ 
is identified with 
the overlattice $\SG$ of $\SzG$ 
by the isomorphism~\eqref{eq:isomS}.
\end{definition}
\begin{proposition}\label{prop:RR}
{\rm (1)}
Suppose that $b/2$ is odd.
Then $\card{A}\equiv 0\;\textrm{or}\; 1\,\bmod 4$ for every $A\in \CCG$.
{\rm(2)}
Suppose that $b/2$ is even.
Then $\card{A}\equiv 0\;\textrm{or}\; 3\,\bmod 4$ for every $A\in \CCG$.
\end{proposition}
\begin{proof}
Let $\KG$ be the canonical divisor of $\XG$.
By Proposition~\ref{prop:KX},
we have $[\KG]=(b/2-3)\,[\AG]$ in $\SG$.
Let $A$ be a word in $\CCG$.
Suppose that $\card{A}$ is even.
Then we have $(A, 0)\in \tCCG$,  and hence the vector 
$$
v:=\frac{1}{2}\sum\sb{P\in A} [\Gamma\sb P]
$$
of $(\SzG)\dual$ is contained  in $\SG$. Since $v\sp 2=-\card{A}/2$ and $v\cdot[\KG]=0$, 
we have 
$$
(v^2-v \cdot[\KG])/2 =-{\card{A}}/{4},
$$
which is an integer by the the Riemann-Roch theorem.
Therefore $\card{A}\equiv 0\,\bmod 4$ holds.
Suppose that $\card{A}$ is odd.
Then we have $(A, 1)\in \tCCG$,  and hence 
$$
w:=\frac{1}{2}\Bigl( \sum\sb{P\in A} [\Gamma\sb P] \;+\;[\AG]\Bigr)
$$
is contained in $\SG$.
From
$$
(w^2-w \cdot [\KG])/2=(7-\card{A}-b)/4\;\in\;\Z,
$$
we have $\card{A}+b\equiv 3\,\bmod 4$.
\end{proof}
\subsection{Geometric realizability of an abstract code.}\label{subsec:geometric}
Let $\tZ$ be a finite set with
$$
\card{\tZ}=n=b^2-3b+3.
$$
The symmetric group $\Sn$ acts on $\tZ$ and $\Pow(\tZ)$.
\begin{definition}
Two codes $\tC$ and $\tC\sprime$ in $\Pow(\tZ)$ are said to be {\it $\Sn$-equivalent}
if there exists $\tau\in \Sn$ such that $\tau(\tC)=\tC\sprime$.
We denote by $[\tC]$
the $\Sn$-equivalence class of codes containing the code $\tC\st\Pow(\tZ)$.
\end{definition}
\begin{definition}
Let $\tC$ be a code in $\Pow(\tZ)$, and 
let $[\tC]$ be the  $\Sn$-equivalence class of codes 
containing  $\tC$.
We say that $[\tC]$ is {\it geometrically realizable}
if there exist $G\in \Utb$ and a bijection $\tZ\isom \ZdG$
that maps $\tC\subset \Pow(\tZ)$  to $\CCG\subset\Pow (\ZdG)$.
\end{definition}
\begin{definition}
Let $[\tC]$ and $[\tC\sprime]$ be two $\Sn$-equivalence classes of codes in $\Pow (\tZ)$.
We write $[\tC]<[\tC\sprime]$
if there exist representatives $\tC\in [\tC]$ and  $\tC\sprime\in [\tC\sprime]$
 such that $\tC \subsetneqq \tC\sprime$.
\end{definition}
Let $[\tC]$ be  a geometrically realizable class of codes.
We put
\begin{eqnarray*}
\UUU\sb{2, b, [\tC]} & := &\set{G\in \Utb}{\tC \cong \CCG \quad\textrm{by some bijection $\tZ \cong \ZdG $}}, \quand\\
\UUU\sb{2, b, \ge [\tC]} & := & \bigsqcup \sb{[\tC\sprime]\ge [\tC]} \UUU\sb{2, b, [\tC\sprime]}.
\end{eqnarray*}
\begin{theorem}\label{thm:closed}
For every $[\tC]$,  the locus ${\,\UUU}\sb{2, b, \ge [\tC]}$ is Zariski closed in $\Utb$.
\end{theorem}
\begin{proof}
Let
$\tUtb\to \Utb$
be the \'etale covering of degree $n!$ over $\Utb$ such that 
each point of $\tUtb$  over $G\in \Utb$ is a pair 
$(G, \tau\sb G)$, where $\tau\sb G$ is a bijection from $\tZ$ to $\ZdG$.
For a word $A\in \Pow (\tZ)$, we put
$$
\tU\sb A:=\set{(G, \tau\sb G)\in \tUtb}{ \tau\sb G (A)\in \CCG}.
$$
Since the specialization homomorphism
of  numerical N\'eron-Severi lattices is injective
for a smooth family of projective varieties,
the locus  $\tU\sb A$ is Zariski closed  in $\tUtb$.
For a geometrically realizable  class $[\tC]$,
the closed subset
$$
\bigcup\sb{\tC\in [\tC]} \;\Bigl(\; \bigcap\sb{A\in \tC} \tU\sb A  \;\Bigr)
$$
of $\tUtb$ is invariant  under the $\Sn$-action on $\tUtb$ over $\Utb$,
and is the  pull-back of the locus ${\,\UUU}\sb{2, b, \ge [\tC]}$.
Therefore ${\,\UUU}\sb{2, b, \ge [\tC]}$ is closed in $\Utb$.
\end{proof}
\begin{corollary}
For every geometrically realizable class
$[\tC]$ of codes,  the locus  ${\,\UUU}\sb{2, b, [\tC]}$ is locally Zariski  closed
in
$\Utb$.
\end{corollary}
\begin{remark}
The \'etale covering   
$\tUtb\to \Utb$ that has appeared in the proof of Theorem~\ref{thm:closed}  is constructed as follows.
Let  $\ZZZ \to \Utb$ be the universal family
$$
\set{ (P, G)\in \Pt \times \Utb}{P\in \ZdG} \;\;\to\;\;\Utb
$$
of $\ZdG$, which is an \'etale covering of degree $n$.
We fix a base point $G\sb 0\in \Utb$,
and let
$$
\map{\mu}{\pione (\Utb, G\sb 0)}{\Aut (Z(dG\sb 0))\cong\Sn}
$$
be 
the monodromy action 
of the algebraic fundamental group of $\Utb$
on the set $Z(dG\sb 0)$.
Let $\widetilde \ZZZ \to\Utb$ be the Galois closure of  $\ZZZ\to\Utb$,
which is an \'etale cover of degree equal to the cardinality of $\Im \mu$.
Then $\tUtb$ is a disjoint union of 
$[\Sn: \Im\mu]$ copies of $\widetilde \ZZZ$.
\end{remark}
\subsection{An algorithm for listing up codes.}\label{subsec:listingupalgorithm}
In this subsection,
we describe an algorithm that will be used in~\S\ref{sec:algorithm},
when we make the complete list of  geometrically realizable  classes of codes
for supersingular $K3$ surfaces in characteristic $2$.
\par
Let $\tZ$ be a finite set with $\card{\tZ}=n$.
Suppose that we are given a subset $\WT$ of $\{0, 1, 2, \dots, n\}$.
\begin{problem}\label{problem:listup}
Make the complete list $L\sb k$ $(k=1, \dots, n)$
of the $\Sn$-equivalence classes $[\tC]$ of codes $\tC\subset\Pow (\tZ)$
with the following properties;
\begin{itemize}
\item[(a)] $\dim \tC=k$,
\item[(b)] $\tZ\in \tC$, and 
\item[(c)] $\card{A}\in \WT$ for every $A\in \tC$.
\end{itemize}
%
%
%Then find all pairs of $[\tC]\in L\sb k$ and $[\tC\sprime]\in L\sb{k\sprime}$
%such that $[\tC]<[\tC\sprime]$.
\end{problem}
First we fix some notation and terminologies.
For a code $\tC\subset \Pow(\tZ)$, we put
$$
\wtenum (\tC):=\sum\sb{A\in \tC} x\sp{\card{A}},
$$
where $x$ is a formal variable.
Let $\sA=(A\sb 0, \dots, A\sb{k-1})$
be a sequence of words $A\sb i \in \Pow (\tZ)$.
We denote by $\gen{\sA}\subset \Pow(\tZ)$
the code generated by $A\sb 0, \dots, A\sb{k-1}$.
A sequence $\sA$ of length $k$ is called {\it linearly independent}
if $\dim\gen{\sA}=k$.
We put
$$
\wt (\sA):= (\card{A\sb 0}, \dots, \card{A\sb{k-1}}).
$$
For another word $A\in \Pow (\tZ)$, we write
$$
(\sA, A):=(A\sb 0, \dots, A\sb{k-1}, A).
$$
For $\tau\in \Sn$, we put
$$
\tau (\sA):=(\tau (A\sb 0), \dots, \tau (A\sb{k-1})).
$$
We define a sequence $\tomg(\sA)$ of length $2^k$ by the following:
\begin{itemize}
\item
If $\sA=(A\sb 0)$, then $\tomg (\sA):=(\tZ, A\sb 0)$.
\item
Suppose that $k>1$.
We put 
$\sA\sprime:=(A\sb 0, \dots, A\sb{k-2})$, and
let the sequence 
$\tomg(\sA\sprime )$ be $(B\sb 1, \dots, B\sb{2\sp{k-1}})$.
Then we define 
$$
\tomg (\sA):=(B\sb 1, \dots, B\sb{2\sp{k-1}}, B\sb 1\cap A\sb{k-1}, \dots, B\sb{2\sp{k-1}}\cap A\sb{k-1}).
$$
\end{itemize}
We then define a sequence $\omg (\sA)$ of non-negative integers  by
$$
\omg (\sA) :=\wt (\tomg (\sA)).
$$
Suppose that we are given $\omg (\sA)$.
Then, for any subsets $I$ and $J$ of $\{0, 1, \dots, k-1\}$, the cardinality
$$
\card{ \bigcap \sb{i\in I} A\sb i \cap \bigcap \sb{j\in J} (\tZ\sm A\sb j) }
$$
can be obtained from $\omg (\sA)$.
Therefore, for  two sequences $\sA$ and $\sA\sprime$,
there exists $\tau\in\Sn$ such that $\tau (\sA)=\sA\sprime$
if and only if $\omg(\sA)=\omg(\sA\sprime)$
holds.
In particular, we have the following:
\begin{proposition}\label{prop:isin}
Let $\sA$ be a  sequence of words,
and let $[\tC\sprime]$ be an $\Sn$-equivalence class of codes
containing $\tC\sprime$.
Then $[\gen{\sA}] \leqq [\tC\sprime]$ holds 
if and only if there exists a sequence $\sA\sprime$ of words of $\tC\sprime$
such that $\omg (\sA)=\omg (\sA\sprime)$.
\end{proposition}
The following subroutine determines 
whether two codes $\gen{\sA}$ and $\gen{\sA\sprime}$
given by   sequences $\sA$ and $\sA\sprime$
 are $\Sn$-equivalent or not.
\begin{subroutine}\label{subroutine:equiv}
First we calculate 
$\dim \gen{\sA}$ and $\dim \gen{\sA\sprime}$.
If they differ, then 
$\gen{\sA}$ and $\gen{\sA\sprime}$ are not $\Sn$-equivalent.
Otherwise, we calculate 
the weight enumerators $\wtenum (\gen {\sA})$ and $\wtenum (\gen{\sA\sprime})$.
If they differ, then 
$\gen{\sA}$ and $\gen{\sA\sprime}$ are not $\Sn$-equivalent.
Otherwise,
we calculate $\omg(\sA)$,
and search for a sequence $\sA\spprime$ of words of  $\gen{\sA\sprime}$
such that $\omg(\sA)=\omg(\sA\spprime)$.
Note that, if  $\sA\spprime$ satisfies $\omg(\sA)=\omg(\sA\spprime)$,
then $\dim \gen{\sA\spprime}=\dim \gen{\sA}=\dim \gen{\sA\sprime}$ holds and
hence $\gen{\sA\spprime}$ coincides with $\gen{\sA\sprime}$.
The codes $\gen{\sA}$ and $\gen{\sA\sprime}$ are $\Sn$-equivalent 
if and only if such a sequence $\sA\spprime$ is found.
\end{subroutine}
We label the elements of $\tZ$ as $\{ \tP\sb 0, \dots, \tP\sb {n-1}\}$,
and represent a word $A$ of $\Pow (\tZ)$ by a bit vector
$$
v(A):=[\alpha\sb 0, \dots, \alpha\sb {n-1}],
$$
where $\alpha\sb i=0$ (resp. $\alpha\sb i=1$)
if $\tP\sb i\notin A$ (resp. $\tP\sb i\in A$).
For a  column bit vector $\vb={}\sp{T}[\beta\sb 0, \dots, \beta\sb{k-1}]$,
we put
$$
\mu (\vb):= 2\sp{k-1} \beta\sb 0 + 2\sp{k-2} \beta\sb 1 + \cdots + 2\beta\sb{k-2} +\beta\sb{k-1}\;\;\in\;\; \Z\sb{\ge 0}.
$$
A sequence $\sA=(A\sb 0, \dots, A\sb{k-1})$
is called {\it $\Sn$-increasing}
if the column vectors
of the $k\times n$ matrix
$$
\begin{bmatrix}
\,v (A\sb 0)\, \\
\, \vdots\, \\
\, v(A\sb{k-1}) \,\\
\end{bmatrix}
=
[\, \vb\sb{0}, \dots, \vb\sb{n-1}\,]
$$
yield an increasing sequence $\mu (\vb\sb 0)\le \dots\le \mu(\vb\sb {n-1})$.
The following proposition is obvious from the definition:
\begin{proposition}\label{prop:seq}
{\rm (1)}
If $\sA=(A\sb 0, \dots, A\sb{k-1})$ is $\Sn$-increasing,
then the subsequence 
$(A\sb 0, \dots, A\sb{m-1})$ of $\sA$ is also $\Sn$-increasing
for any $m\le k$,
\par
{\rm (2)}
For any sequence  $\sA=(A\sb 0, \dots, A\sb{k-1})$,
there exists $\tau\in \Sn$
such that $\tau(\sA)$ is $\Sn$-increasing.
\par
{\rm (3)}
Suppose that $\sA=(A\sb 0, \dots, A\sb{k-1})$ is $\Sn$-increasing,
and let $A\in \Pow (\tZ)$ be an arbitrary word.
Then there exists $\tau\in\Sn$ such that $\tau(\sA)$ coincides with $\sA$ and that $(\sA, \tau(A))$
is $\Sn$-increasing.
\end{proposition}
\begin{example}
The sequence given by the first three  row vectors of the  matrix $M$ below is $\SSSS\sb 7$-increasing,
while the sequence of length $4$ 
given by all the  row vectors of $M$ is not $\SSSS\sb 7$-increasing.
By applying  transpositions $\tP\sb 3\leftrightarrow \tP\sb 4$ and $\tP\sb 5\leftrightarrow \tP\sb 6$
to $M$,
we obtain the  matrix $M\sprime$,
which yields the $\SSSS\sb 7$-increasing sequence of length $4$.
$$
M:=\begin{bmatrix}
0 & 0 & 0 & 1 & 1 & 1 & 1 \\
0 & 1 & 1 & 0 & 0 & 1 & 1 \\
0 & 0 & 1 & 0 & 0 & 1 & 1 \\
0 & 1 & 1 & 1 & 0 & 1 & 0
\end{bmatrix},
\qquad
M\sprime:=\begin{bmatrix}
0 & 0 & 0 & 1 & 1 & 1 & 1 \\
0 & 1 & 1 & 0 & 0 & 1 & 1 \\
0 & 0 & 1 & 0 & 0 & 1 & 1 \\
0 & 1 & 1 & 0 & 1 & 0 & 1
\end{bmatrix}.
$$
\end{example}
Let $[\tC]$ be an $\Sn$-equivalence class satisfying the conditions (a), (b) and (c) in Problem~\ref{problem:listup}.
Then there exists a sequence $\sA=(A\sb 0, \dots, A\sb{k-1})$
of length $k$  with the following properties;
\begin{itemize}
\item $\sA$ is linearly independent, and $\gen{\sA}\in [\tC]$,
\item $\sA$ is $\Sn$-increasing,
\item $A\sb 0=\tZ$, and $\card{A\sb i}\le n/2$ for $i=1, \dots, k-1$.
\end{itemize}
Indeed, we have a linearly independent sequence $\sA\sprime=(A\sb 0\sprime, \dots, A\sb{k-1}\sprime)$ that is 
a basis of a code $\tC\in [\tC]$ with $A\sb 0\sprime=\tZ$.
If there is a word $A\sprime\sb i\;(i>0)$
with $\card{A\sprime\sb i}>n/2$, then we replace $A\sprime\sb i$ by $\tZ+A\sprime\sb i$
so that we can assume $\card{A\sprime\sb i}\le n/2$ for $i=1, \dots, k-1$.
By applying a suitable permutation $\tau\in \Sn$,
the sequence 
$\sA:=\tau (\sA\sprime)$
becomes $\Sn$-increasing,
which is a basis of the  code $\tau (\tC)$ in the class $[\tC]$.
\begin{definition}\label{definition:stdbasis}
A sequence $\sA$ with these properties is called a {\it standard basis}
of the $\Sn$-equivalence class $[\tC]$.
\end{definition}
The complete list $L\sb k$ that we want to make will be given as a set
$$
\LL\sb k=\{ \sA\spar{1}, \dots, \sA\spar{N}\}
$$
of standard bases of length $k$.
\begin{proposition}\label{prop:algo}
Suppose that the  complete list $L\sb k$ $(k\ge 1)$ has been given
as a set $\LL\sb k$ of standard bases of length $k$.
Then Algorithm~\ref{algo:listup} below produces  a set $\LL\sb{k+1}$
of standard bases of length $k+1$ that gives the complete list $L\sb{k+1}$.
\end{proposition}
\begin{algorithm}\label{algo:listup}
{\it Step 1.}
For each basis $\sA\spar{i}\in \LL\sb k$,
we make the  list $\AAA\spar{i}$ of words  $A\in \Pow (\tZ)$
with the following properties;
\begin{itemize}
\item[(i)] $\card{A}\le n/2$,
\item[(ii)] $(\sA\spar{i}, A)$ is $\Sn$-increasing, and 
\item[(iii)] for any $B\in \gen{\sA\spar{i}}$, $\card{B+A}\ne 0$ and  $\card{B+A}\in\WT$.
\end{itemize}
In other words,
$\AAA\spar{i}$ is the list of all $A\in \Pow (\tZ)$ such that $(\sA\spar{i}, A)$
is a standard basis of an $\Sn$-equivalence class of $(k+1)$-dimensional codes
satisfying the conditions (b) and (c)  in Problem~\ref{problem:listup}.
\par
\noindent
{\it Step 2.}
Set $\LL\sb{k+1}$ to be an empty set.
\par
\noindent
{\it Step 3.}
For each pair of $\sA\spar{i}\in \LL\sb k$
and $A\in \AAA\spar{i}$,
we check whether there exists
$\sA\sprime\in \LL\sb{k+1}$ such that
$\gen{\sA\sprime}$ and $\gen{(\sA\spar{i}, A)}$
are $\Sn$-equivalent
by using Subroutine~\ref{subroutine:equiv}.
If there are no such $\sA\sprime$,
then we put $(\sA\spar{i}, A)$ in $\LL\sb{k+1}$.
\end{algorithm}
\begin{proof}%[Proof of Proposition~\ref{prop:algo}]
It is obvious that,
if $\sA\in \LL\sb{k+1}$, then $\gen{\sA}$ is a $(k+1)$-dimensional code 
satisfying (b) and (c).
It is also obvious that,
if $\sA$ and $\sA\sprime$ are distinct standard bases  in  $\LL\sb{k+1}$,
then $\gen{\sA}$ and $\gen{\sA\sprime}$ are not $\Sn$-equivalent.
Therefore it is enough to show that,
for an arbitrary  $(k+1)$-dimensional code $\tC$
satisfying (b) and (c),
there exists an element of  $\LL\sb{k+1}$ that is a standard basis of $[\tC]$.
\par
Let $\sA=(A\sb 0, \dots, A\sb k)$ be a standard basis of $[\tC]$.
We put
$\sA\sprime :=(A\sb 0, \dots, A\sb {k-1})$.
Then $\gen{\sA\sprime}$ is a $k$-dimensional code satisfying (b) and (c).
Hence there exists a standard basis  $\sA\spar{i}\in \LL\sb{k}$ 
of the $\Sn$-equivalence class $[\gen{\sA\sprime}]$.
Let $\tau\in \Sn$ be an element that maps
the code $\gen{\sA\sprime}$ to $\gen{\sA\spar{i}}$.
We have
$$
\gen{ (\sA\spar{i}, \tau (A\sb k))}=\tau(\gen{(\sA\sprime, A\sb k)})= \tau (\gen{\sA})\in [\tC].
$$
Because $\sA\spar{i}$
is $\Sn$-increasing,
there exists $\sigma\in \Sn $ such that $\sigma (\sA\spar{i})=\sA\spar{i}$
and that
$$
\sigma ((\sA\spar{i}, \tau (A\sb k)))=(\sA\spar{i}, \sigma\tau (A\sb k))
$$
is $\Sn$-increasing.
Note that the sequence  $(\sA\spar{i}, \sigma\tau (A\sb k))$ is linearly independent,
because 
the code $\gen{(\sA\spar{i}, \sigma\tau (A\sb k))}=\sigma\tau ( \gen{ \sA})$
is of dimension $k+1$.
Note also  that $\card{\sigma\tau (A\sb k)}=\card{A\sb k}\le n/2$,
because $\sA=(\sA\sprime, A\sb{k})$ is a standard basis.
Therefore $(\sA\spar{i}, \sigma\tau (A\sb k))$ is a standard basis of 
the $\Sn$-equivalence class 
$$
[\gen{ (\sA\spar{i}, \sigma\tau (A\sb k))} ]=[\sigma\tau ( \gen{ \sA}) ]=[\tC].
$$
In other words, the word $ \sigma\tau ( A\sb k)$ appears in  $\AAA\spar{i}$.
Therefore we have a hoped-for standard basis  in $\LL\sb{k+1}$. 
\end{proof}
Starting with $\LL\sb{1}=\{ (\tZ ) \}$,
we can make the lists $\LL\sb{k}$ inductively.
\begin{remark}
By Proposition~\ref{prop:isin},
we can make the list of pairs $\sA\in \LL\sb k$ and $\sA\sprime\in \LL\sb{k\sprime}$ such that
$[\gen{\sA}]< [\gen{\sA\sprime}]$.
\end{remark}
\section{Geometry of splitting curves}\label{sec:splitting}
In this section, we assume $p=2$,
and fix a polynomial   $G\in \Utb$,
where $b$ is an even integer $\ge 4$.
\subsection{Definition of splitting curves and associated code words.}\label{subsec:defsplitting}
Let $C\subset\Pt$ be a reduced irreducible curve,
and let $\PT{C}$ be the proper transform of $C$ in $\XG$.
Since $\phiG:\XG\to\Pt$ is purely inseparable of degree $2$,
either one of the following holds;
\begin{itemize}
\item[(i)]
$\PT{C}$ is reduced and irreducible, or
\item[(ii)]
$\PT{C}=2\HPT{C}$,
where $\HPT{C}$ is a reduced irreducible curve on  $\XG$ birational to $C$ via $\phiG$.
\end{itemize}
\begin{definition}
We say that a reduced irreducible plane curve $C\st\Pt$ is {\it
splitting in $\XG$}
if (ii) above holds.
A reduced (but not necessarily irreducible) curve is said to be {\it splitting in $\XG$}
if every irreducible component of $C$ is splitting in $\XG$.
\end{definition}
\begin{definition}
Let $C\st\Pt$ be a reduced curve
 splitting in $\XG$.
We denote by $\HPT{C}$ the reduced divisor of $\XG$ such that
$2\HPT{C}$ is the proper transform of $C$ in $\XG$,
and by $\wG (C)\in \CCG$ the image of
the numerical equivalence class $[\HPT{C}]\in \SG$ by
$$
\SG \; \maprightsp{}\; \SG/\SG\sp 0=\tCCG\ \;\maprightsp{\rho\sb{\ZdG}}\; \CCG.
$$
\end{definition}
Let $C\subset \Pt$ be a reduced  curve splitting in $\XG$.
For a point $P\in \ZdG$,
let $m\sb P (C)$ denote the multiplicity of  $C$ at $P$.
Then we have
\begin{equation}\label{eq:Fm}
[\HPT{C}]=\frac{1}{2}\Bigl( -\sum\sb{P\in \ZdG} m\sb P (C) [\Gamma\sb P]\;+\; (\deg C)\, [\AG] \Bigr)
\end{equation}
in $\SG$.
Hence we have
\begin{equation}\label{eq:wm}
\wG (C) =\set{P\in \ZdG}{ m\sb P (C) \equiv 1 \,\bmod 2}.
\end{equation}
Suppose that $C$ is a union $C\sb 1\cup C\sb2$ of two splitting curves $C\sb1$ and $C\sb2$
that have no common irreducible components.
From~\eqref{eq:wm}, we have
\begin{equation}\label{eq:sum}
\wG (C\sb1 \cup C\sb2)=\wG(C\sb1) +\wG(C\sb2).
\end{equation}
\subsection{The general member of the linear system $|\IZdG (b-1)|$.}\label{subsec:stdsplit}
\begin{proposition}\label{prop:stdsplit}
The general member $C$ of $|\IZdG (b-1)|$ is splitting in $\XG$.
\end{proposition}
\begin{proof}
Recall that $C$ is reduced and irreducible by Corollary~\ref{cor:1secglobal}. 
By Proposition~\ref{prop:1secglobal}, 
there exist an affine part $U$ of $\Pt$ containing $\ZdG$ and affine coordinates $(x\sb 0, x\sb 1)$
on $U$ such that $C$ is defined by
$$
\varphi\sb{dG} (\theta\sb 0)=0,
$$
where $\theta\sb 0\in \Hz{\Theta(-1)}$ is given by $\theta\sb 0 |U=\der{}{x\sb 0}\otimes e\sb{-1}$.
If $G$ is written on $U$ in terms of $(x\sb 0, x\sb 1)$ as
$$
g(x\sb 0, x\sb 1)= \gamma\sb{00}(x\sb 0, x\sb 1)^2 +x\sb 0 \gamma\sb{10}(x\sb 0, x\sb 1)^2 + x\sb 1 \gamma\sb{01}(x\sb 0, x\sb 1)^2
+x\sb 0 x\sb 1 \gamma\sb{11}(x\sb 0, x\sb 1)^2,
$$
then $C$ is defined by
$$
\gamma\sb{10}^2 + x\sb1 \gamma\sb{11}^2 =0, 
$$
and $\ZdG$ is defined by
$$
\gamma\sb{10}^2 + x\sb 1 \gamma\sb{11}^2=\gamma\sb{01}^2 +x\sb 0 \gamma\sb{11}^2 =0.
$$
Note that $\gamma\sb{11}|\sb{C}$ is not zero,
because $\ZdG$ is reduced.
Hence we obtain
$$
g|\sb{C}= (\gamma\sb{00}^2 + x\sb1 \gamma\sb{01}^2)|\sb{C}=
\Bigl(\gamma\sb{00} +\frac{\gamma\sb{10}}{\gamma\sb{11}}\gamma\sb{01} \Bigr)^2\Big|\sb{C}.
$$
We put $\delta\sb C:=(\gamma\sb{00} +\gamma\sb{10}\gamma\sb{01}/\gamma\sb{11})|\sb C$.
The inverse image in $\XG$ of the generic point of $C$ is therefore isomorphic  to
$$
\Spec k(C) [w]/ (w + \delta\sb C)^2,
$$
which is not reduced.
Therefore $C$ is splitting in $\XG$.
\end{proof}
\begin{corollary}\label{cor:containsZdG}
The code $\CCG\subset \Pow (\ZdG)$ contains the word  $\ZdG$.
\end{corollary}
\begin{proof}
Because  $\ZdG$  is reduced,
the general member $C$ of $|\IZdG (b-1)|$ is smooth at each point of $\ZdG$.
Therefore we have $\wG(C)=\ZdG$ by~\eqref{eq:wm}.
\end{proof}
\begin{corollary}\label{cor:onSG}
The lattice $\SG$ is a $2$-elementary hyperbolic lattice of type I.
It is even if and only if $ b/2$ is odd.
\end{corollary}
\begin{proof}
The fact that $\SG$ is $2$-elementary and hyperbolic
follows from the fact that $\SG$ is an overlattice of $\SG\sp 0$.
Because $\ZdG\in\CCG$,
the lattice $\SG$ is of type I by Proposition~\ref{prop:typeIcode}.
(Note that $n=\card{\ZdG}$ is odd.)
Suppose that $b/2$ is odd.
Then $\SG$ is even by Propositions~\ref{prop:evencode} and~\ref{prop:RR}.
Suppose that $b/2$ is even.
Then $\card{\ZdG}\equiv 3\,\bmod 4$.
Because $\ZdG\in\CCG$,
the lattice  $\SG\cong \tS\sb{\ZdG} (\CCG)$ is not even
by Proposition~\ref{prop:evencode}.
\end{proof}
\subsection{Splitting curves with mild singularities.}\label{subsec:mild}
Let $C\st\Pt$ be a reduced (not necessarily irreducible) curve,
and $P$ a point of $C$.
Let $(\xi, \eta)$ be a formal parameter system of $\Pt$ at $P$.
\begin{definition}
Let $(a, b)$ be a pair of integers such that $a>b>1$
and that $a$ and $b$ are prime to each other.
We say that
$P$ is a {\it  cusp of $C$ of type $(a, b)$}
if $C$ is defined by $\xi^a+\eta^b=0$
locally at $P$ under a suitable choice of $(\xi, \eta)$.
A cusp of type $(3,2)$ is called an {\it ordinary cusp}.
Note that, if $P$ is a cusp of type $(a, b)$, then $C$ is locally irreducible at $P$.
\end{definition}
\begin{definition}
Let $m$ be a positive integer.
We say that $P$ is a {\it tacnode
of $C$ with tangent multiplicity $m$}
if  $C$ is defined by $\eta (\eta+ \xi^m)=0$
locally at $P$ under a suitable choice of $(\xi, \eta)$.
A tacnode with tangent multiplicity $1$ is called an {\it ordinary node}.
\end{definition}
\begin{proposition}\label{prop:mildsing}
Let $C\st\Pt$ be a reduced  curve
splitting in $\XG$,
and let $P$ be a point of $C$.
\par
{\rm (1)}
Suppose that $P\in C$ is a cusp of type $(a, b)$.
Then $P\in \ZdG$ if and only if $a+b \equiv 0 \,\bmod 2$.
\par
{\rm (2)}
Suppose that $P\in C$ is a tacnode with tangent multiplicity $m$. 
Then $P\in \ZdG$ if and only if $m \equiv 1 \,\bmod 2$.
\end{proposition}
\begin{proof}
Let $(\xi,\eta)$ be a formal parameter system of $\Pt$ at $P$.
We fix a global section $e\sb{b/2}$ of the line bundle $\MMM\cong\OPt(b/2)$
that is not zero at $P$.
The global section $G$ of $\LLL=\MMM\sp{\otimes 2}$
is given  by
$$
\gamma (\xi, \eta) \cdot e\sb{b/2}\sp{\otimes 2}
$$
locally at $P$, where   $\gamma (\xi, \eta)$ is a formal power series of $\xi$ and $\eta$,
which we write as 
$$
\gamma (\xi, \eta)=\sum  c\sb{ij} \xi\sp i \eta\sp j\qquad (c\sb{ij}\in k).
$$
The subscheme $\ZdG$ is defined by
$$
\Der{\gamma}{\xi}=\Der{\gamma}{\eta}=0
$$
locally at $P$.
\par
(1)
We choose $(\xi, \eta)$ in such a way that $C$ is defined by $\xi^a+\eta^b=0$
locally at $P$.
Then
$$
t\mapsto (\xi, \eta)=(t^b, t^a)
$$
is a normalization of $C$ at $P$.
Since $C$ is splitting in $\XG$,
the formal power series $\gamma (t^b, t^a)$
has a square root in the ring $k[[t]]$
of formal power series of $t$.
Suppose that $a+b$ is even.
Then both  $a$ and $b$ are odd,
because $a$ and $b$ are prime to each other.
Looking at the coefficients 
of $t^a$ and $t^b$ in $\gamma (t^b, t^a)$,
we obtain $c\sb{10}=c\sb{01}=0$.
Hence $P\in \ZdG$.
Suppose that $a+b$ is odd.
Looking at the coefficient of $t^{a+b}$
in $\gamma (t^b, t^a)$,
we obtain $c\sb{11}=0$.
If $P\in \ZdG$,
then $c\sb{11}=0$ implies that $\ZdG$ would fail to be  reduced of dimension $0$  at $P$.
Hence $P\notin \ZdG$.
\par
(2)
We choose $(\xi, \eta)$ in such a way that $C$ is defined by $\eta (\eta+\xi^m)=0$
locally at $P$.
Since $C$ is splitting in $\XG$,
both $\gamma (t, 0)$ and $\gamma (t, t^m)$
have square roots  in $k[[t]]$.
From $\sqrt{\gamma (t, 0)}\in k[[t]]$, 
we obtain $c\sb{10}=0$.
Suppose that $m$ is odd.
Then we also obtain $c\sb{m0}=0$
from $\sqrt{\gamma (t, 0)}\in k[[t]]$.
Looking at the coefficient of $t^m$ in $\gamma (t, t^m)$,
we have $c\sb{m0}+c\sb{01}=0$.
Therefore we have $P\in \ZdG$.
Suppose that $m$ is even.
Then we  obtain $c\sb{m+1, 0}=0$
from $\sqrt{\gamma (t, 0)}\in k[[t]]$.
Looking at the coefficient of $t^{m+1}$ in $\gamma (t, t^m)$,
we have $c\sb{m+1, 0}+c\sb{11}=0$.
Therefore $c\sb{11}=0$ follows and hence $P\notin \ZdG$.
\end{proof}
\begin{corollary}\label{cor:nodeandcusp}
Let $C\st\Pt$ be a reduced curve splitting in $\XG$.
If $P\in C$ is  an ordinary node,
then $P\in \ZdG$.
If $P\in C$ is an ordinary cusp, then $P\notin \ZdG$.
\end{corollary}
\begin{proposition}\label{prop:normalize}
Let $C\st \Pt $ be a reduced irreducible curve splitting in $\XG$.
Suppose that $C$ has
ordinary nodes and ordinary cusps 
as its only singularities.
Then the morphism $\phiG|\sb{\HPT{C}} : \HPT{C} \to C$ is the  normalization of $C$.
\end{proposition}
\begin{proof}
Suppose that $P\in C$ is an ordinary node.
Then $P\in \ZdG$ by Corollary~\ref{cor:nodeandcusp}.
The curve $\HPT{C}$ intersects $\Gamma\sb P$ at distinct two points,
and $\HPT{C}$ is smooth at each of these points.
\par
Suppose that $P\in C$ is an ordinary cusp.
Since $P\notin \ZdG$ by  Corollary~\ref{cor:nodeandcusp},
there exists a unique point $Q$ of $\XG$
such that $\phiG  (Q)=P$.
We  choose a formal parameter system $(\xi, \eta)$ of $\Pt$ at $P$ 
 so that $C$ is defined by $\xi^3+\eta^2=0$
locally at $P$,
and let $\gamma (\xi, \eta)$ be the formal power series introduced  in the proof of Proposition~\ref{prop:mildsing}.
Then $\XG$ is defined by
$$
w^2= \gamma (\xi, \eta)
$$
locally at $Q$,
where $w$ is a fiber coordinate of $\MMM$.
Since $\sqrt{\gamma (t^2, t^3)}\in k[[t]]$,
we have 
$$
\Der{\gamma}{\eta} (0, 0)=c\sb{01}=0.
$$
Therefore
the pair $(w-w(Q), \eta)$ is  a formal parameter system
of $\XG$ at $Q$.
Moreover,
we have $c\sb{10}\ne 0$
because $P\notin \ZdG$.
We put
$$
\beta (t):=\sqrt{\gamma(t^2, t^3)}=b\sb 0 + b\sb1 t + \dots.
$$
The curve $\HPT{C}$ is given by
$w=\beta (t)$ and $\eta=t^3$ at $Q$.
Since $c\sb{10}\ne 0$,
we have $b\sb 1\ne 0$,
which implies that $\HPT{C}$ is smooth at $Q$.
\end{proof}
\begin{proposition}\label{prop:num}
Let $C$ be a reduced {\rm (}possibly reducible{\rm )}
curve of degree $d$ that is splitting in $\XG$.
Suppose that $C$ has only ordinary nodes and ordinary cusps as its singularities.
Then we have
\begin{equation}\label{eq:num}
\card{\wG (C)}=d(b-d)+4\kappa,
\end{equation}
where $\kappa$ is the number of ordinary cusps on $C$.
\end{proposition}
\begin{proof}
Let $\Nodes (C)$
denote the set of ordinary nodes of $C$.
By~\eqref{eq:Fm},~\eqref{eq:wm} and Corollary~\ref{cor:nodeandcusp},
the assumption on the singularities of $C$ implies that
\begin{eqnarray}
&&\wG (C)=\set{P\in C\cap\ZdG}{\hbox{$C$ is smooth at $P$}},\label{eq:num1}\\
&& C\cap\ZdG =\wG (C)\;\sqcup\;  \Nodes (C),  \quand \label{eq:num2}\\
&& [\HPT{C}]=\frac{1}{2}(- \sum\sb{P\in \wG(C)} [\Gamma\sb P] - 2 \sum\sb{P\in \Nodes (C)} [\Gamma\sb P] \;+\;
d\,[\AG]). \label{eq:num3}
\end{eqnarray}
We prove~\eqref{eq:num} by induction on the number of  irreducible components of $C$.
Suppose that $C$ is irreducible.
Since $\HPT{C}$ is the normalization of $C$  by Proposition~\ref{prop:normalize},
the geometric genus of $C$ is given by
\begin{equation}\label{eq:genus}
\frac{1}{2}(d-1)(d-2)-\kappa-\card{\Nodes (C)} =\frac{1}{2}\HPT{C} (\HPT{C}+\KG) +1,
\end{equation}
where $\KG$ is the canonical divisor of $\XG$.
By~Proposition~\ref{prop:KX} and~\eqref{eq:num1},~\eqref{eq:num3},
we obtain~\eqref{eq:num}.
Suppose that $C$ is a union of two splitting curves
$C\sb 1$ and $C\sb 2$ that have  no common irreducible components.
Let $d\sb i$ be the degree of $C\sb i$,
and $\kappa\sb i$ the number of ordinary cusps of $C\sb i$.
We have $d=d\sb 1+d\sb 2$
and $\kappa=\kappa\sb 1+\kappa\sb 2$.
By the induction hypothesis,
we have
$\card{\wG (C\sb i)} = d\sb i (b-d\sb i) + 4\kappa\sb i$ for $i=1, 2$.
%
%
%\begin{equation}\label{eq:num4}
%\card{\wG (C\sb i)} = d\sb i (b-d\sb i) + 4\kappa\sb i\qquad\textrm{for $i=1, 2$}.
%\end{equation}
%
%
By~\eqref{eq:sum},
we have
\begin{equation}\label{eq:num5}
\card{\wG (C)} = \card{\wG(C\sb 1)}+\card{\wG (C\sb 2)}- 2\card{\wG(C\sb 1)\cap \wG (C\sb 2) }.
\end{equation}
Suppose that $P\in \wG(C\sb 1)\cap \wG (C\sb 2) $.
Then $P\in C\sb 1\cap C\sb 2$ by~\eqref{eq:num1}.
Suppose that $P\in C\sb 1\cap C\sb 2$.
Then $P$ is an ordinary node of $C$
and hence is contained in $\ZdG$ by Corollary~\ref{cor:nodeandcusp}.
Therefore
$P$ is contained in $\wG(C\sb 1)\cap \wG(C\sb 2)$ by~\eqref{eq:num1}.
Thus we obtain
$$
\wG(C\sb 1)\cap \wG (C\sb 2)=C\sb 1\cap C\sb 2,
$$
which implies $\card{\wG(C\sb 1)\cap \wG (C\sb 2)}=d\sb1  d\sb 2$.
Putting this into~\eqref{eq:num5}
and using the induction hypothesis, %~\eqref{eq:num4},
we obtain~\eqref{eq:num}.
\end{proof}
\begin{remark}\label{rem:cusp}
Let $G\in \Utb$ be chosen generally.
Then the general member of the linear system $|\IZdG (b-1)|$
has $(b-2)^2/4$ ordinary cusps as its only singularities.
Indeed,
we choose homogeneous coordinates $[X\sb0, X\sb 1, X\sb 2]$
generally so that the member $C$ of $|\IZdG (b-1)|$
defined by $\der{G}{X\sb 2}=0$ is general.
We write $G$ as
$$
X\sb 0 ^2 \Gamma\sb{00}^2 +
X\sb 1 ^2 \Gamma\sb{11}^2 +
X\sb 2 ^2 \Gamma\sb{22}^2 +
X\sb 0 X\sb 1 \Gamma\sb{01}^2 +
X\sb 1 X\sb 2 \Gamma\sb{12}^2 +
X\sb 2 X\sb 0 \Gamma\sb{20}^2, 
$$
where $\Gamma\sb{ij}$ are homogeneous polynomials of degree $(b-2)/2$.
Then $C$ is defined by 
$$
X\sb 1 \Gamma\sb{12}\sp 2 + X\sb 0 \Gamma\sb{20}\sp 2 =0.
$$
Since $G$ and  $[X\sb0, X\sb 1, X\sb 2]$ are general,
the homogeneous polynomials 
$\Gamma\sb{12}$ and $\Gamma\sb{20}$ are also general.
Hence $\Sing (C)$ consists of $(b-2)^2/4$ ordinary cusps
located at the intersection points of the curves 
defined by $\Gamma\sb{12}=0$ and $\Gamma\sb{20}=0$.
The equality~\eqref{eq:num} becomes
$$
n=b-1 + (b-2)^2
$$
in this case.
The linear system $|\IZdG (b-1)|$
gives a generalization of Serre's example~\cite[Chapter 3, Section 10, Exercise 10.7]{Hartshorne}
of  linear systems of plane curves  with moving singularities in positive characteristics.
\end{remark}
\subsection{Splitting curves with only ordinary nodes.}
\begin{proposition}\label{prop:GC}
Let $G\sb C$ and $G\sb D$ be homogeneous polynomials 
defining plane curves $C$ and $D$
such that $\deg G\sb C +\deg G\sb D=b$.
Suppose that $G\sb C G\sb D$ is a polynomial contained  in $k\sptimes G +\Vtb$.
Then the following hold; %{;}
\begin{itemize}
\item[(i)]
$C$ and $D$ are reduced and have no common irreducible components,
\item[(ii)]
$C\cup D$ has only  ordinary nodes as its singularities,
\item[(iii)]
$C$ and $D$ are splitting in $\XG$, and
\item[(iv)]
$\wG (C)=\wG (D)=C\cap D$.
\end{itemize}
\end{proposition}
\begin{proof}
The assertions~(i) and (ii) follow from Proposition~\ref{prop:Upb}.
The assertion~(iii) is obvious because $X\sb{G\sb C G\sb D}$ is isomorphic to $\XG$ over $\Pt$.
By Corollary~\ref{cor:nodeandcusp},
we have $C\cap D\st \ZdG$.
Since $C$ and $D$ are smooth at each point of $C\cap D$,
we have $C\cap D\subset\wG (C)$ and $C\cap D\subset\wG (D)$
by~\eqref{eq:wm}.
From Proposition~\ref{prop:num},
we have 
$$
\card{\wG (C)}=\card{\wG (D)}=\deg C\cdot\deg D= \card{C\cap D}.
$$
Therefore (iv) holds.
\end{proof}
The converse of Proposition~\ref{prop:GC} is also true:
\begin{proposition}\label{prop:convGC}
Let $C$ be a curve defined by $G\sb C=0$.
Suppose that $C$ is reduced,
has only  ordinary nodes as its singularities,
and is splitting in $\XG$.
Then there exists a homogeneous polynomial $G\sb D$
of degree $b-\deg G\sb C$
such that $G\sb C G\sb D$ is  contained  in $k\sptimes G +\Vtb$.
\end{proposition}
\begin{proof}
First note that the degree of $G\sb C$ is $\le b$ by Proposition~\ref{prop:num}.
Let $\Nodes (C)$ denote the set of  ordinary nodes of $C$,
and let $\nu :\wtC\to C$ be the normalization of $C$,
that is, $\wtC$ is the disjoint union of normalizations of
irreducible components of $C$.
For $P\in \Nodes (C)$,
let $P\sb 1$ and $P\sb 2$ denote the points of $\wtC$  that are mapped to  $P$ by $\nu$.
Consider the following commutative diagram:
$$
\begin{matrix}
\Hz{\MMM} &\maprightsp{\res} & H\sp 0 (C, \MMM|\sb C) & \maprightsp{\nu\sp *\sb{\MMM}} & H\sp 0 (\wtC, \nu\sp * \MMM|\sb C) &\\
\mapdown & & \mapdown & & \mapdown & \\
\Hz{\LLL} &\maprightsp{\res} & H\sp 0 (C, \LLL|\sb C) & \maprightsp{\nu\sp *\sb{\LLL}} & H\sp 0 (\wtC, \nu\sp * \LLL|\sb C), &
\end{matrix}
$$
where the left horizontal arrows are restrictions,
the right horizontal arrows are the pull-backs by $\nu$,
and the vertical arrows are the squaring map $f\mapsto f^2$.
For each $P\in \Nodes (C)$,
we have canonical isomorphisms of $1$-dimensional vector spaces
\begin{equation}\label{eq:canonicalisom}
\nu\sp * \MMM|\sb C \otimes k(P\sb 1)\cong\nu\sp * \MMM|\sb C \otimes k(P\sb 2),
\quad
\nu\sp * \LLL|\sb C \otimes k(P\sb 1)\cong\nu\sp * \LLL|\sb C \otimes k(P\sb 2),
\end{equation}
where $k (P\sb i)$ is the residue field of $\OOO\sb{\wtC}$ at $P\sb i\in \wtC$.
The homomorphisms $\nu\sp *\sb{\MMM}$ and $\nu\sp *\sb{\LLL}$
are injective,
and their images coincide with the spaces of all sections $f$  that
satisfy
$f(P\sb 1)=f(P\sb 2)$
for every $P\in\Nodes (C)$,
where $f(P\sb 1)$ and $f(P\sb 2)$ are compared by the canonical isomorphisms~\eqref{eq:canonicalisom}.
Consider the images $g\in H\sp 0 (C, \LLL|\sb C)$
and $\tilde g \in H\sp 0 (\wtC, \nu\sp *\LLL |\sb C)$ of $G\in\Hz{\LLL}$.
We have
\begin{equation}\label{eq:tlg}
\tilde g (P\sb 1)=\tilde g(P\sb 2)
\qquad\textrm{for any $P\in \Nodes (C)$}.
\end{equation}
Because $C$ is splitting in $\XG$,
there exists a global section $\tilde h\in H\sp 0 (\wtC, \nu\sp *\MMM |\sb C)$
such that $\tilde h^2 =\tilde g$.
By~\eqref{eq:tlg},
we have
$\tilde h (P\sb 1)=\tilde h (P\sb 2)$
for each $P\in\Nodes (C)$.
Hence there exists $h\in H\sp 0 (C, \MMM|\sb C)$
such that $\nu\sb{\MMM}\sp * (h)=\tilde h$.
Then we have $g=h^2$ because $\nu\sb{\LLL}\sp *$ is injective.
Since the restriction homomorphism 
$\Hz{\MMM} \to H\sp 0 (C, \MMM|\sb C)$
is surjective,
there exists $H\in \Hz{\MMM}$
such that $(G+H^2)|\sb C=0$.
Then the polynomial  $G+H^2$ is divisible  by $G\sb C$.
\end{proof}
\subsection{Splitting lines and splitting smooth conics}\label{subsec:linesandconics}
\begin{proposition}\label{prop:convnum}
{\rm (1)}
Let $L\subset \Pt$ be a line.
If $\card{L\cap \ZdG}> (b-2)/2$,
then $L$ is splitting in $\XG$.
{\rm (2)}
Let $Q\subset \Pt$ be a smooth conic.
If $\card{Q\cap \ZdG}> b-1$,
then $Q$ is splitting in $\XG$.
\end{proposition}
\begin{proof}
(1)
We choose a general line $ l \sb{\infty}\st\Pt$,
and fix affine coordinates $(x\sb 0, x\sb 1)$
on $U:=\Pt\sm l \sb{\infty}$
such that $L$ is defined by $x\sb 1=0$.
Let us consider $x\sb 0$ as an affine parameter of $L$.
We express $G$ on $U$ by
\begin{equation}\label{eq:GU}
\gamma\sb{00} (x\sb 0, x\sb 1)\sp2 
+ x\sb 0 \gamma\sb{10} (x\sb 0, x\sb 1)\sp2 
+ x\sb 1 \gamma\sb{01} (x\sb 0, x\sb 1)\sp2 
+ x\sb 0 x\sb 1 \gamma\sb{11} (x\sb 0, x\sb 1)\sp2.
\end{equation}
Then $L\cap \ZdG$ is defined on $L$ by
$$
\gamma\sb{10} (x\sb 0, 0)^2 =\gamma\sb{01}(x\sb 0, 0)^2 + x\sb 0 \gamma\sb{11} (x\sb 0, 0)^2 =0.
$$
Note that the degree of $\gamma\sb{10}$ is at most $(b-2)/2$. 
Hence the assumption
$\card{L\cap \ZdG}>(b-2)/2$ implies that
$\gamma\sb{10} (x\sb 0, 0)$ is constantly equal to zero.
Therefore  $\gamma\sb{10} (x\sb 0, x\sb 1)$ can be written as
$x\sb 1 \delta\sb{10} (x\sb 0, x\sb 1)$.
Then $G$ is equal to 
$$
\gamma\sb{00}^2 + x\sb 1 ( x\sb 0 x\sb 1 \delta\sb{10}^2 +\gamma\sb{01}^2 + x\sb 0 \gamma\sb{11}^2)
$$
on $U$.
Hence $L$ is splitting in $\XG$.
\par
(2)
Let $ l \sb{\infty}$
be a general tangent line to $Q$,
and let $(x\sb 0, x\sb 1)$ be affine coordinates on $U=\Pt\sm l \sb{\infty}$
such that $Q$ is defined by $x\sb 1 + x\sb 0^2=0$.
We consider $x\sb 0$ as an affine parameter of $Q$.
Again we write $G$ on $U$ as in~\eqref{eq:GU}.
Then $Q\cap\ZdG$ is defined on $Q$ by
$$
\gamma\sb{10} (x\sb 0, x\sb 0^2)^2 + x\sb 0^2 \gamma\sb{11}(x\sb 0, x\sb 0^2)^2
=
\gamma\sb{01}(x\sb 0, x\sb 0^2)^2 + x\sb 0 \gamma\sb{11}(x\sb 0, x\sb 0^2)^2=0.
$$
Since the degrees of $\gamma\sb{10}$ and $\gamma\sb{11}$ are at most $(b-2)/2$,
the number of the roots of
\begin{equation*}\label{eq:gammasq}
\gamma\sb{10} (x\sb 0, x\sb 0^2)^2 + x\sb 0^2 \gamma\sb{11}(x\sb 0, x\sb 0^2)^2=
(\gamma\sb{10} (x\sb 0, x\sb 0^2) + x\sb 0 \gamma\sb{11}(x\sb 0, x\sb 0^2))^2
\end{equation*}
is at most $b-1$.
Consequently 
the assumption
$\card{Q\cap \ZdG}>b-1$ implies that
$(\gamma\sb{10}+x\sb 0 \gamma\sb{11})|\sb Q=0$.
Then $G|\sb Q$ is written as
$$
\gamma\sb{00}(x\sb 0, x\sb 0^2)^2 + x\sb 0^2 \gamma\sb{01}(x\sb 0, x\sb 0^2)^2,
$$
which is the  square of $(\gamma\sb{00}+x\sb 0\gamma\sb{01})|\sb Q$.
Therefore $Q$ is splitting.
\end{proof}
\begin{corollary}\label{cor:convnum}
{\rm (1)}
If $L\subset \Pt$ is a line, then 
$\card{L\cap \ZdG}$ is either $\le  (b-2)/2$ or $b-1$.
{\rm (2)}
If $Q\subset \Pt$ is a smooth conic, then 
$\card{Q\cap \ZdG}$ is either $\le  b-1$ or $2 (b-2)$.
\end{corollary}
\begin{example}\label{example:GDK}
Let $q=2\sp\nu$ be a power of $2$.
We put $b:=q+2$, and consider the homogeneous polynomial
$$
G\sb{\DK, q}=X\sb 0 X\sb 1 X\sb 2 (X\sb 0\sp{q-1} + X\sb 1 \sp{q-1} + X\sb 2 \sp{q-1})
$$
of degree $b$,
which is a generalization of Dolgachev-Kondo's polynomial~\eqref{eq:DK} of degree $6$.
It is easy to see that $Z(dG\sb{\DK, q})$
consists of all $\F\sb q$-rational points of $\Pt$.
Because $n=b^2-3b+3 = q^2+q+1$ is equal to the cardinality of $\Pt(\F\sb q)$,
the polynomial $G\sb{\DK, q}$ is a member of $\Utb$.
Every $\F\sb q$-rational line contains $q+1=b-1$ points of $Z(dG\sb{\DK, q})$,
and hence is splitting in $X\sb{G\sb{\DK, q}}$. 
\end{example}
\section{Known facts about  $K3$ surfaces} \label{sec:ssk3}
\subsection{The Artin-Rudakov-Shafarevich theory}
Let $p$ be an arbitrary prime integer, and 
$X$  a supersingular $K3$ surface in characteristic $p$.
Artin~\cite{Artin74} showed that 
the discriminant  of the numerical N\'eron-Severi lattice
$\NS\sb X$ of $X$ is equal to  $-p\sp{2\sigma}$,
where $\sigma$ is a positive integer $\le 10$.
This integer $\sigma$ is called the {\it Artin invariant} of $X$.
\begin{proposition}[Artin~\cite{Artin74},  Rudakov-Shafarevich~\cite{RS_char2}, Shioda~\cite{Shioda78}]\label{prop:exist}
For any pair $(p, \sigma)$ of a prime integer $p$ and a positive integer $\sigma\le 10$,
there exists a supersingular $K3$ surface in characteristic $p$
with Artin invariant $\sigma$.
\end{proposition}
For an integer $\sigma$ with $1\le\sigma\le 10$,
let $\RSts$ denote the lattice with the following properties;
\begin{itemize}
\item[]\hskip -15pt (RS1) even, hyperbolic, and  of rank $22$,
\item[]\hskip -15pt (RS2) $2$-elementary of type I, and 
\item[]\hskip -15pt (RS3) $\disc \RSts=-2^{2\sigma}$.
\end{itemize}
\begin{proposition}[Rudakov-Shafarevich~\cite{RS}]\label{prop:RSlattice}
The conditions {\rm (RS1)-(RS3)} determine the lattice $\RSts$ uniquely up to isomorphisms.
\end{proposition}
\begin{proposition}[Rudakov-Shafarevich~\cite{RS}]\label{prop:RSnef}
Let $X$ be a supersingular $K3$ surface in characteristic $2$
with Artin invariant $\sigma$.
Then the lattice $\NS\sb X$ is isomorphic to $\RSts$.
More precisely,
let $v\in \RSts$ be a vector with $v^2>0$.
Then there exists an isometry
$\phi$ from $\RSts$ to $\NS\sb X$ such that $\phi (v)$ is the class $[H]$ of a nef line bundle $H$ of $X$.
\end{proposition}
\subsection{$K3$ surfaces as sextic double planes}
Let $T$ be a negative definite even lattice.
A vector $v\in T$ is called a {\it root} if $v^2=-2$.
We put
$$
\Roots (T):=\set{v\in T}{v^2=-2}.
$$
It is well-known that $\Roots (T)$ forms a root system of type $ADE$~(\cite{B, E}).
\begin{definition}
A pair $(X, H)$ of a $K3$ surface $X$ and a  line bundle   $H$ of $X$ with $H^2=2$ and $|H|\ne\emptyset$
is called a {\it sextic double plane}
if the complete linear system $|H|$ is fixed component free.
If $(X, H)$ is a sextic double plane, then  $|H|$  defines
a generically finite morphism
$$
\map{\Phi\sb{|H|}}{X}{\Pt}
$$
of degree $2$.
\end{definition}
For a sextic double plane $(X, H)$, we denote by
$$
X\;\to\; Y\sb{|H|}\;\to\;\Pt
$$
the Stein factorization of $\Phi\sb{|H|}$.
The normal $K3$ surface $\YH$ has only rational double points
as its singularities.
We denote by $R(X, H)$
the $ADE$-type of the singular points of $\YH$,
that is,
$R(X, H)$ is the type of  the $ADE$-configuration of 
$(-2)$-curves that are contracted by $X\to\YH$.
\begin{remark}
Let $(X, H)$ be a sextic double plane.
We have
\begin{equation}\label{eq:stein}
\YH:=\SSpec \Phi\sb{|H| *} \OOO\sb X \;\cong\; \Proj \left( \bigoplus\sb{m=0}\sp{\infty} H\sp 0 (X,  H \sp{\otimes m}) 
\right).
\end{equation}
Indeed,
let $s $ be a non-zero element of $\Hz{\OPt(1)}$,
and
let $s\sb X$ be the  global section $\Phi\sb{|H|}\sp * (s)$ of $H$.
We put $U:=\{ s\ne 0\}\st\Pt$.
Then
the module 
$\Gamma (U, \Phi\sb{|H| *} \OOO\sb X )$
of sections of $\OOO\sb X$ over $\Phi\sb{|H|}\inv (U)=\{s\sb X\ne 0\}\st X$ 
is canonically isomorphic to the degree $0$ part of the graded ring
$$
\bigoplus\sb{m=0}\sp{\infty} H\sp 0 (X, H\sp{\otimes m })\left[\frac{1}{s\sb X}\right].
$$
Hence the isomorphism~\eqref{eq:stein} holds.
\par
The graded ring $\oplus\sb{m=0}\sp{\infty} H\sp 0 (X, H\sp{\otimes m})$
is generated by elements $X\sb 0, X\sb 1, X\sb 2$ of degree $1$ and an element $w$ of degree $3$,
and the relations are generated by 
\begin{equation}\label{eq:CG}
w^2 + C (X\sb 0, X\sb 1, X\sb 2) w + G(X\sb 0, X\sb 1, X\sb 2)=0.
\end{equation}
where $C$ and $G$ are homogeneous polynomials of degree $3$ and $6$, respectively.
Hence $\YH$ is defined by~\eqref{eq:CG}
in the weighted projective space $\P (3,1,1,1)$.
\end{remark}
%
%
%\begin{remark}
%If $p\ne 2$, then we can make $C=0$ by the coordinate change $w\mapsto w+C/2$,
%and $R(X, H)$ coincides with the $ADE$-type of the singularities of the plane curve defined by $G=0$.
%\end{remark}
%
%
%
%
\begin{proposition}[Urabe~\cite{U}, Nikulin~\cite{N2}]\label{prop:nu}
Let $X$ be a $K3$ surface and $H$ a line bundle on $X$ with $H^2=2$.
\par
{\rm (1)}
The pair $(X, H)$ is a sextic double plane if and only if $H$ is nef and the set
$\shortset{u\in \NS\sb X}{u^2=0, u\cdot[H]=1}$
is empty.
\par
{\rm (2)}
Suppose that $(X, H)$ is a sextic double plane.
Then $R(X, H)$ coincides with the $ADE$-type
of the root system $\Roots ([H]\sperp)$,
where $[H]\sperp\subset \NS\sb X$
is the orthogonal complement of $[H]$ in $\NS\sb X$.
More precisely,
the classes of  $(-2)$-curves contracted by $X\to\YH$ form 
a simple root system of $\Roots ([H]\sperp)$.
\end{proposition}
Proposition~\ref{prop:nu} is true in any characteristic.
Indeed,
the proof of Proposition~1.7 in~\cite{U} can be transplanted in any characteristic except for the use of
the Kawamata-Vieweg vanishing theorem,
which can be replaced by Proposition~0.1 in~\cite{N2}.
\subsection{Purely inseparable sextic double planes}
The following is obvious:
\begin{proposition}
If $G$ is a polynomial in $ \Uts$, then $(\XG, \AG)$ is a sextic double plane,
and $R(\XG, \AG)=21 A\sb 1$  holds.
\end{proposition}
Conversely, we have the following:
\begin{proposition}[\cite{Shimada2003}]\label{prop:21A1}
Let $(X, H)$ be a sextic double plane.
If $R(X, H)=21 A\sb 1$, then $p=2$ and the morphism
$\Phi\sb{|H|} : X\to\Pt$ is purely inseparable.
\end{proposition}
Let $(X, H)$ be a sextic double plane such that $R(X, H)=21 A\sb 1$.
Then there exists a homogeneous polynomial $G(X\sb 0, X\sb 1, X\sb 2)$ of degree $6$ such that $\YH$ is defined by
$w^2=G$.
Since $\YH$ has rational double points of type $21A\sb 1$ as its only singularities,
we have $G\in \Uts$.
\begin{corollary}\label{cor:21A1}
If $(X, H)$ is a sextic double plane 
with  $R(X, H)=21A\sb 1$,
then there exists $G\in \Uts$ such that $X=\XG$, $H=\AG$, $\YH=\YG$ and $\Phi\sb{|H|}=\phiG$.
\end{corollary}
\section{The list of geometrically realizable classes of codes}\label{sec:codelist}
In this section,
we study the case where $p=2$ and $b=6$.
\subsection{A characterization of geometrically realizable classes of codes.}
\begin{theorem}\label{thm:wt}
Let $\tZ$ be a set with $\card{\tZ}=21$,
and let $\tC\subset\Pow (\tZ)$ be a code.
The $\St$-equivalence class $[\tC]$ containing $\tC$ 
is geometrically realizable if and only if the following hold: %{:}
\begin{itemize}
\item[(a)] $\dim\tC\le 10$,
\item[(b)] $\tZ\in \tC$, and
\item[(c)] $\card{A}\in \{ 0,5,8,9,12,13,16,21\}$ for any $A\in \tC$.
\end{itemize}
\end{theorem}
\begin{proof}
Suppose that $[\tC]$ is geometrically realizable,
and let $G$ be a polynomial  in  $\Uts$ such that $\tC\cong \CCG$ by 
some bijection $\tZ\cong \ZdG$.
We have
$$
|\disc\SG|=2^{22-2\dim{\tCCG}}=2^{22-2\dim\tC}.
$$
Since the Artin invariant of $\XG$ is positive,
we have $\dim \tC\le 10$.
By Corollary~\ref{cor:containsZdG},
we have $\ZdG\in\CCG$,
and hence $\tZ\in \tC$.
By Proposition~\ref{prop:RR},
$\card{A} \,\bmod 4$ is either $0$ or $1$ for any $A\in \CCG$.
Therefore,
in order to show that $\tC$ satisfies (c),
it is enough to show that $\card{A}\notin \{1,4,17,20\}$ for any $A\in \CCG$.
Suppose that there is an element $A\in\CCG$ with $\card{A}=1$.
Then there exists $P\in \ZdG$ such that $(\{P\}, 1)$ is  contained  in the lift 
$\tCCG=\CCG\lift$ of $\CCG$.
Hence the vector
$$
v:=\frac{1}{2} (-[\Gamma\sb P]\;+\; [\AG])
$$
is  contained  in $\SG$.
Because $v\cdot [\AG]=1$ and $v^2=0$,
we see from Proposition~\ref{prop:nu}
that $(X\sb G, \AG)$ is not a sextic double plane,
which is absurd.
Suppose that there is a word  $A\in\CCG$ with $\card{A}=4$.
Then  $(A, 0)$ is a word  in  the lift $\CCG\lift$ of $\CCG$.
Hence the vector
$$
v:=\frac{1}{2} (\sum\sb{P\in A} [\Gamma\sb P])
$$
is   contained  in $\SG$.
Because $v\cdot [\AG]=0$ and $v^2=-2$,
the vector $v$ is an element of $\Roots ([\AG]\sperp)$.
However,
we see from Proposition~\ref{prop:nu}
that every vector in $\Roots ([\AG]\sperp)$ is written as a linear combination of $[\Gamma\sb P]$ $(P\in \ZdG)$
with {\it integer} coefficients.
Thus we get a contradiction.
Suppose that there is a word  $A\in\CCG$ with $\card{A}=17$ or $20$.
Then $Z(dG)+A\in \CCG$ is of weight $4$ or $1$,
which is impossible as has been shown above.
Therefore the code $\tC$ satisfies (a), (b) and (c).
\par
Suppose that $\tC$ satisfies (a), (b) and (c).
We put 
$$
\sigma:=11-\dim\tC.
$$
By Proposition~\ref{prop:evencode} and the property (c), the submodule
$\tS\sb{\tZ} (\tC)=\pr\sb{\tS\sp0\sb{\tZ}}\inv (\tC\lift)$
of $(\tS\sp0\sb{\tZ})\dual $
corresponding to the lift $\tC\lift\subset \DG(\tS\sp0\sb{\tZ})$ of $\tC$
is an even overlattice of $\tS\sp0\sb{\tZ}$.
\begin{claim}\label{claim:1}
The even overlattice   $\tS\sb{\tZ} (\tC)$ of $\tS\sb{\tZ}\sp 0$ 
is isomorphic to $\RSts$.
\end{claim}
\begin{proof}[Proof of Claim~\ref{claim:1}]
By Proposition~\ref{prop:RSlattice},
it is enough to show that $\tS\sb{\tZ} (\tC)$ satisfies 
the conditions (RS1), (RS2) and (RS3).
It is obvious that $\tS\sb{\tZ} (\tC)$ is $2$-elementary and hyperbolic.
By Proposition~\ref{prop:typeIcode},
 the property (b) implies that $\tS\sb{\tZ} (\tC)$ is of type I.
By~\eqref{eq:disc}, we have $|\disc (\tS\sb{\tZ} (\tC))|=2^{2\sigma}$.
\end{proof}
By Proposition~\ref{prop:exist},
there exists a 
supersingular $K3$ surface  $X$ in characteristic $2$
with Artin invariant $\sigma$.
In  $\tS\sb{\tZ} (\tC)$,
we have a vector $\tth$ with $\tth^2=2$.
By Proposition~\ref{prop:RSnef},
there exists an isometry
$$
\mapisom{\phi}{ \tS\sb{\tZ} (\tC)}{\NS\sb X}
$$
such that $\phi (\tth)$ is the class $[H]$ of a nef line bundle $H$ on $X$ with $H^2=2$.
\begin{claim}\label{claim:2}
The pair $(X, H)$ is a sextic double plane with $R(X, H)=21 A\sb 1$.
\end{claim}
\begin{proof}[Proof of Claim~\ref{claim:2}]
By Proposition~\ref{prop:nu}
and the isometry $\phi$,
it is enough to show that the set
\begin{equation}\label{eq:theset}
\set{ u\in \tS\sb{\tZ} (\tC)}{u^2=0, \, u\tth=1}
\end{equation}
is empty,
and that the $ADE$-type of the root system $\Roots (\tth\sperp)$ 
is $21 A\sb 1$,
where $\tth\sperp$ is the orthogonal complement of $\tth$ in $\tS\sb{\tZ} (\tC)$.
Suppose that a vector
$$
u=\frac{1}{2} (\sum\sb{\tP\in\tZ} a\sb {\tP} \te\sb {\tP} + b\tth ) \qquad(a\sb{\tP}\in \Z,\;\; b\in \Z)
$$ 
of $\tS\sb{\tZ} (\tC)$ is contained in the set~\eqref{eq:theset}.
Because $u\tth=1$, we have $b=1$.
Because $u^2=0$, we have
$\sum a\sb {\tP}^2=1$.
Hence $u$ is of the form $(\tth \pm \te\sb{\tP})/2$.
Its image in $\tC\lift$ by the natural projection $\tS\sb{\tZ} (\tC)\to \tS\sb{\tZ} (\tC)/\tS\sp0\sb{\tZ}$
therefore yields an element $(\{\tP\}, 1)\in \Pow (\tZ) \oplus \F\sb 2$.
This contradicts the property (c).
Let 
$$
r=\frac{1}{2} (\sum\sb{\tP\in\tZ} a\sb {\tP} \te\sb {\tP} + b\tth )\qquad(a\sb {\tP}\in \Z, \;\; b\in \Z)
$$ 
be a root of $\tth\sperp$.
Because $r\tth=0$, we have $b=0$.
Because $u^2=-2$, we have
$\sum a\sb {\tP}^2=4$.
Hence $r$ is 
either
$$
\pm \te\sb{\tP}, \quad\textrm{or}\quad \frac{1}{2}\sum\sb{\tP\in A} (\pm \te\sb{\tP}) \;\;\textrm{with}\;\; \card{A}=4.
$$
By the property (c) of $\tC$, the latter cannot occur.
Hence $\Roots (\tth\sperp)$ is equal to $\shortset{\pm \te\sb{\tP}}{\tP\in \tZ}$,
and its $ADE$-type is $21 A\sb 1$.
\end{proof}
By Corollary~\ref{cor:21A1},
there exists $G\in \Uts$
such that $X=\XG$, $H=\AG$ and $\Phi\sb{|H|}=\phiG$.
Note that the isometry
$$
\mapisom{\phi}{\tS\sb{\tZ} (\tC)}{\NS\sb X\cong \SG}
$$
maps $\Roots (\tth\sperp)$ to $\Roots ([\AG]\sperp)$
bijectively.
Composing the isometry  $\phi$
with reflections with respect to some $\te\sb{\tP}$ if necessary,
we can assume that $\phi$ maps each $\te\sb{\tP}\, (\tP\in \tZ)$ to $[\Gamma\sb{P}]$ for some $P\in \ZdG$.
The correspondence $\te\sb{\tP} \mapsto \Gamma\sb{P}$
gives us a bijection $\tZ\cong \ZdG$ that induces $\tC\cong \CCG$.
Hence the class  $[\tC]$ is geometrically realizable.
\end{proof}
\subsection{From the code to the configuration of splitting curves}
In this subsection,
we fix a polynomial  $G\in \Uts$ and 
show how to read from $\CCG$ the configuration of plane curves of degree $\le 3$  splitting in $\XG$.
\begin{definition}
For a word $A\in \Pow (\ZdG)$ with $\card{A}\in \{ 5, 8, 9\}$,
we put
$$
\deg A:=
\begin{cases}
1 & \textrm{ if $\card{A}=5$}, \\
2 & \textrm{ if $\card{A}=8$}, \\
3 & \textrm{ if $\card{A}=9$}. \\
\end{cases}
$$
We say that a word $A$ of $\CCG$ is {\it reducible in $\CCG$}
if there exist
words $A\sb1$ and $A\sb 2$ of $\CCG$ with $\card{A\sb 1}, \card{A\sb 2}\in \{5,8,9\}$
such that  $A=A\sb 1+A\sb 2$
and  $\deg A=\deg A\sb 1+\deg A\sb 2$ hold.
We say that  $A$ is  {\it irreducible in $\CCG$}
if $A$ is not reducible in $\CCG$.
\end{definition}
A word of $\CCG$ with weight $5$ is always irreducible in $\CCG$.
\begin{proposition}\label{prop:line}
The correspondence
$L\mapsto L\cap \ZdG$ gives a bijection from the set of lines $L\st\Pt$
splitting in $\XG$ to the set of words $A\in \CCG$ of weight  $5$.
\end{proposition}
\begin{proof}
Suppose that a line $L$ is splitting in $\XG$.
Then we have $\wG (L)= L\cap \ZdG$   by~\eqref{eq:wm}  and $\card{\wG (L)}=5$
by Proposition~\ref{prop:num}.
\par
Conversely suppose that a word $A\in \CCG$ with $\card{A}=5$ is given.
A line $L$ satisfying  $L\cap\ZdG=A$ is, if exists, obviously unique.
Because $(A, 1)$ is a word in  the lift $\CCG\lift$ of $\CCG$, we have a vector
\begin{equation*}\label{eq:linevec}
u:=\frac{1}{2} (-\sum\sb{P\in A} [\Gamma\sb P]\;+\; [\AG])
\end{equation*}
in $\SG$.
Because $u^2=-2$ and $u\cdot [\AG]=1$,
the class $u$ is represented by  an effective divisor $D$ of $\XG$.
Since $D\AG=1$,
there exists a reduced  irreducible component $D\sb 0$ of $D$
such that 
$\phiG:\XG\to\Pt$ induces a birational morphism from $D\sb 0$ to a line $L\st\Pt$.
Moreover $D-D\sb 0$ is 
a linear combination of  the curves $\Gamma\sb P$
with non-negative integer coefficients.
Since the proper transform of $L$ in $\XG$ is $2D\sb 0$, the line $L$ is splitting in $\XG$, 
and $\HPT{L}=D\sb 0$ holds.
Since  $u-[D\sb 0]$ is in $\SzG$,
we have
$$
(A, 1)= u \;\bmod \SzG = [D\sb 0] \;\bmod \SzG = [\HPT{L}] \; \bmod\SzG.
$$
Therefore 
we obtain  $A=\wG(L)=L\cap \ZdG$.
\end{proof}
\begin{remark}\label{rem:intersection_lines}
Let $L\sb 1$ and $L\sb 2$ be  distinct splitting lines.
By Corollary~\ref{cor:nodeandcusp},
we see that 
$\wG(L\sb 1)\cap \wG(L\sb 2)$ consists of one point,
which is the intersection point of  $L\sb 1$ and $L\sb 2$,
and the word $\wG(L\sb 1\cup L\sb 2)=\wG(L\sb 1)+ \wG(L\sb 2)$ is of weight $8$.
\end{remark}
\begin{remark}\label{rem:intersection_threelines}
Let $L\sb 1$, $L\sb 2$ and $L\sb 3$ be  distinct splitting lines.
The word 
$$
\wG(L\sb 1\cup L\sb 2\cup L\sb 3)=\wG(L\sb 1)+ \wG(L\sb 2)+\wG(L\sb 3)
$$
is of weight $9$
if $L\sb 1\cup L\sb 2\cup L\sb 3$ has only ordinary nodes as its singularities,
while this word is of weight $13$ if
$L\sb 1\cap L\sb 2\cap L\sb 3$ is non-empty.
\end{remark}
\begin{proposition}\label{prop:conic}
The correspondence
$Q\mapsto Q\cap \ZdG$ gives a bijection from the set of smooth conics  $Q\st\Pt$
splitting in $\XG$ to the set of  words $A\in \CCG$ of weight  $8$
irreducible in $\CCG$.
\end{proposition}
\begin{proof}
Suppose that   a smooth conic $Q$ is splitting in $\XG$.
Then the word $\wG (Q)= Q\cap \ZdG$
of $\CCG$ is of weight $8$ by Proposition~\ref{prop:num}.
If $\wG (Q)$ were reducible in $\CCG$,
then $Q\cap \ZdG$ would be  written as $A\sb 1 + A\sb 2$,
where
$A\sb 1$ and $A\sb 2$ are words of $\CCG$ with weight $5$.
By Proposition~\ref{prop:line},
the points in $A\sb i$ $(i=1, 2)$ are collinear,
and hence $Q$ would contain two sets of four collinear points,
which contradicts the assumption that $Q$ is smooth.
Hence the word $\wG (Q)$ is irreducible in $\CCG$.
\par
Suppose that $A\in \CCG$ is a word of weight $8$ that is irreducible in $\CCG$.
Since $(A, 0)\in \CCG\lift$,
the vector
\begin{equation*}\label{eq:conicvec1}
u:=\frac{1}{2} (-\sum\sb{P\in A} [\Gamma\sb P]\;+\; 2\,[\AG])
\end{equation*}
of $(\SzG)\dual$ 
is contained in $\SG$.
Because $u^2=-2$ and $u\cdot [\AG]=2$,
the vector $u$ is the class of an effective divisor $D$ on $\XG$.
Let $D\sb 0$ be the union of irreducible components of $D$
whose image by $\phiG$ are of dimension $1$.
Since  $u-[D\sb 0]$ is a linear combination of the classes $[\Gamma\sb P]$
with non-negative integer coefficients,
we have
\begin{equation*}\label{eq:conicvec3}
[D\sb 0] \,\bmod \SG\sp 0 =(A, 0)\qquad\textrm{in $\CCG\lift$}.
\end{equation*}
Because $D\sb 0 \AG=2$,
the plane curve $\phiG(D\sb 0)$ with the reduced structure 
is either a line or a conic.
Suppose that $\phiG(D\sb 0)$ is a line $L$.
If $L$ is not splitting in $\XG$,
then the morphism $\shortmap{\phiG|\sb{D\sb 0}}{D\sb 0}{ L}$ is  of degree $2$,
while 
if $L$ is  splitting,
then $D\sb 0$ is  $2\HPT{L}$.
In either case,
$D\sb 0$ is the proper transform of $L$
and hence $[D\sb 0]$ is  contained in $\SG\sp0$.
This is absurd because $A\ne 0$.
Therefore $\phiG(D\sb 0)$ is a conic $Q$.
Since $\shortmap{\phiG|\sb{D\sb 0}}{D\sb 0}{Q}$ is of degree $1$,
the  conic $Q$ is splitting,  and $D\sb 0=\HPT{Q}$ holds.
From~\eqref{eq:conicvec3},
we have
$A=\wG(Q)$.
If $Q$ is
 a union of two lines $L\sb 1$ and $L\sb 2$,
then both  $L\sb 1$ and $L\sb 2$ are splitting and 
$A=\wG(L\sb 1)+\wG(L\sb 2)$ holds from~\eqref{eq:sum},
which contradicts the irreducibility of the word $A$ in $\CCG$.
(See Remark~\ref{rem:intersection_lines}.)
Therefore $Q$ is a smooth conic.
Because
$\wG (Q)=Q\cap\ZdG$ by~\eqref{eq:wm},
we obtain $A=Q\cap\ZdG$.
\end{proof}
\begin{remark}\label{rem:intersection_line_and_conic}
Let $L$ be a splitting line,
and $Q$ a splitting smooth conic.
If $L$ intersects $Q$ transversely,
then $\wG (L)\cap \wG (Q)$ consists of the two intersection points of $L$ and $Q$,
and $\wG(L\cup Q)= \wG (L) +\wG(Q)$ is of weight $9$.
If $L$ is tangent to  $Q$,
then $\wG (L)\cap \wG (Q)$ is empty,
and $\wG(L\cup Q)= \wG (L) +\wG(Q)$ is of weight $13$.
\end{remark}
\begin{remark}\label{rem:conicintersection}
Let $Q\sb 1$ and $Q\sb 2$ be  distinct splitting smooth conics.
Let us investigate the intersection of $Q\sb 1$ and $Q\sb 2$.
Because
$$
\card{\wG(Q\sb 1\cup Q\sb 2)}= \card{\wG(Q\sb 1) + \wG(Q\sb 2)}=16-2\card{\wG(Q\sb 1) \cap \wG(Q\sb 2)}
$$
is in $\{0,5,8,9,12,13,16,21\}$
by Theorem~\ref{thm:wt},
$\card{\wG(Q\sb 1) \cap \wG(Q\sb 2)}$ is $4$, $2$ or $0$.

Suppose that  $\card{\wG(Q\sb 1) \cap \wG(Q\sb 2)}=4$.
Then
$Q\sb 1$ and $Q\sb 2$ intersect transversely.
Let $G\sb {Q\sb 1}$ and $G\sb{Q\sb 2}$ be homogeneous polynomials of degree $2$
defining $Q\sb 1$ and $Q\sb 2$, respectively.
Since $Q\sb 1\cup Q\sb 2$ is a  splitting curve with only ordinary nodes, 
Proposition~\ref{prop:convGC} implies that 
there exists a homogeneous polynomial $G\sb{Q\sb 3}$ of degree $2$
such that $G\sb{Q\sb 1} G\sb{Q\sb 2} G\sb{Q\sb 3}$
is a member of  $k\sptimes G +\Vts$.
Then the conic  $Q\sb 3$ defined by $G\sb{Q\sb 3}=0$ is splitting in $\XG$,
and
$\wG (Q\sb 3)=\wG(Q\sb 1) +\wG (Q\sb 2)$ holds.

Suppose that  $\card{\wG(Q\sb 1) \cap \wG(Q\sb 2)}=2$.
By Proposition~\ref{prop:mildsing},
we have the following two possibilities of intersection of $Q\sb 1$ and $Q\sb 2$;
\begin{itemize}
\item
transverse at two points, and with multiplicity $2$ at one point, or 
\item
transverse at one point, and with multiplicity $3$ at one point.
\end{itemize}

Suppose that  $\card{\wG(Q\sb 1) \cap \wG(Q\sb 2)}=0$.
Then
$Q\sb 1$ and $Q\sb 2$ intersect
either
 with multiplicity $2$ at two points, or 
with multiplicity $4$ at one  point.
\end{remark}
\begin{corollary}\label{cor:collinear}
A word $A\in \CCG$ of weight $8$ or $9$ is irreducible in $\CCG$ if and only if no three points of $A$
are collinear.
\end{corollary}
\begin{proof}
Suppose that $A$ is reducible in $\CCG$.
Then $A$ is written as $A\sb 1+ A\sb 2$,
where 
$A\sb 1$ and $A\sb 2$ are words of $\CCG$ such that 
$(\card{A}, \card{A\sb 1}, \card{A\sb 2})$ 
is either $(8, 5, 5)$ or $(9,5,8)$.
Note that $A\cap A\sb 1=A\sb 1\sm (A\sb 1\cap A\sb 2)$
is of weight $\ge 3$,
because $\card{A\sb 1\cap A\sb 2}= (\card{A\sb 1}+\card{A\sb 2}-\card{A})/2$ is $\le 2$.
Since the points of $A\sb 1$ are collinear by Proposition~\ref{prop:line},
three points  of $A$ are collinear.
Suppose that three points of $A$ are on a line $L$.
By Proposition~\ref{prop:convnum},
the line $L$ is splitting in $\XG$.
We put
$A\sprime :=A+\wG (L) \in \CCG$.
The weight
$$
\card{A\sprime}=\card{A}+5- 2 \card{A\cap \wG(L)}
$$
of $A\sprime$ is among the set $\{0,5,8,9,12,13,16,21\}$
by Theorem~\ref{thm:wt}.
Because $\wG(L)=L\cap\ZdG$ and $A\st\ZdG$,
we have $A\cap \wG (L)= A\cap L$ and hence 
$\card{A\cap \wG(L)}$ is $\ge 3$.
Therefore
the triple $(\card{A}, \card{A\cap \wG (L)}, \card{A\sprime})$ is either $(8, 4, 5)$
or $(9, 3, 8)$.
In either case, $A= A\sprime +\wG (L)$ is reducible in $\CCG$.
\end{proof}
\begin{definition}
A pencil $\EEE=\{E\sb t\}$ 
of cubic curves $E\sb t\st\Pt$ is called {\it regular}
if the base locus $\Bs (\EEE)$ of $\EEE$ consists of distinct $9$ points
and every singular member of $\EEE$ is an irreducible nodal curve.
\end{definition}
Note that the  general member of a regular
pencil $\EEE$ of cubic curves is smooth.
Indeed,
the general member of $\EEE$ is reduced because $\card{\Bs(\EEE)} =9$.
If the   general member of $\EEE$ is singular,
then it  must have an ordinary cusp
(\cite{Tate, Shimada92}),
and hence any singular member cannot be an irreducible nodal curve.
\begin{lemma}\label{lem:reg}
Let $\EEE$ be a regular pencil of cubic curves.
\par
{\rm (1)} The pencil $\EEE$ coincides with $|\III\sb{\Bs (\EEE)} (3)|$.
\par
{\rm (2)} There are no three collinear points in $\Bs (\EEE)$.
\end{lemma}
\begin{proof}
In order to prove  (1),
it is enough to show that $\dim |\III\sb{\Bs (\EEE)} (3)|\le 1$.
If $\dim |\III\sb{\Bs (\EEE)} (3)|> 1$,
then there would be  eight points in $\Bs(\EEE)$ on a conic,
or five points in $\Bs(\EEE)$ on a line.
(See for example~\cite[p.715]{GH}.)
In either case,
we get a contradiction to B\'ezout's theorem.
Suppose that there exists a subset of $\Bs(\EEE)$ of weight $3$  that is on a line $L$.
We put $B\sprime := \Bs(\EEE)\cap L$,
and let $\III\sb{B\sprime\st L}\subset \OOO\sb L$ be the ideal sheaf of $B\sprime$ on $L$.
From the exact sequence
$$
\Hz {\III\sb{\Bs(\EEE)\sm B\sprime} (2)}\;\to\;
\Hz {\III\sb{\Bs(\EEE)} (3)}\;\to\;
H\sp 0 (L, \III\sb{B\sprime\st L} (3)),
$$
we see that a union of $L$ and a conic is a member of  $\EEE=|\III\sb{\Bs (\EEE)} (3)|$,
which contradicts the regularity of $\EEE$.
\end{proof}
\begin{definition}
A pencil $\EEE$ 
of cubic curves  is called {\it splitting in $\XG$}
if every member of $\EEE$ is reduced and splitting in $\XG$.
\end{definition}
\begin{proposition}\label{prop:cubic}
The correspondence
$\EEE\mapsto \Bs(\EEE)$ gives a bijection from the set of 
regular pencils of cubic curves 
splitting in $\XG$ to the set of  irreducible words $A\in \CCG$ of weight  $9$.
The inverse map is given by $A\mapsto |\III\sb{A} (3)|$.
\end{proposition}
\begin{proof}
Let $\EEE$ be a regular pencil of cubic curves splitting in $\XG$,
and let  $E$ and $E\sprime$ be  members of $\EEE$
that span  $\EEE$.
Each of $E$ and $E\sprime$ is a  smooth or irreducible nodal cubic curve
splitting in $\XG$.
Let $E\sp o$ and $E\sp{\prime o}$ be the smooth parts of $E$ and $E\sprime$,
respectively.
Then we have
\begin{equation}\label{eq:cubic1}
\wG (E) =E\sp o \cap \ZdG
\quand
\wG (E\sprime) =E\sp {\prime o} \cap \ZdG
\end{equation}
by~\eqref{eq:wm},
and
\begin{equation}\label{eq:cubic2}
\card{\wG (E)} =\card{\wG (E\sprime)} =9
\end{equation}
by Proposition~\ref{prop:num}.
On the other hand,
the base locus $\Bs (\EEE)$
of $\EEE$ is equal to $E\sp o \cap E\sp{\prime o}$,
and is contained in the set of ordinary nodes
of the reducible splitting curve $E\cup E\sprime$.
Hence 
\begin{equation}\label{eq:cubic3}
\Bs (\EEE) = E\sp o \cap E\sp{\prime o}\subset \ZdG
\end{equation}
holds by Corollary~\ref{cor:nodeandcusp}.
Comparing~\eqref{eq:cubic1},~\eqref{eq:cubic2} and~\eqref{eq:cubic3},
we obtain
$$
\wG(E)=\wG(E\sprime)=\Bs (\EEE).
$$
In particular, $\Bs(\EEE)$ is a word in $\CCG$.
From Lemma~\ref{lem:reg} and Corollary~\ref{cor:collinear},
the word $\Bs(\EEE)$ is irreducible in $\CCG$.
\par
Suppose that an irreducible word $A$ of $\CCG$ with weight $9$ is given.
A splitting regular pencil $\EEE$ with $\Bs (\EEE)=A$ is, if exists,
equal to $|\III\sb A (3)|$ 
by Lemma~\ref{lem:reg},
and hence is unique.
Let us prove the existence of such a pencil $\EEE$.
Since $(A, 1)\in\CCG\lift$,
we have a vector
\begin{equation*}\label{eq:cubicvec}
u:=\frac{1}{2} (-\sum\sb{P\in A} [\Gamma\sb P]\;+\; 3\,[\AG])
\end{equation*}
in $\SG$.
Because $u^2=0$ and $u\cdot [\AG]=3$,
the vector $u$ is the class of an effective divisor $D$ on $\XG$.
Let $D\sb 0$ be the union of irreducible components of $D$
whose image by $\phiG$ are of dimension $1$.
Because $D-D\sb 0$ is a sum of the curves $\Gamma\sb P$
with non-negative integer coefficients, we have 
\begin{equation*}\label{eq:cubicvec2}
[D\sb 0] \,\bmod \SG\sp 0 =(A, 1)\qquad\textrm{in $\CCG\lift$}.
\end{equation*}
Because $D\sb 0 \AG=3$,
there are three possibilities;
\begin{itemize}
\item
there exists a splitting line $L$
such that $D\sb 0 =3 \HPT{L}$, 
\item
there exist distinct lines $L$ and $L\sprime$ such that
$L$ is splitting and that
$D\sb 0$ is the union of $\HPT{L}$ and the proper transform of $L\sprime$, or
\item
there exists a reduced cubic curve $E$ splitting in $\XG$ such that $D\sb 0=\HPT{E}$.
\end{itemize}
In the first or the second case, we have
$(A, 1)=[\HPT{L}] \,\bmod \SG\sp 0$, 
and hence $\card{A}=\card{\wG(L)}=5$,
which contradicts the assumption.
Therefore
there exists a 
reduced splitting cubic curve $E$ such that $D\sb 0=\HPT{E}$.
In particular,
we have $A=\wG(E)$.
If $E$ were reducible,
then
the word $A$ would be  also reducible in $\CCG$ by
Remarks~\ref{rem:intersection_threelines}~and~\ref{rem:intersection_line_and_conic}.
Hence $E$ is irreducible.
If $E$ had an ordinary  cusp,
then $A=\wG(E)$ would be of weight $13$ by Proposition~\ref{prop:num}.
Therefore $E$ is a smooth or irreducible nodal cubic curve.
Let $G\sb E$ be a homogeneous polynomial of degree $3$
such that $E$ is defined by $G\sb E=0$.
By Proposition~\ref{prop:convGC},
there exists another homogeneous cubic polynomial $G\sb{E\sprime}$
such that $G\sb E G\sb{E\sprime}\in k\sptimes G+\Vts$.
For $t\in k$,
we put
$$
G\sb{E\sb t}:= G\sb{E\sprime} + t G\sb E.
$$
Then we have
$$
G\sb E G\sb{E\sb t}\in k\sptimes G+ \Vts
$$
for any $t\in k$.
Let $E\sb t$ denote the cubic curve defined by $G\sb{E\sb t}=0$,
and let $\EEE$ be the pencil $\shortset{ E\sb t}{t\in k\cup \{\infty\}}$.
By Proposition~\ref{prop:GC},
every  member $E\sb t$ is a reduced curve
with only ordinary nodes as its singularities,
and is splitting in $\XG$.
Moreover, the cubic curves 
$E$ and $E\sb{t}$ intersect transversely and
$$
\wG (E)=\wG (E\sb t) =E\cap E\sb t = \Bs (\EEE).
$$
Hence $\EEE$ is a pencil splitting in $\XG$ such that $\Bs (\EEE)=A$.
If a member $E\sb{t\sb 0}$ of $\EEE$ were reducible,
then
the word $A=\wG(E\sb{t\sb 0})$ would  also  be reducible in $\CCG$.
Hence $\EEE$ is regular.
\end{proof}
\begin{corollary}
The word $\Bs (\EEE)$ of $\CCG$ corresponding to 
a regular splitting pencil $\EEE$ of cubic curves 
is equal to $\wG (E)$,
where $E$ is an arbitrary  member of $\EEE$.
\end{corollary}
\begin{corollary}\label{cor:EE}
Let $A\in \CCG$ be an irreducible word of weight $9$.
If the $2$-dimensional vector space $\Hz {\III\sb A (3)}$
is generated by $G\sb E$ and $G\sb{E\sprime}$,
then the homogeneous polynomial $G\sb E G\sb{E\sprime} $ of degree $6$ is contained in $k\sptimes G + \Vts$.
\end{corollary}
\begin{remark}\label{rem:12}
It is known that a regular pencil $\EEE$ of cubic curves has
exactly $12$
singular members $\{E\sb1, \dots, E\sb{12}\}$.
Suppose that the  regular pencil $\EEE$ is splitting in $\XG$.
The ordinary node $P\sb i$ of a singular member $E\sb i$ is a point of  $\ZdG$ by Corollary~\ref{cor:nodeandcusp}.
By assigning $P\sb i$ to the singular member $E\sb i$,
we obtain a bijection
$$
\{E\sb1, \dots, E\sb{12}\}\;\cong\; \ZdG\sm \Bs(\EEE).
$$
\end{remark}
\begin{remark}
The decomposition of a reducible word $A\in \CCG$ of weight $9$
into a sum of irreducible words is {\it not} unique.
For example,
let $G\sb 1$ and  $G\sb 1\sprime$ be  general homogeneous polynomials of degree $1$,
and 
let $G\sb 2$ and  $G\sb 2\sprime$ be  general homogeneous polynomials of degree $2$.
Then $G:=G\sb 1 G\sb 1\sprime G\sb 2 G\sb 2\sprime$ is  contained  in $\Uts$.
(See Example~\ref{example:nodal2211}.)
The lines $L:=\{G\sb 1=0\}$, $L\sprime := \{ G\sb 1\sprime=0\}$
and the smooth conics $Q:=\{G\sb 2=0\}$, $Q\sprime := \{ G\sb 2\sprime=0\}$
are splitting in $\XG$ by Proposition~\ref{prop:GC}.
We have two decompositions of the word
$$
\wG(L)+\wG(Q)= \wG(L\sprime)+\wG(Q\sprime)
$$
of weight $9$,
which is equal to $\wG (E)$,
where $E$ is an arbitrary member of the splitting (non-regular)
pencil of cubic curves spanned by $L\cup Q$ and $L\sprime\cup Q\sprime$.
\end{remark}
\begin{remark}
Let $\EEE$ be a regular splitting pencil of cubic curves.
\par
Let $L$ be a splitting line.
Because 
$$
\card{\Bs(\EEE)+\wG (L)}=14 - 2 \card{\Bs(\EEE)\cap\wG (L)},
$$
the weight of 
$\Bs(\EEE)\cap\wG (L)$ is either $1$ or $3$.
By Corollary~\ref{cor:collinear},
$\card{\Bs(\EEE)\cap\wG (L)}$ cannot be $3$.
Let $E\sb t$ be the general member of $\EEE$.
Suppose that $E\sb t$ intersects $L$ transversely at a point $P$.
Then $P$ is an ordinary node of the reducible splitting curve $E\sb t \cup L$,
and hence $P\in \ZdG$ by Corollary~\ref{cor:nodeandcusp}.
In particular, $P$ is contained in $ \Bs(\EEE)\cap\wG (L)$.
Therefore the restriction $\EEE|L$ of $\EEE$ to $L$
consists of one fixed point
and a moving non-reduced point of multiplicity $2$.
\par
Let $Q$ be a smooth splitting conic.
Then 
$\card{\Bs(\EEE)\cap\wG (Q)}$ is either $2$ or $4$ or $6$.
Suppose that $\card{\Bs(\EEE)\cap\wG (Q)}=6$,
and let $P$ be a point of $\wG (Q) \sm (\Bs(\EEE)\cap\wG (Q))$.
There exists a member $E\sb P$ of $\EEE$ that has an ordinary node at $P$
by Remark~\ref{rem:12}.
Then $Q$ must be contained in $E\sb P$,
which contradicts the regularity of $\EEE$.
Hence $\card{\Bs(\EEE)\cap\wG (Q)}$ is  $2$ or $4$.
When $\card{\Bs(\EEE)\cap\wG (Q)}=2$ (resp.  $4$), 
the restriction $\EEE|Q$ of $\EEE$ to $Q$
consists of two  (resp. four) fixed points
and  moving non-reduced points of total multiplicity $4$ (resp. $2$).
 \end{remark}
\begin{remark}
Let $A\in \CCG$ be a word of weight $13$.
Then one of the following holds:
\par
(i) There are three splitting lines $L\sb 1, L\sb 2, L\sb 3$
meeting at a point such that
$A=\wG(L\sb 1)+\wG(L\sb 2)+\wG(L\sb 3)$.
\par
(ii) 
There are a splitting line $L$ and a splitting smooth conic $Q$
such that $L$ is tangent to $Q$ and that
$A=\wG(L)+\wG(Q)$.
\par
(iii)
There exists a cuspidal cubic curve $C$ splitting in $\XG$ such that $A=\wG (C)$.
\par
We put
$G\sb Q:= X\sb 0^2 +X\sb 1 X\sb 2$,
and let $G\sb 4$ be a general homogeneous polynomial of degree $4$.
Then $G\sb Q G\sb 4$ is a polynomial  in  $\Uts$,
and the smooth conic $Q$ defined by $G\sb Q=0$ is splitting in $X\sb{G\sb Q G\sb 4}$.
Let $C$ be the cubic curve defined by
$\der{G\sb 4}{ X\sb 0}=0$.
It is easy to see that $C$ has one ordinary cusp as its only singularities,
and is splitting in $X\sb{G\sb Q G\sb 4}$.
Moreover, the word $w\sb{G\sb Q G\sb 4} (C)$ coincides with $Z(dG\sb Q G\sb 4)\sm w\sb{G\sb Q G\sb 4} (Q)$.
\end{remark}
Since $\CCG$ is generated by $\ZdG\in \CCG$
and irreducible  words of weight $5$, $8$ and $9$,
we obtain the following:
\begin{corollary}\label{cor:generator}
The lattice $\SG$ is generated by the following vectors;
\begin{itemize}
\item $[\AG]$ and $[\Gamma\sb P]$ $(P\in \ZdG)$, 
\item $[\HPT{C}]$, where $C$ is the  general member of $|\III\sb{\ZdG} (5)|$, 
\item $[\HPT{L}]$, where $L$ runs through the set of  splitting lines, 
\item $[\HPT{Q}]$, where $Q$ runs through the set of splitting smooth conics,
\item $[\HPT{E}]$, where $E$ runs through the set of  members of  regular splitting pencils of cubic curves.
\end{itemize}
\end{corollary}
Main Theorem in Introduction has now been proved by
Propositions~\ref{prop:stdsplit}, \ref{prop:line}, \ref{prop:conic}, \ref{prop:cubic} and Corollary~\ref{cor:generator}. 
\subsection{The list.}
Using Theorem~\ref{thm:wt} and Algorithm~\ref{algo:listup},
we  make the complete list of geometrically realizable 
classes of codes.
In the list below,
the following data are recorded.
\begin{itemize}
\item $\sigma$: The Artin invariant $11-\dim \tC$ of the  corresponding supersingular $K3$ surfaces.
For each $\sigma$, 
the number $r(\sigma)$ of geometrically realizable classes with Artin invariant $\sigma$ is also given.
\item {\tt std}: A standard basis of the $\St$-equivalence class $[\tC]$. 
(See Definition~\ref{definition:stdbasis}.)
A word is  expressed by a bit vector, and 
a bit vector $[\alpha\sb 0, \dots, \alpha\sb {20}]$ is expressed by the integer
$2^{20} \alpha\sb 0+ \cdots + 2\alpha\sb{19} +\alpha\sb{20}$.
Since $[1, \dots, 1]=2^{21}-1$ corresponding to the word $\tZ$ is always in  standard bases by definition,
it is omitted.
\item {\tt l}: The number of words of weight $5$; that is,
the number of splitting lines.
\item {\tt q}: The number of irreducible words of weight $8$; that is,
the number of splitting smooth conics.
\item {\tt e}: The number of irreducible words of weight $9$; that is,
the number of splitting regular pencils of cubic curves.
\end{itemize}
There are several pairs of  classes of codes
with identical   $(\sigma, {\tt l}, {\tt q}, {\tt e})$.
(For example, the classes No.134\,-\,No.136.
See Examples~\ref{example:792} and~\ref{example:pascal}.)
By trial and error,
we have found that the following added data are sufficient to distinguish all the 
geometrically realizable  classes of codes.
\begin{itemize}
\item {\tt tl}: The number of triples $\sset{L\sb 1, L\sb 2, L\sb 3}$ of splitting lines such that 
$L\sb 1\cap L\sb 2\cap L\sb 3$ consists of one point;
that is, the number of triples $\sset{A\sb 1, A\sb 2, A\sb 3}$
of distinct words  of weight $5$ 
satisfying   $\card{A\sb 1\cap A\sb 2\cap A\sb 3}=1$.
\item {\tt lq}: The number of pairs $(L, Q)$ 
of a splitting line $L$ and a splitting smooth conic $Q$ such that $L$ is tangent to $Q$;
that is, the number of pairs $(A, B)$ of words such that
$\card{A}=5$, $\card{B}=8$, $B$ is irreducible, and $A\cap B=\emptyset$.
\item {\tt qq}: The number of pairs $\sset{Q, Q\sprime}$
of  splitting smooth conics such that
there exist exactly two  points of $Q\cap Q\sprime$ 
at which $Q$ and $Q\sprime$ intersect 
with  odd intersection multiplicity;
that is, the number of pairs $\sset{A, A\sprime}$
of irreducible words of weight $8$ such that $\card{A\cap A\sprime}=2$.
See Figure~\ref{fig:qq}.
\begin{figure}
\includegraphics[height=2cm]{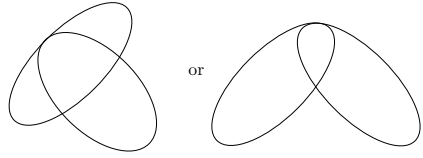}
\caption{The configurations of smooth conics for {\tt qq}}\label{fig:qq}
\end{figure}
\item {\tt tq1}: 
The number of triples $\sset{Q\sb 1, Q\sb 2, Q\sb 3}$ of smooth splitting conics with 
the configuration as in Figure~\ref{fig:qqq};
that is,
the number of triples $\sset{A\sb 1, A\sb 2, A\sb 3}$ of
irreducible words of weight $8$
such that $\card{A\sb i\cap A\sb j}=4$ for each $i\ne j$ and $\card{A\sb 1\cap A\sb 2\cap A\sb 3}=3$.
\begin{figure}
\includegraphics[height=2.4cm]{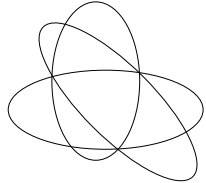}
\caption{The configuration of smooth conics for {\tt tq1}}
\label{fig:qqq}
\end{figure}
\item {\tt tq2}: 
The number of triples $\sset{Q\sb 1, Q\sb 2, Q\sb 3}$ of smooth splitting conics 
such that, for each $i, j$ with $i\ne j$,
there exist exactly two  points of $Q\sb i\cap Q\sb j$ 
at which $Q\sb i$ and $Q\sb j$ intersect 
with  odd intersection multiplicity;
that is,
the number of triples $\sset{A\sb 1, A\sb 2, A\sb 3}$ of
irreducible words of weight $8$
such that $\card{A\sb i\cap A\sb j}=2$ for $i\ne j$.
See Figure~\ref{fig:qqqsprime}.
\begin{figure}
\includegraphics[height=2.4cm]{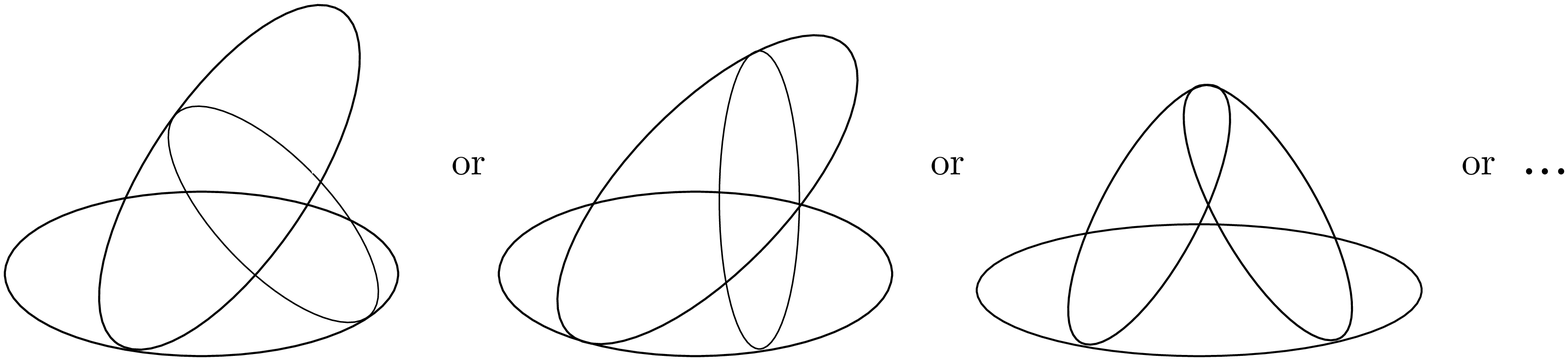}
\caption{The configurations of smooth conics for {\tt tq2}}
\label{fig:qqqsprime}
\end{figure}
\end{itemize}
\newcommand{\cparbox}[2]{\parbox[t]{#1}{\hfil#2\hfil}}
\newcommand{\codebox}[9]{%
\framebox {%
		\hbox {%
			\cparbox{.7cm}{#1}%
			\vrule
			\cparbox{.6cm}{#2}% sigma
   \vrule \hskip 3pt
			\parbox[t]{5.35cm}{#3\hfil}
			\vrule
			\cparbox{.6cm}{#4}% l
			\cparbox{.6cm}{#5}% q
			\cparbox{.6cm}{#6}% e
			\vrule
			\cparbox{.6cm}{#7}% tl
   \vrule
			\cparbox{.6cm}{#8}% lq
			\vrule
			\cparbox{2.5cm}{#9}% qq tq1 tq2
		}%
	}%
\vskip -1.4pt
}
\par\medskip\noindent
\begin{center}
{\bf The complete list of geometrically  realizable  classes of codes}
\end{center}
\par\medskip\noindent
%
%%%%%%%%%%%%%%%%%%%%%%%%%%
% Here starts the list.
%%%%%%%%%%%%%%%%%%%%%%%%%%
%
{\small
\codebox{No.}{$\sigma$}{\hfil\hbox{\tt std}}{\hbox{\tt l}}%
{\hbox{\tt q}} {\hbox{\tt e}} {\hbox{\tt tl}}  {\hbox{\tt lq}}%
  {\hbox{\tt qq} \hfil\hbox{\tt tq1}\hfil \hbox{\tt tq2}}\par\noindent
\smallskip\par\noindent
$\sigma=10. \qquad r(10)=1$.\smallskip\par\noindent
\codebox{0}{10}{}{0}{0}{0}{0}{0}{0,\hfil 0, \hfil 0}\par\noindent
\smallskip\par\noindent
$\sigma=9. \qquad r(9)=3$.\smallskip\par\noindent
\codebox{1}{9}{31}{1}{0}{0}{0}{0}{0,\hfil 0, \hfil 0}\par\noindent
\codebox{2}{9}{255}{0}{1}{0}{0}{0}{0,\hfil 0, \hfil 0}\par\noindent
\codebox{3}{9}{511}{0}{0}{1}{0}{0}{0,\hfil 0, \hfil 0}\par\noindent
\smallskip\par\noindent
$\sigma=8. \qquad r(8)=8$.\smallskip\par\noindent
\codebox{4}{8}{31, 481}{2}{0}{0}{0}{0}{0,\hfil 0, \hfil 0}\par\noindent
\codebox{5}{8}{31, 8160}{1}{2}{0}{0}{2}{0,\hfil 0, \hfil 0}\par\noindent
\codebox{6}{8}{31, 2019}{1}{1}{0}{0}{0}{0,\hfil 0, \hfil 0}\par\noindent
\codebox{7}{8}{31, 8161}{1}{0}{2}{0}{0}{0,\hfil 0, \hfil 0}\par\noindent
\codebox{8}{8}{255, 3855}{0}{3}{0}{0}{0}{0,\hfil 0, \hfil 0}\par\noindent
\codebox{9}{8}{255, 16131}{0}{2}{1}{0}{0}{1,\hfil 0, \hfil 0}\par\noindent
\codebox{10}{8}{255, 7951}{0}{1}{2}{0}{0}{0,\hfil 0, \hfil 0}\par\noindent
\codebox{11}{8}{511, 32263}{0}{0}{3}{0}{0}{0,\hfil 0, \hfil 0}\par\noindent
\smallskip\par\noindent
$\sigma=7. \qquad r(7)=21$.\smallskip\par\noindent
\codebox{12}{7}{31, 8160, 481}{3}{1}{0}{1}{3}{0,\hfil 0, \hfil 0}\par\noindent
\codebox{13}{7}{31, 2019, 2301}{3}{0}{0}{0}{0}{0,\hfil 0, \hfil 0}\par\noindent
\codebox{14}{7}{31, 8160, 516193}{2}{2}{0}{0}{2}{0,\hfil 0, \hfil 0}\par\noindent
\codebox{15}{7}{31, 2019, 6244}{2}{2}{0}{0}{0}{0,\hfil 0, \hfil 0}\par\noindent
\codebox{16}{7}{31, 8161, 253987}{2}{1}{1}{0}{0}{0,\hfil 0, \hfil 0}\par\noindent
\codebox{17}{7}{31, 8160, 123360}{1}{6}{0}{0}{6}{0,\hfil 0, \hfil 0}\par\noindent
\codebox{18}{7}{31, 8160, 25059}{1}{4}{0}{0}{2}{2,\hfil 0, \hfil 0}\par\noindent
\codebox{19}{7}{31, 2019, 63533}{1}{3}{0}{0}{0}{3,\hfil 0, \hfil 1}\par\noindent
\codebox{20}{7}{31, 2019, 14565}{1}{3}{0}{0}{0}{0,\hfil 0, \hfil 0}\par\noindent
\codebox{21}{7}{31, 8160, 123361}{1}{2}{4}{0}{2}{0,\hfil 0, \hfil 0}\par\noindent
\codebox{22}{7}{31, 8161, 25062}{1}{2}{2}{0}{0}{1,\hfil 0, \hfil 0}\par\noindent
\codebox{23}{7}{31, 8161, 254178}{1}{1}{4}{0}{0}{0,\hfil 0, \hfil 0}\par\noindent
\codebox{24}{7}{255, 3855, 13107}{0}{7}{0}{0}{0}{0,\hfil 0, \hfil 0}\par\noindent
\codebox{25}{7}{255, 3855, 28951}{0}{6}{1}{0}{0}{3,\hfil 4, \hfil 0}\par\noindent
\codebox{26}{7}{255, 3855, 62211}{0}{5}{2}{0}{0}{4,\hfil 0, \hfil 0}\par\noindent
\codebox{27}{7}{255, 3855, 127249}{0}{4}{3}{0}{0}{3,\hfil 0, \hfil 0}\par\noindent
\codebox{28}{7}{255, 16131, 115471}{0}{3}{4}{0}{0}{3,\hfil 0, \hfil 1}\par\noindent
\codebox{29}{7}{255, 3855, 29491}{0}{3}{4}{0}{0}{0,\hfil 0, \hfil 0}\par\noindent
\codebox{30}{7}{255, 16131, 50973}{0}{2}{5}{0}{0}{1,\hfil 0, \hfil 0}\par\noindent
\codebox{31}{7}{255, 7951, 123187}{0}{1}{6}{0}{0}{0,\hfil 0, \hfil 0}\par\noindent
\codebox{32}{7}{511, 32263, 233016}{0}{0}{7}{0}{0}{0,\hfil 0, \hfil 0}\par\noindent
\smallskip\par\noindent
$\sigma=6. \qquad r(6)=43$.\smallskip\par\noindent
\codebox{33}{6}{31, 8160, 123360, 1966081}{5}{0}{0}{10}{0}{0,\hfil 0, \hfil 0}\par\noindent
\codebox{34}{6}{31, 8160, 25059, 28385}{4}{1}{0}{1}{3}{0,\hfil 0, \hfil 0}\par\noindent
\codebox{35}{6}{31, 2019, 6244, 8637}{4}{1}{0}{0}{0}{0,\hfil 0, \hfil 0}\par\noindent
\codebox{36}{6}{31, 8160, 25059, 105991}{3}{5}{0}{1}{7}{0,\hfil 0, \hfil 0}\par\noindent
\codebox{37}{6}{31, 8160, 25059, 26215}{3}{5}{0}{1}{3}{4,\hfil 0, \hfil 0}\par\noindent
\codebox{38}{6}{31, 8161, 253987, 319591}{3}{3}{1}{0}{0}{0,\hfil 1, \hfil 0}\par\noindent
\codebox{39}{6}{31, 8160, 25059, 238049}{3}{3}{0}{1}{3}{0,\hfil 0, \hfil 0}\par\noindent
\codebox{40}{6}{31, 8160, 25059, 42497}{3}{3}{0}{0}{2}{1,\hfil 0, \hfil 0}\par\noindent
\codebox{41}{6}{31, 8160, 516193, 582560}{2}{6}{0}{0}{6}{0,\hfil 0, \hfil 0}\par\noindent
\codebox{42}{6}{31, 8160, 25059, 100324}{2}{6}{0}{0}{4}{6,\hfil 0, \hfil 0}\par\noindent
\codebox{43}{6}{31, 8160, 25059, 44583}{2}{6}{0}{0}{2}{6,\hfil 2, \hfil 2}\par\noindent
\codebox{44}{6}{31, 2019, 63533, 68551}{2}{6}{0}{0}{0}{12,\hfil 0, \hfil 8}\par\noindent
\codebox{45}{6}{31, 2019, 6244, 27049}{2}{6}{0}{0}{0}{0,\hfil 0, \hfil 0}\par\noindent
\codebox{46}{6}{31, 8160, 25059, 492257}{2}{4}{2}{0}{2}{2,\hfil 0, \hfil 0}\par\noindent
\codebox{47}{6}{31, 8161, 253987, 271302}{2}{4}{2}{0}{0}{5,\hfil 0, \hfil 2}\par\noindent
\codebox{48}{6}{31, 8161, 253987, 288708}{2}{4}{2}{0}{0}{2,\hfil 0, \hfil 0}\par\noindent
\codebox{49}{6}{31, 8160, 123360, 419424}{1}{14}{0}{0}{14}{0,\hfil 0, \hfil 0}\par\noindent
\codebox{50}{6}{31, 8160, 25059, 241184}{1}{10}{0}{0}{6}{12,\hfil 16, \hfil 0}\par\noindent
\codebox{51}{6}{31, 8160, 25059, 124512}{1}{10}{0}{0}{6}{12,\hfil 0, \hfil 0}\par\noindent
\codebox{52}{6}{31, 8160, 25059, 492069}{1}{8}{0}{0}{2}{12,\hfil 4, \hfil 4}\par\noindent
\codebox{53}{6}{31, 8160, 25059, 42605}{1}{8}{0}{0}{2}{6,\hfil 0, \hfil 0}\par\noindent
\codebox{54}{6}{31, 8160, 123360, 419425}{1}{6}{8}{0}{6}{0,\hfil 0, \hfil 0}\par\noindent
\codebox{55}{6}{31, 8160, 25059, 99948}{1}{6}{4}{0}{2}{8,\hfil 0, \hfil 4}\par\noindent
\codebox{56}{6}{31, 8160, 25059, 238119}{1}{6}{4}{0}{2}{8,\hfil 0, \hfil 0}\par\noindent
\codebox{57}{6}{31, 8161, 25062, 99051}{1}{6}{2}{0}{0}{9,\hfil 0, \hfil 4}\par\noindent
\codebox{58}{6}{31, 8161, 25062, 42602}{1}{6}{2}{0}{0}{3,\hfil 4, \hfil 0}\par\noindent
\codebox{59}{6}{31, 8160, 25059, 239201}{1}{4}{8}{0}{2}{2,\hfil 0, \hfil 0}\par\noindent
\codebox{60}{6}{31, 8161, 25062, 229998}{1}{4}{6}{0}{0}{6,\hfil 0, \hfil 4}\par\noindent
\codebox{61}{6}{31, 8161, 25062, 501288}{1}{4}{6}{0}{0}{3,\hfil 0, \hfil 0}\par\noindent
\codebox{62}{6}{255, 3855, 13107, 21845}{0}{15}{0}{0}{0}{0,\hfil 0, \hfil 0}\par\noindent
\codebox{63}{6}{255, 3855, 28951, 46881}{0}{13}{2}{0}{0}{12,\hfil 32, \hfil 0}\par\noindent
\codebox{64}{6}{255, 3855, 28951, 492145}{0}{11}{4}{0}{0}{16,\hfil 16, \hfil 0}\par\noindent
\codebox{65}{6}{255, 3855, 62211, 208947}{0}{9}{6}{0}{0}{18,\hfil 0, \hfil 6}\par\noindent
\codebox{66}{6}{255, 3855, 28951, 233577}{0}{9}{6}{0}{0}{15,\hfil 8, \hfil 3}\par\noindent
\codebox{67}{6}{255, 3855, 13107, 116021}{0}{9}{6}{0}{0}{12,\hfil 0, \hfil 0}\par\noindent
\codebox{68}{6}{255, 3855, 127249, 405606}{0}{7}{8}{0}{0}{12,\hfil 0, \hfil 4}\par\noindent
\codebox{69}{6}{255, 3855, 28951, 111147}{0}{7}{8}{0}{0}{9,\hfil 4, \hfil 3}\par\noindent
\codebox{70}{6}{255, 3855, 13107, 54613}{0}{7}{8}{0}{0}{0,\hfil 0, \hfil 0}\par\noindent
\codebox{71}{6}{255, 16131, 115471, 412723}{0}{5}{10}{0}{0}{10,\hfil 0, \hfil 10}\par\noindent
\codebox{72}{6}{255, 3855, 127249, 144998}{0}{5}{10}{0}{0}{7,\hfil 0, \hfil 3}\par\noindent
\codebox{73}{6}{255, 3855, 62211, 79157}{0}{5}{10}{0}{0}{4,\hfil 0, \hfil 0}\par\noindent
\codebox{74}{6}{255, 16131, 115471, 396597}{0}{3}{12}{0}{0}{3,\hfil 0, \hfil 1}\par\noindent
\codebox{75}{6}{255, 3855, 29491, 230741}{0}{3}{12}{0}{0}{0,\hfil 0, \hfil 0}\par\noindent
\smallskip\par\noindent
$\sigma=5. \qquad r(5)=58$.\smallskip\par\noindent
\codebox{76}{5}{31, 8160, 25059, 238049, 3618}{6}{0}{0}{10}{0}{0,\hfil 0, \hfil 0}\par\noindent
\codebox{77}{5}{31, 2019, 6244, 8637, 19179}{6}{0}{0}{0}{0}{0,\hfil 0, \hfil 0}\par\noindent
\codebox{78}{5}{31, 8160, 25059, 105991, 26232}{5}{8}{0}{10}{8}{0,\hfil 0, \hfil 0}\par\noindent
\codebox{79}{5}{31, 8160, 25059, 105991, 147041}{5}{4}{0}{2}{8}{0,\hfil 0, \hfil 0}\par\noindent
\codebox{80}{5}{31, 8160, 25059, 42605, 26781}{5}{4}{0}{1}{3}{3,\hfil 0, \hfil 0}\par\noindent
\codebox{81}{5}{31, 8161, 253987, 288708, 894990}{4}{7}{2}{0}{0}{0,\hfil 8, \hfil 0}\par\noindent
\codebox{82}{5}{31, 8160, 25059, 238119, 25661}{4}{7}{0}{1}{7}{4,\hfil 6, \hfil 0}\par\noindent
\codebox{83}{5}{31, 8160, 25059, 42605, 98704}{4}{7}{0}{1}{5}{8,\hfil 3, \hfil 0}\par\noindent
\codebox{84}{5}{31, 8160, 25059, 492069, 534498}{4}{7}{0}{0}{4}{10,\hfil 4, \hfil 4}\par\noindent
\codebox{85}{5}{31, 8160, 25059, 105991, 394851}{3}{13}{0}{1}{15}{24,\hfil 0, \hfil 0}\par\noindent
\codebox{86}{5}{31, 8160, 25059, 105991, 42605}{3}{13}{0}{1}{15}{0,\hfil 0, \hfil 0}\par\noindent
\codebox{87}{5}{31, 8160, 25059, 238119, 377379}{3}{13}{0}{1}{11}{28,\hfil 32, \hfil 8}\par\noindent
\codebox{88}{5}{31, 8160, 25059, 105991, 434281}{3}{13}{0}{1}{7}{32,\hfil 16, \hfil 24}\par\noindent
\codebox{89}{5}{31, 8160, 25059, 42605, 2724}{3}{13}{0}{1}{3}{12,\hfil 0, \hfil 0}\par\noindent
\codebox{90}{5}{31, 8161, 253987, 271302, 901198}{3}{9}{3}{0}{0}{27,\hfil 3, \hfil 27}\par\noindent
\codebox{91}{5}{31, 8160, 25059, 42605, 100414}{3}{9}{2}{0}{2}{13,\hfil 6, \hfil 6}\par\noindent
\codebox{92}{5}{31, 8160, 25059, 238119, 49277}{3}{9}{1}{0}{4}{17,\hfil 5, \hfil 7}\par\noindent
\codebox{93}{5}{31, 8160, 25059, 105991, 140901}{3}{9}{0}{1}{7}{8,\hfil 0, \hfil 0}\par\noindent
\codebox{94}{5}{31, 8160, 25059, 238119, 1736}{3}{9}{0}{1}{3}{18,\hfil 4, \hfil 6}\par\noindent
\codebox{95}{5}{31, 8160, 25059, 492069, 106180}{3}{9}{0}{0}{6}{15,\hfil 4, \hfil 6}\par\noindent
\codebox{96}{5}{31, 8160, 25059, 124512, 951009}{3}{9}{0}{0}{6}{9,\hfil 0, \hfil 0}\par\noindent
\codebox{97}{5}{31, 8160, 25059, 238119, 1869504}{2}{14}{0}{0}{8}{36,\hfil 22, \hfil 18}\par\noindent
\codebox{98}{5}{31, 8160, 25059, 492069, 1615373}{2}{14}{0}{0}{4}{42,\hfil 24, \hfil 32}\par\noindent
\codebox{99}{5}{31, 8160, 25059, 42605, 101942}{2}{14}{0}{0}{4}{30,\hfil 24, \hfil 16}\par\noindent
\codebox{100}{5}{31, 8160, 25059, 241184, 370273}{2}{10}{4}{0}{6}{12,\hfil 16, \hfil 0}\par\noindent
\codebox{101}{5}{31, 8160, 25059, 492069, 101592}{2}{10}{4}{0}{4}{24,\hfil 4, \hfil 20}\par\noindent
\codebox{102}{5}{31, 8160, 25059, 238119, 884843}{2}{10}{4}{0}{4}{18,\hfil 0, \hfil 0}\par\noindent
\codebox{103}{5}{31, 8160, 25059, 238119, 888353}{2}{10}{4}{0}{2}{24,\hfil 6, \hfil 18}\par\noindent
\codebox{104}{5}{31, 8161, 253987, 288708, 622825}{2}{10}{4}{0}{0}{30,\hfil 0, \hfil 32}\par\noindent
\codebox{105}{5}{31, 8161, 253987, 288708, 796873}{2}{10}{4}{0}{0}{24,\hfil 0, \hfil 16}\par\noindent
\codebox{106}{5}{31, 8161, 253987, 288708, 567406}{2}{10}{4}{0}{0}{12,\hfil 16, \hfil 0}\par\noindent
\codebox{107}{5}{31, 8160, 123360, 419424, 699040}{1}{30}{0}{0}{30}{0,\hfil 0, \hfil 0}\par\noindent
\codebox{108}{5}{31, 8160, 25059, 124512, 494240}{1}{22}{0}{0}{14}{56,\hfil 128, \hfil 0}\par\noindent
\codebox{109}{5}{31, 8160, 25059, 124512, 396941}{1}{18}{0}{0}{6}{60,\hfil 48, \hfil 32}\par\noindent
\codebox{110}{5}{31, 8160, 25059, 124512, 166317}{1}{18}{0}{0}{6}{54,\hfil 68, \hfil 24}\par\noindent
\codebox{111}{5}{31, 8160, 25059, 124512, 43685}{1}{18}{0}{0}{6}{36,\hfil 0, \hfil 0}\par\noindent
\codebox{112}{5}{31, 8160, 123360, 419424, 699041}{1}{14}{16}{0}{14}{0,\hfil 0, \hfil 0}\par\noindent
\codebox{113}{5}{31, 8160, 25059, 238119, 828508}{1}{14}{8}{0}{6}{40,\hfil 32, \hfil 24}\par\noindent
\codebox{114}{5}{31, 8160, 25059, 238119, 372292}{1}{14}{8}{0}{6}{40,\hfil 0, \hfil 16}\par\noindent
\codebox{115}{5}{31, 8160, 25059, 492069, 124520}{1}{14}{4}{0}{2}{48,\hfil 16, \hfil 44}\par\noindent
\codebox{116}{5}{31, 8160, 25059, 238119, 885801}{1}{14}{4}{0}{2}{42,\hfil 20, \hfil 28}\par\noindent
\codebox{117}{5}{31, 8160, 25059, 42605, 101044}{1}{14}{4}{0}{2}{24,\hfil 32, \hfil 12}\par\noindent
\codebox{118}{5}{31, 8160, 25059, 124512, 436897}{1}{10}{16}{0}{6}{12,\hfil 0, \hfil 0}\par\noindent
\codebox{119}{5}{31, 8160, 25059, 238119, 296165}{1}{10}{12}{0}{2}{26,\hfil 4, \hfil 20}\par\noindent
\codebox{120}{5}{31, 8160, 25059, 42605, 477857}{1}{10}{12}{0}{2}{20,\hfil 0, \hfil 12}\par\noindent
\codebox{121}{5}{31, 8161, 25062, 99051, 427305}{1}{10}{10}{0}{0}{30,\hfil 0, \hfil 30}\par\noindent
\codebox{122}{5}{31, 8161, 25062, 99051, 173347}{1}{10}{10}{0}{0}{24,\hfil 8, \hfil 18}\par\noindent
\codebox{123}{5}{255, 3855, 28951, 492145, 538402}{0}{25}{6}{0}{0}{60,\hfil 240, \hfil 0}\par\noindent
\codebox{124}{5}{255, 3855, 28951, 492145, 564498}{0}{21}{10}{0}{0}{66,\hfil 128, \hfil 14}\par\noindent
\codebox{125}{5}{255, 3855, 28951, 492145, 558755}{0}{21}{10}{0}{0}{60,\hfil 80, \hfil 0}\par\noindent
\codebox{126}{5}{255, 3855, 28951, 492145, 110650}{0}{17}{14}{0}{0}{58,\hfil 48, \hfil 30}\par\noindent
\codebox{127}{5}{255, 3855, 28951, 492145, 623923}{0}{17}{14}{0}{0}{52,\hfil 48, \hfil 24}\par\noindent
\codebox{128}{5}{255, 3855, 28951, 233577, 893570}{0}{13}{18}{0}{0}{42,\hfil 16, \hfil 34}\par\noindent
\codebox{129}{5}{255, 3855, 13107, 116021, 415508}{0}{13}{18}{0}{0}{42,\hfil 0, \hfil 30}\par\noindent
\codebox{130}{5}{255, 3855, 28951, 492145, 570411}{0}{13}{18}{0}{0}{36,\hfil 16, \hfil 24}\par\noindent
\codebox{131}{5}{255, 3855, 28951, 111147, 398693}{0}{9}{22}{0}{0}{24,\hfil 4, \hfil 28}\par\noindent
\codebox{132}{5}{255, 3855, 127249, 144998, 284986}{0}{9}{22}{0}{0}{24,\hfil 0, \hfil 20}\par\noindent
\codebox{133}{5}{255, 3855, 62211, 208947, 87381}{0}{9}{22}{0}{0}{18,\hfil 0, \hfil 6}\par\noindent
\smallskip\par\noindent
$\sigma=4. \qquad r(4)=41$.\smallskip\par\noindent
\codebox{134}{4}{31, 8160, 25059, 238119, 1736, 1867799}{7}{7}{0}{11}{9}{0,\hfil 0, \hfil 0}\par\noindent
\codebox{135}{4}{31, 8160, 25059, 105991, 394851, 139649}{7}{7}{0}{7}{21}{0,\hfil 0, \hfil 0}\par\noindent
\codebox{136}{4}{31, 8160, 25059, 105991, 434281, 614571}{7}{7}{0}{3}{9}{12,\hfil 0, \hfil 0}\par\noindent
\codebox{137}{4}{31, 8160, 25059, 238119, 884843, 418183}{6}{12}{0}{3}{15}{24,\hfil 30, \hfil 6}\par\noindent
\codebox{138}{4}{31, 8160, 25059, 42605, 2724, 987586}{6}{12}{0}{2}{6}{18,\hfil 18, \hfil 0}\par\noindent
\codebox{139}{4}{31, 8160, 25059, 492069, 534498, 1812520}{6}{12}{0}{0}{12}{30,\hfil 40, \hfil 0}\par\noindent
\codebox{140}{4}{31, 8160, 25059, 238119, 372292, 29575}{5}{24}{0}{10}{24}{96,\hfil 192, \hfil 64}\par\noindent
\codebox{141}{4}{31, 8160, 25059, 105991, 26232, 43689}{5}{24}{0}{10}{24}{0,\hfil 0, \hfil 0}\par\noindent
\codebox{142}{4}{31, 8160, 25059, 238119, 884843, 1058259}{5}{16}{0}{2}{16}{44,\hfil 40, \hfil 24}\par\noindent
\codebox{143}{4}{31, 8160, 25059, 238119, 884843, 7297}{5}{16}{0}{2}{16}{20,\hfil 48, \hfil 0}\par\noindent
\codebox{144}{4}{31, 8160, 25059, 238119, 49277, 516264}{5}{16}{0}{1}{11}{53,\hfil 44, \hfil 44}\par\noindent
\codebox{145}{4}{31, 8160, 25059, 238119, 884843, 1409677}{4}{19}{2}{0}{8}{74,\hfil 64, \hfil 74}\par\noindent
\codebox{146}{4}{31, 8160, 25059, 238119, 884843, 52788}{4}{19}{0}{1}{13}{70,\hfil 71, \hfil 58}\par\noindent
\codebox{147}{4}{31, 8160, 25059, 238119, 884843, 1474759}{4}{19}{0}{1}{9}{66,\hfil 43, \hfil 36}\par\noindent
\codebox{148}{4}{31, 8160, 25059, 238119, 49277, 984106}{4}{19}{0}{0}{12}{78,\hfil 58, \hfil 86}\par\noindent
\codebox{149}{4}{31, 8160, 25059, 238119, 372292, 103644}{3}{29}{0}{1}{23}{152,\hfil 272, \hfil 152}\par\noindent
\codebox{150}{4}{31, 8160, 25059, 105991, 394851, 696425}{3}{29}{0}{1}{15}{184,\hfil 224, \hfil 272}\par\noindent
\codebox{151}{4}{31, 8160, 25059, 238119, 377379, 950861}{3}{29}{0}{1}{15}{160,\hfil 272, \hfil 192}\par\noindent
\codebox{152}{4}{31, 8160, 25059, 238119, 49277, 281774}{3}{21}{4}{0}{6}{111,\hfil 64, \hfil 174}\par\noindent
\codebox{153}{4}{31, 8160, 25059, 238119, 884843, 1475209}{3}{21}{4}{0}{6}{87,\hfil 96, \hfil 98}\par\noindent
\codebox{154}{4}{31, 8160, 25059, 238119, 884843, 1451537}{3}{21}{2}{0}{10}{95,\hfil 74, \hfil 104}\par\noindent
\codebox{155}{4}{31, 8160, 25059, 238119, 884843, 1352755}{3}{21}{0}{1}{15}{72,\hfil 0, \hfil 0}\par\noindent
\codebox{156}{4}{31, 8160, 25059, 105991, 42605, 141990}{3}{21}{0}{1}{15}{48,\hfil 128, \hfil 0}\par\noindent
\codebox{157}{4}{31, 8160, 25059, 238119, 372292, 699489}{3}{21}{0}{1}{7}{104,\hfil 64, \hfil 144}\par\noindent
\codebox{158}{4}{31, 8160, 25059, 238119, 1869504, 475241}{2}{30}{0}{0}{12}{186,\hfil 276, \hfil 244}\par\noindent
\codebox{159}{4}{31, 8160, 25059, 238119, 1869504, 1902665}{2}{30}{0}{0}{12}{162,\hfil 276, \hfil 180}\par\noindent
\codebox{160}{4}{31, 8160, 25059, 238119, 884843, 321232}{2}{22}{8}{0}{8}{110,\hfil 90, \hfil 150}\par\noindent
\codebox{161}{4}{31, 8160, 25059, 238119, 884843, 167565}{2}{22}{8}{0}{4}{122,\hfil 72, \hfil 192}\par\noindent
\codebox{162}{4}{31, 8160, 25059, 238119, 888353, 1355336}{2}{22}{8}{0}{4}{122,\hfil 64, \hfil 200}\par\noindent
\codebox{163}{4}{31, 8160, 25059, 124512, 494240, 700700}{1}{46}{0}{0}{30}{240,\hfil 1280, \hfil 0}\par\noindent
\codebox{164}{4}{31, 8160, 25059, 124512, 396941, 662065}{1}{38}{0}{0}{14}{240,\hfil 720, \hfil 192}\par\noindent
\codebox{165}{4}{31, 8160, 25059, 238119, 372292, 955584}{1}{30}{16}{0}{14}{176,\hfil 256, \hfil 192}\par\noindent
\codebox{166}{4}{31, 8160, 25059, 238119, 372292, 442537}{1}{30}{8}{0}{6}{192,\hfil 272, \hfil 256}\par\noindent
\codebox{167}{4}{31, 8160, 25059, 238119, 372292, 950861}{1}{30}{8}{0}{6}{192,\hfil 208, \hfil 240}\par\noindent
\codebox{168}{4}{31, 8160, 25059, 238119, 372292, 829089}{1}{22}{24}{0}{6}{120,\hfil 48, \hfil 176}\par\noindent
\codebox{169}{4}{31, 8160, 25059, 238119, 296165, 591468}{1}{22}{20}{0}{2}{128,\hfil 64, \hfil 220}\par\noindent
\codebox{170}{4}{255, 3855, 28951, 492145, 564498, 42406}{0}{45}{18}{0}{0}{270,\hfil 1440, \hfil 90}\par\noindent
\codebox{171}{4}{255, 3855, 28951, 492145, 564498, 722490}{0}{37}{26}{0}{0}{246,\hfil 640, \hfil 210}\par\noindent
\codebox{172}{4}{255, 3855, 28951, 492145, 564498, 1127602}{0}{29}{34}{0}{0}{190,\hfil 224, \hfil 266}\par\noindent
\codebox{173}{4}{255, 3855, 28951, 233577, 893570, 308270}{0}{21}{42}{0}{0}{126,\hfil 56, \hfil 238}\par\noindent
\codebox{174}{4}{255, 3855, 13107, 116021, 415508, 714818}{0}{21}{42}{0}{0}{126,\hfil 0, \hfil 210}\par\noindent
\smallskip\par\noindent
$\sigma=3. \qquad r(3)=13$.\smallskip\par\noindent
\codebox{175}{3}{31, 8160, 25059, 238119, 884843, 1474759, 475241}{9}{18}{0}{20}{18}{0,\hfil 0, \hfil 0}\par\noindent
\codebox{176}{3}{31, 8160, 25059, 238119, 884843, 418183, 1451537}{9}{18}{0}{16}{30}{48,\hfil 96, \hfil 16}\par\noindent
\codebox{177}{3}{31, 8160, 25059, 238119, 884843, 418183, 57025}{9}{18}{0}{9}{27}{63,\hfil 102, \hfil 0}\par\noindent
\codebox{178}{3}{31, 8160, 25059, 238119, 884843, 418183, 699489}{7}{31}{0}{5}{35}{182,\hfil 374, \hfil 228}\par\noindent
\codebox{179}{3}{31, 8160, 25059, 238119, 884843, 1409677, 1058259}{7}{31}{0}{3}{33}{204,\hfil 368, \hfil 288}\par\noindent
\codebox{180}{3}{31, 8160, 25059, 238119, 372292, 29575, 955584}{5}{56}{0}{10}{56}{576,\hfil 2176, \hfil 1152}\par\noindent
\codebox{181}{3}{31, 8160, 25059, 238119, 884843, 1451537, 699489}{5}{40}{0}{2}{32}{324,\hfil 688, \hfil 608}\par\noindent
\codebox{182}{3}{31, 8160, 25059, 238119, 884843, 1451537, 1474759}{5}{40}{0}{1}{27}{357,\hfil 628, \hfil 804}\par\noindent
\codebox{183}{3}{31, 8160, 25059, 238119, 372292, 442537, 934222}{3}{61}{0}{1}{39}{744,\hfil 2640, \hfil 1800}\par\noindent
\codebox{184}{3}{31, 8160, 25059, 238119, 884843, 1451537, 167565}{3}{45}{6}{0}{18}{495,\hfil 774, \hfil 1476}\par\noindent
\codebox{185}{3}{31, 8160, 25059, 238119, 884843, 167565, 1352755}{3}{45}{0}{1}{15}{504,\hfil 672, \hfil 1520}\par\noindent
\codebox{186}{3}{31, 8160, 25059, 124512, 396941, 662065, 700700}{1}{78}{0}{0}{30}{1008,\hfil 6720, \hfil 1536}\par\noindent
\codebox{187}{3}{31, 8160, 25059, 238119, 372292, 442537, 955584}{1}{62}{16}{0}{14}{816,\hfil 2624, \hfil 2112}\par\noindent
\smallskip\par\noindent
$\sigma=2. \qquad r(2)=3$.\smallskip\par\noindent
\codebox{188}{2}{31, 8160, 25059, 238119, 884843, 418183, 1451537, 699489}{13}{28}{0}{46}{60}{96,\hfil 416, \hfil 0}\par\noindent
\codebox{189}{2}{31, 8160, 25059, 238119, 884843, 418183, 699489, 152785}{9}{66}{0}{12}{90}{864,\hfil 3672, \hfil 2448}\par\noindent
\codebox{190}{2}{31, 8160, 25059, 238119, 372292, 442537, 934222, 1844576}{5}{120}{0}{10}{120}{2880,\hfil 21120, \hfil 13440}\par\noindent
\smallskip\par\noindent
$\sigma=1. \qquad r(1)=1$.\smallskip\par\noindent
\codebox{191}{1}{31, 8160, 25059, 238119, 884843, 418183, 1451537, 699489, 929948}{21}{0}{0}{210}{0}{0,\hfil 0, \hfil 0}\par\noindent

} 
\par
\medskip
\begin{remark}
Using Proposition~\ref{prop:isin},
we have also made the complete list of   pairs $([\tC], [\tC\sprime])$
of geometrically realizable classes of codes
satisfying  $[\tC]< [\tC\sprime]$.
\end{remark}
\subsection{Proof of Corollaries.}
In this subsection,
we prove Corollaries~\ref{cor:nopsirred},~\ref{cor:le9} and~\ref{cor:equalto1}
that are stated in Introduction.
We denote by 
$\CCL{\nu}$   the geometrically realizable class 
of No.\,$\nu$ in the list.
\begin{proof}[Proof of Corollary~\ref{cor:nopsirred}]
Note that
$$
\UUU\sb{\sigma}=\bigsqcup\sb{11-\dim \tC=\sigma } {\,\UUU}\sb{2, 6,[\tC]}.
$$
Let $\tU\sb{\sigma}$ be the pull-back of ${\,\UUU}\sb{\sigma}$
by the \'etale covering $\tUts\to\Uts$ constructed in the proof of Theorem~\ref{thm:closed}.
The code 
$\tau\sb G\inv (\CCG)$ in  $\Pow (\tZ)$
does not vary when $(G, \tau\sb G)$ moves on
an irreducible component of $\tU\sb{\sigma}$.
Hence each irreducible component of ${\,\UUU}\sb{\sigma}$ is contained 
in a unique $ {\,\UUU}\sb{2, 6,[\tC]}$  with $\dim\tC=11-\sigma$.
Therefore the number of the irreducible components of ${\,\UUU}\sb{\sigma}$ is greater than or equal to
 the number $r(\sigma)$
of  geometrically realizable  classes $[\tC]$ of codes with $\dim\tC=11-\sigma$.
\end{proof}
\begin{proof}[Proof of Corollary~\ref{cor:le9}]
Let $G$ be a polynomial  in  $\Uts$.
The Artin invariant of $\XG$ is $<10$
if and only if there exists a reduced  irreducible curve of degree $\le 2$
splitting in $\XG$,
or there exists a regular pencil of cubic curves splitting in $\XG$.
If there is a line (resp. a smooth conic)
splitting in $\XG$,
then $G\in {\,\UUU}[51]$ (resp. $G\in {\,\UUU}[42]$)
by Proposition~\ref{prop:convGC}.
If there is a regular  pencil of cubic curves splitting in $\XG$,
then $G\in {\,\UUU}[33]$ by Corollary~\ref{cor:EE}.
\par
It is obvious that the loci ${\,\UUU}[51]$, ${\,\UUU}[42]$ and ${\,\UUU}[33]$ are irreducible.
Because the locus $k\sptimes G +\Vts$ is closed in $\Uts$ for any $G\in \Uts$,
these loci  are Zariski closed in $\Uts$.
Because of the existence of
the geometrically realizable class  $\CCL{0}$,
Proposition~\ref{prop:GC}
implies that ${\,\UUU}[51]$, ${\,\UUU}[42]$ and ${\,\UUU}[33]$  are proper subsets of $\Uts$.
Therefore it remains to show that the codimension of these loci in $\Uts$ is $\le 1$.
\par
Let $\tUts\to \Uts$ be 
the \'etale covering   
 that has appeared in the proof of Theorem~\ref{thm:closed}.
We choose  six elements $\tP\sb 1$, \dots, $\tP\sb 6$ of $\tZ$, 
and consider the locus
\begin{equation}\label{eq:six}
\Biggl\{ 
\,
(G, \tau\sb G)\in \tUts
\;\;
\Bigg|
\;\;
\vcenter{
\hbox{
\vbox{
\hbox{there exists a smooth conic passing through}
\hbox{$\tau\sb G (\tP\sb 1)$, \dots, $\tau\sb G(\tP\sb 6)$}
}
}
}
\Biggr\}
\end{equation}
of $\tUts$.
Because of the existence of
the geometrically realizable class  $\CCL{2}$,
for example,
the locus~\eqref{eq:six} is non-empty.
Since $\dim |\OPt (2)|=5$,
the locus~\eqref{eq:six} is of codimension  $\le 1$ in $\tUts$.
If $(G, \tau\sb G)$ is in the locus~\eqref{eq:six},
then there exists a smooth conic splitting in $\XG$ by Proposition~\ref{prop:convnum},
and hence 
 $G$ is contained in ${\,\UUU}[42]$ by  Proposition~\ref{prop:convGC}. 
Therefore  the codimension of ${\,\UUU}[42]$ in $\Uts$ is also  $\le 1$.
The fact that ${\,\UUU}[51]\st \Uts$ is of codimension $1$ is proved in a similar way.
\par
Because of the existence of
the geometrically realizable class  $\CCL{3}$,
if $G$ is a general point of ${\,\UUU}[33]$,
then there exists only one 
regular pencil of cubic curves splitting in $\XG$.
Consider the morphism
$$
\map{\varrho}{\Hz {\OPt(3)}\times \Hz {\OPt(3)}\times k\sptimes \times  \Hz {\OPt(3)}
}{\Hz {\OPt(6)}}
$$
defined by
$$
(G\sb 3, G\sprime\sb 3, c, H)\;\mapsto\; c\, G\sb 3 G\sprime\sb 3 +H^2.
$$
Let $G\sb 3$ and $G\sb 3\sprime$ be general homogeneous polynomials of degree $3$.
Suppose that
$$
\varrho(G\sb 3, G\sprime\sb 3, 1, 0)=\varrho(\Gamma\sb 3, \Gamma\sprime\sb 3, c, H).
$$
Then the pencil of cubic curves spanned by
the curves defined by $G\sb 3=0$ and $G\sprime\sb 3=0$
coincides with 
the pencil  spanned by
the curves defined by $\Gamma\sb 3=0$ and $\Gamma\sprime\sb 3=0$.
Hence there exists an invertible matrix
$$
\begin{pmatrix} s & t \\ u & v \end{pmatrix}
$$
such that
$$
G\sb 3 = s\Gamma\sb 3 + t\Gamma\sprime\sb 3
\quand
G\sprime \sb 3 = u\Gamma\sb 3 + v\Gamma\sprime\sb 3
$$
hold.
Then we have
$$
c=sv+tu
\quand H=\sqrt{su} \,\Gamma\sb 3  +\sqrt{tv} \, \Gamma\sprime\sb 3.
$$
Hence we have
$$
\dim {\,\UUU}[33]= 3 \hz{\OPt (3)} +1 -\dim \GL (2)=27=\dim \Uts-1.
$$
Therefore ${\,\UUU}[33]$ is a hypersurface of $\Uts$.
\end{proof}
\begin{proof}[Proof of Corollary~\ref{cor:equalto1}]
Let $G\sb{\DK}$ be the Dolgachev-Kondo  polynomial~\eqref{eq:DK}.
Note that  $Z (dG\sb{\DK})$ coincides with  the set $\Pt (\F\sb 4)$ of $\F\sb 4$-rational points of $\Pt$,
and hence the set of lines splitting in $X\sb{G\sb{\DK}}$ is equal to 
the set $(\Pt)\dual (\F\sb 4)$
of $\F\sb 4$-rational lines of $\Pt$.
\par
Let $G$ be a polynomial in $\Uts$ such that the Artin invariant of  $\XG$ is $1$.
It is enough to show that,
if we choose homogeneous coordinates of $\Pt$ appropriately,
then $G$ is contained in $k\sptimes G\sb{\DK} +\Vts$.
Let $\LLL\sb{G}$
be the set of lines splitting in $\XG$.
Since there exists only one geometrically realizable  class $\CCL{191}$ 
with Artin invariant $1$,
the configuration $(\LLL\sb G, \ZdG)$ of lines and points
is isomorphic as abstract configurations (see ~\cite{DolCon}) to $((\Pt)\dual (\F\sb 4), \Pt(\F\sb 4))$.
In particular,
for any two points $P, Q\in \ZdG$,
the line $\overline{PQ}$ passing through  $P$ and $Q$ is in $\LLL\sb G$.
By choosing suitable homogeneous coordinates $[X\sb 0, X\sb 1, X\sb 2]$
and numbering the lines $\LLL\sb G=\{L\sb 0, \dots, L\sb{20}\}$ appropriately,
we can assume that
\begin{eqnarray*}
L\sb 0=\{X\sb 2=0\}, \quad 
L\sb 1=\{X\sb 1=0\}, \quad
L\sb 2=\{X\sb 1=X\sb 2\}, \quad
L\sb 3=\{X\sb 0=0\}, \\
L\sb 4=\{X\sb 0=X\sb 2\}, \quad
L\sb 5=\{X\sb 0=X\sb 1\}, \quad
L\sb 6=\{X\sb 0 +X\sb 1 + X\sb 2=0\}.
\end{eqnarray*}
The following points are in $\ZdG$:
$$
P\sb 0 :=L\sb 0\wedge L\sb 1 =[1,0,0], \quad
P\sb 1 :=L\sb 0\wedge L\sb 3 =[0,1,0], \quad
P\sb 2 :=L\sb 3\wedge L\sb 1 =[0,0,1].
$$
There exists a point $Q\sb 0 :=[\alpha, 0, 1]$ in $L\sb 1\cap \ZdG$ with $\alpha\ne 0, 1$.
Then we have
\begin{eqnarray*}
L\sb 7 &:= & \overline{P\sb 1 Q\sb 0 }=\{X\sb 0= \alpha  X\sb 2\}\in \LLL\sb G, \\
Q\sb 1&:= & L\sb 5\wedge L\sb 7=[\alpha,\alpha,1]\in \ZdG, \\
L\sb 8&:= & \overline{P\sb 0 Q\sb 1}=\{X\sb 1= \alpha  X\sb 2\} \in \LLL\sb G, \\
Q\sb 2 &:= &L\sb 6\wedge L\sb 8 =[1+\alpha, \alpha, 1]\in\ZdG, \\
L\sb 9&:= &\overline{P\sb 1 Q\sb 2}=\{X\sb 0= (1+\alpha) X\sb 2\}\in \LLL\sb G.
\end{eqnarray*}
\begin{figure}
\begin{center}
\setlength{\unitlength}{1.1pt}
\begin{picture}(200, 170)
\thinlines
\put(40,10){\line(0, 1){125}}
\put(120,10){\line(0, 1){125}}
\put(20,30){\line(1, 0){125}}
\put(20,110){\line(1,0){125}}
\put(20,10){\line(1,1){120}}
\put(20,130){\line(1,-1){120}}
\put(60,10){\line(0, 1){125}}
\put(100,10){\line(0, 1){125}}
\put(20,50){\line(1, 0){125}}
\thicklines
\put(35, 10){\line(1,4){31}}
\put(40, 144){\makebox(0,0){$L\sb 3$}}
\put(60, 144){\makebox(0,0){$L\sb 7$}}
\put(77, 134){\makebox(0,0){$L\sb {10}$}}
\put(100, 144){\makebox(0,0){$L\sb 9$}}
\put(120, 144){\makebox(0,0){$L\sb 4$}}
\put(158, 110){\makebox(0,0){$L\sb 2$}}
\put(158, 50){\makebox(0,0){$L\sb 8$}}
\put(158, 30){\makebox(0,0){$L\sb 1$}}
\put(148, 130){\makebox(0,0){$L\sb 5$}}
\put(148, 10){\makebox(0,0){$L\sb 6$}}
\put(60, 30){\circle*{4}}
\put(60, 110){\circle*{4}}
\put(64, 20){{$Q\sb 0$}}
\put(64, 100){{$R$}}
\end{picture}
\end{center}
\caption{Lines in $\LLL\sb G$}
\end{figure}
The five  points consisting  $L\sb 9\cap \ZdG$ are  therefore
\begin{multline*}
P\sb 1=[0,1,0], \quad Q\sb 2=[1+\alpha, \alpha, 1], \quad L\sb 2\wedge L\sb 9=[1+\alpha, 1, 1],\\
 L\sb 5\wedge L\sb 9=[1+\alpha, 1+\alpha, 1],  \quand L\sb 1\wedge L\sb 9=[1+\alpha, 0, 1].
\end{multline*}
On the other hand,
the point
$R:=L\sb 7\wedge L\sb 2=[\alpha,1,1]$ is contained in $\ZdG$,
and hence a line 
$$
L\sb{10}:=\overline{P\sb 2 R}=\{ X\sb 0+ \alpha X\sb 1=0\}
$$
is an element of  $\LLL\sb G$.
The point 
$$
L\sb{10}\wedge L\sb 9=[\alpha^2+\alpha, \alpha+1, \alpha]
$$
is therefore among the five points above.
Because $\alpha\ne 0, 1$,
this point must be $Q\sb 2$, and $\alpha$ is a root of $t^2+t+1=0$.
Then we can show that all points of $\ZdG$ are $\F\sb4$-rational,
and hence $\ZdG=Z(dG\sb{\DK})$ holds.
By the uniqueness assertion of Theorem~\ref{thm:mainsecglobal},
we have $dG= c \cdot dG\sb {DK}$, where $c$ is a  non-zero constant.
Since $\Vts$ is the kernel of the linear homomorphism $G\mapsto dG$,
we have $G\in  k\sptimes G\sb{\DK}+ \Vts$.
\end{proof}
\section{The algorithm}\label{sec:algorithm}
\subsection{The description of the algorithm.}
We present an algorithm that calculates
the code $\CCG$ from a given homogeneous  polynomial $G\in\Uts$.
From the results in the previous sections, we obtain the following:
\begin{corollary}
Let $G$ be a polynomial  in  $\Uts$.
\par
{\rm (1)}
A subset $B\st\ZdG$ of weight $5$ is contained in $\CCG$ if and only if
the points of $B$ are collinear.
\par
{\rm (2)}
 Let $B\st\ZdG$ be a subset of weight $8$ such that no three points of $B$ are collinear.
Then $B$ is contained in $\CCG$ if and only if there exists a  conic containing $B$.
{\rm (}Note that,  if such a conic exists,
then it must be smooth because no three points of $B$ are collinear.{\rm )} 
\end{corollary}
\begin{corollary}\label{cor:EEEalgo}
Let $G$ be a polynomial  in  $\Uts$, and 
let $B\st\ZdG$ be a subset of weight $9$ such that no three points of $B$ are collinear.
Then $B$ is contained in $\CCG$ if and only if
the following hold; %{;}
{\rm (i)}
the linear system $|\III\sb{B} (3)|$ of cubic curves containing $B$ is of dimension $1$,
and
{\rm (ii)} if $\Hz{\III\sb{B} (3)}$ is generated by $G\sb E$ and $G\sb{E\sprime}$,
 then $G\sb E G\sb {E\sprime}$ is contained in $ k\sptimes G+ \Vts$.
\end{corollary}
\begin{proof}
If $B\in \CCG$, then
(i) and (ii) hold by Proposition~\ref{prop:cubic} and Corollaries~\ref{cor:collinear} and~\ref{cor:EE}.
Suppose that  (i) and (ii) hold,
and let $E$ and $E\sprime$ be the cubic curves defined by $G\sb E=0$ and $G\sb{E\sprime}=0$.
Then  $E$ and $E\sprime$ are splitting in $\XG$, and 
$$
B=E\cap E\sprime=\wG (E)=\wG(E\sprime)
$$
holds by Proposition~\ref{prop:GC}.
Hence  $B$ is contained in $\CCG$.
\end{proof}
\begin{remark}
In Corollary~\ref{cor:EEEalgo},
the condition (i) alone  is not enough for $B$ to be contained in $\CCG$.
See Example~\ref{example:1416}.
\end{remark}
\begin{algorithm}\label{algo:main}
Suppose that we are given a homogeneous polynomial $G\in \Uts$.
This algorithm  outputs a  set $\tGen=\{ A\sb 0, \dots, A\sb {k-1}\}\subset \Pow(\ZdG)$ 
that generates $\CCG$,
and the Artin invariant of $\XG$.
\par
{\it Step 1.}
Set $\tGen$ to be $\emptyset$.
\par
{\it Step 2.}
Calculate the coordinates of the points $P\sb 0, \dots, P\sb{20}$ of $\ZdG$
by solving 
$$
\Der{G}{X\sb 0}=\Der{G}{X\sb 1}= \Der{G}{X\sb 2}=0.
$$
\par
{\it Step 3.}
Put the word $\ZdG=\{ P\sb 0, \dots, P\sb{20}\}$ in $\tGen$.
\par
{\it Step 4.}
Make the list $\tCol$ of all triples $\{P\sb i, P\sb j, P\sb k\}$
of points of $\ZdG$ 
that are collinear.
\par
{\it Step 5.}
Using $\tCol$,
list up  all $5$-tuples $\{P\sb{i\sb 1}, \dots, P\sb {i\sb 5}\}$ that are collinear,
and put them in $\tGen$.
By Proposition~\ref{prop:convnum},
every triple in  $\tCol$ must extend to a collinear $5$-tuple.
\par
{\it Step 6.}
For each $8$-tuple $B=\{P\sb{i\sb 1}, \dots, P\sb {i\sb 8}\}$ of  points of $\ZdG$,
check whether there exist   collinear three  points of $B$  by using $\tCol$.
If there are no such three points,
then check whether there exists a conic that passes through  the points of $B$.
If such a conic exists, then   put $B$ in $\tGen$.
\par
{\it Step 7.}
For each $9$-tuple $B=\{P\sb{i\sb 1}, \dots, P\sb {i\sb 9}\}$,
check whether there exist   collinear three  points of $B$  by using $\tCol$.
If there are no such three points,
then calculate $\dim |\III\sb B (3)|$.
If $\dim |\III\sb B (3)|=1$,
choose polynomials   $G\sb E$ and $G\sb{E\sprime}$ that span  $\Hz{\III\sb B (3)}$,
and check whether $G\sb E G\sb {E\sprime}$ is contained in $k\sptimes G+ \Vts$ or not by using
the method described in Remark~\ref{rem:equivalgo}.
If $G\sb E G\sb {E\sprime}\in  k\sptimes G+ \Vts$, then  put $B$ in $\tGen$.
\par
{\it Step 8.}
Calculate the code $\CCG$ generated by the words in $\tGen$.
The Artin invariant of $\XG$ is $11-\dim\CCG$.
\end{algorithm}
\subsection{Examples.}
\begin{example}\label{example:792}
The code $\CCG$ of the polynomial $G$ in Example~\ref{example:1598} is in the class $\CCL{135}$.
Let us consider the polynomial
\begin{multline*}
G\sprime:={X_{{0}}}^{5}X_{{2}}+{X_{{0}}}^{4}X_{{1}}X_{{2}}+{X_{{0}}}^{3}{X_{{1}}
}^{2}X_{{2}}+{X_{{0}}}^{2}{X_{{1}}}^{3}X_{{2}}+\\+X_{{0}}{X_{{1}}}^{4}X_{
{2}}+X_{{0}}{X_{{1}}}^{3}{X_{{2}}}^{2}+X_{{0}}X_{{1}}{X_{{2}}}^{4}.
\end{multline*}
The points of $Z (dG\sprime)$ are defined over $\F\sb{2^{24}}$.
Under the Frobenius morphism over $\F\sb 2$,
they are decomposed into six orbits,
the cardinalities of which are $1$, $1$, $3$, $4$, $4$, $8$.
The set of curves of degree $\le 3$ splitting in $X\sb{G\sprime}$
consists of  seven lines,
which are decomposed into four   Frobenius orbits of cardinalities $1$, $1$, $1$,  $4$, 
 and seven smooth conics,
which are  decomposed into three  Frobenius orbits of cardinalities $1$, $2$, $4$.
The class $[\CCC\sb{G\sprime}]$ 
is $\CCL{134}$.
\end{example}
\begin{example}\label{example:643}
Consider the polynomial 
$$
G:={X_{{0}}}^{4}X_{{1}}X_{{2}}+{X_{{0}}}^{3}{X_{{1}}}^{3}+X_{{0}}{X_{{1}}
}^{4}X_{{2}}+X_{{0}}X_{{1}}{X_{{2}}}^{4}.
$$
The subscheme $\ZdG$ is reduced of dimension $0$,
and each point is defined over $\F\sb{2^4}$.
The class of the code $\CCG$ is $\CCL{190}$.
In particular, the Artin invariant of $\XG$ is $2$.
\end{example}
\begin{example}\label{example:1416}
We will give an example of $\XG$ with Artin invariant $3$.
Consider the polynomial 
\begin{multline*}
G:={X_{{0}}}^{5}X_{{2}}+{X_{{0}}}^{4}X_{{1}}X_{{2}}+{X_{{0}}}^{3}{X_{{1}}
}^{3}+{X_{{0}}}^{3}{X_{{1}}}^{2}X_{{2}}+{X_{{0}}}^{3}{X_{{2}}}^{3}+{X_
{{0}}}^{2}{X_{{1}}}^{3}X_{{2}}+\\+X_{{0}}{X_{{1}}}^{3}{X_{{2}}}^{2}+X_{{0
}}X_{{1}}{X_{{2}}}^{4}+{X_{{1}}}^{5}X_{{2}}.
\end{multline*}
Let $\alpha$ be a root of the irreducible polynomial 
$$
{t}^{6}+{t}^{5}+{t}^{3}+{t}^{2}+1 \;\in\; \F\sb 2[t].
$$
Then $\ZdG$ consists of the points in Table~\ref{table:1416ZdG}.
\begin{table}
\begin{eqnarray*}
&&P_{{0}}=[{\alpha}^{5}+{\alpha}^{3}+\alpha+1,{\alpha}^{3}+{\alpha}^{2}+
\alpha+1,1], 
 \\
&&P_{{\nu}}=\Frob \sp{\nu} (P_0)\quad (\nu=1, \dots, 5), \\
&&P_{{6}}=[1,1,1], \quad P_{{7}}=[1,0,1], \\
&&P_{{8}}=[{\alpha}^{4}+{\alpha}^{3}+{\alpha}^{2}+\alpha,\alpha+1,1], \\
&&P_{{8+\nu}}=\Frob \sp{\nu} (P_8)\quad (\nu=1, \dots, 5),
\\
&& P_{{14}}=[0,0,1], \\
&&P_{{15}}=[{\alpha}^{5}+{\alpha}^{4}+{\alpha}^{3}+{\alpha}^{2}+1,{
\alpha}^{5}+{\alpha}^{4}+{\alpha}^{3}+{\alpha}^{2}+\alpha,1], \\
&&P_{{15+\nu}}=\Frob \sp{\nu} (P_{15})\quad (\nu=1, \dots, 5).
\end{eqnarray*}
\vskip 3pt
\caption{Points of $\ZdG$ in Example~\ref{example:1416}} \label{table:1416ZdG}
\end{table}
The words of weight $5$ in $\CCG$ are
$$
\{0, 3, 6, 16, 19\}, \quad  \{1, 4, 6, 17, 20\}, \quad \{2, 5, 6, 15, 18\},
$$
which form one  Frobenius orbit,
where
the set $\{ P\sb{i\sb 1}, \dots, P\sb{i\sb k}\}$ is simply denoted by $\{i\sb 1, \dots, i\sb k\}$.
There are $45$ irreducible words of weight $8$  in $\CCG$.
The cardinalities of Frobenius orbits are
$$
1,6,6,2,2,6,6,3,6,6,1.
$$
There are no irreducible words of weight $9$  in $\CCG$.
The class  $[\CCG]$ is  $\CCL{185}$.
In particular, the Artin invariant of $\XG$ is $3$.
\par
Consider the following word of weight $9$;
$$
A:=\{0,1,2,3,7,8,9,15,20\}.
$$
Note that 
no three points of $A$ are collinear.
There exists a  pencil of cubic curves whose base locus is $A$,
which is spanned by
{
\newcommand{\salpha}{\hskip -.5pt\alpha}
\begin{multline*}
{X_{{0}}}^{2}X_{{1}}+
\left ({\salpha}^{4}+{\salpha}^{3}+{\salpha}^{2}\right ){X_{{0}}}^{2}X_{{2}}+
\left ({\salpha}^{5}+{\salpha}^{4}+{\salpha}^{2}\right )X_{{0}}{X_{{1}}}^{2} +\\+ 
\left ({\salpha}^{5}+{\salpha}^{4}+\salpha+1\right )X_{{0}}{X_{{2}}}^{2}+
\left ({\salpha}^{4}+{\salpha}^{3}+1\right ){X_{{1}}}^{3}+
\left ({\salpha}^{4}+{\salpha}^{3}+\salpha\right ){X_{{1}}}^{2}X_{{2}}+\\+
\left ({\salpha}^{4}+{\salpha}^{3}+1\right)X_{{1}}{X_{{2}}}^{2}+
\left ({\salpha}^{5}+{\salpha}^{3}+{\salpha}^{2}+\salpha+1\right){X_{{2}}}^{3}=0, 
\end{multline*}
and 
\begin{multline*}
{X_{{0}}}^{3}+\left ({\salpha}^{4}+\salpha\right ){X_{{0}}}^{2}X_{{2}}+
\left ({\salpha}^{5}+{\salpha}^{3}\right )X_{{0}}{X_{{1}}}^{2}+
\left ({\salpha}^{3}+{\salpha}^{2}+1\right )X_{{0}}{X_{{2}}}^{2} +\\
+{\salpha}^{3}{X_{{1}}}^{3}+
\left ({\salpha}^{5}+{\salpha}^{4}+{\salpha}^{2}+\salpha+1\right ){X_{{1}}}^{2}X_{{2}}+
\left ({\salpha}^{3}+1\right )X_{{1}}{X_{{2}}}^{2}+\\
+\left ({\salpha}^{4}+{\salpha}^{3}+{\salpha}^{2}+\salpha\right ){X_{{2}}}^{3}=0.
\end{multline*}
}
However this pencil is not splitting in $\XG$.
\end{example}
\subsection{Irreducibility of  ${\,\UUU}\sb{2, 6, \CCL{}}$ for some $\CCL{}$.}
For some geometrically realizable classes $\CCL{}$,
we can prove the irreducibility of the locus ${\,\UUU}\sb{2, 6, \CCL{}}$,
and give a homogeneous polynomial $G$ that corresponds to the generic point of  ${\,\UUU}\sb{2, 6, \CCL{}}$.
\begin{definition}
For a non-increasing sequence $[a\sb 1 \dots a\sb k]$
of positive integers with $a\sb 1+\cdots +a\sb k=6$,
we denote by ${\,\UUU}[a\sb1\dots a\sb k]$
the locus of $G\in\Uts$ such that
there exist homogeneous polynomials
$G\sb{a\sb 1}$, \dots, $G\sb{a\sb k}$
of degrees $a\sb 1$, \dots, $a\sb k$
satisfying
$$
G\sb{a\sb 1} \cdots G\sb{a\sb k} \in k\sptimes G +\Vts.
$$
\end{definition}
It is obvious that ${\,\UUU}[a\sb1\dots a\sb k]$ is an irreducible Zariski closed subset of $\Uts$.
\begin{example}\label{example:nodal2211}
Let $G$ be a point of ${\,\UUU}[2211]$.
By Proposition~\ref{prop:GC},
there exist splitting lines $L\sb 1$, $L\sb 2$
and splitting smooth conics $Q\sb 1$, $Q\sb 2$
such that the union $L\sb 1\cup L\sb 2\cup Q\sb 1\cup Q\sb 2$
has only ordinary nodes as its singularities.
Hence $\CCG$ contains words
$A\sb 1$, $A\sb 2$, $B\sb 1$, $B\sb 2$
satisfying the following:
\begin{itemize}
\item $\card{A\sb 1}=\card{A\sb 2}=5$, $\card{B\sb 1}=\card{B\sb 2}=8$,
\item $B\sb 1$ and $B\sb 2$ are irreducible in $\CCG$,
\item $\card{A\sb i\cap B\sb j}=2$ for $i, j=1, 2$, and $\card{B\sb 1\cap B\sb 2}=4$,
\item$\card{A\sb 1\cap A\sb 2\cap B\sb j}= \card{A\sb i\cap B\sb 1\cap B\sb 2}=0$ for $i, j=1, 2$.
\end{itemize}
Conversely,
suppose that
the code $\CCG$ of  a polynomial  $G\in \Uts$ contains words 
$A\sb 1$, $A\sb 2$, $B\sb 1$, $B\sb 2$
satisfying the conditions above.
By Propositions~\ref{prop:line} and~\ref{prop:conic},
there exist lines $L\sb 1$, $L\sb 2$ and smooth conics $Q\sb 1$, $Q\sb 2$
splitting $\XG$ 
such that 
$L\sb i\cap\ZdG=A\sb i$
and $Q\sb j\cap \ZdG=B\sb j$ hold.
By Remarks~\ref{rem:intersection_lines},~\ref{rem:intersection_line_and_conic} and~\ref{rem:conicintersection}, 
the union $L\sb 1\cup L\sb 2\cup Q\sb 1\cup Q\sb 2$
has only ordinary nodes as its singularities.
Hence, by Proposition~\ref{prop:convGC},
$G$ is a point of ${\,\UUU}[2211]$.
\par
If
$[\CCG]=\CCL{15}$,
then 
$\CCG$ 
contains words 
$A\sb 1$, $A\sb 2$, $B\sb 1$, $B\sb 2$
satisfying the conditions above.
Conversely,
from the complete list of geometrically realizable classes of codes,
we see that if $\CCG$ 
contains words 
$A\sb 1$, $A\sb 2$, $B\sb 1$, $B\sb 2$
satisfying the conditions above, then 
$\CCL{15}\leqq [\CCG]$ holds.
Hence we have
$$
{\,\UUU}\sb{2, 6, \CCL{15}}\quad\st\quad {\,\UUU}[2211]\quad\st\quad {\,\UUU}\sb{2, 6, \ge \CCL{15}}.
$$
Therefore ${\,\UUU}\sb{2, 6, \CCL{15}}$
is irreducible and its generic point
coincides with the generic point of ${\,\UUU}[2211]$.
\end{example}
By the same argument,
we obtain Table~\ref{table:generic}
of the pairs of $\CCL{\nu}$ and $[a\sb 1\dots a\sb k]$
such that 
${\,\UUU}\sb{2, 6, \CCL{\nu}}$
is irreducible,
and 
that 
the generic point of ${\,\UUU}\sb{2, 6, \CCL{\nu}}$
coincides with the generic point of ${\,\UUU}[a\sb 1\dots a\sb k]$.
\begin{table}
$$
{
\renewcommand{\arraystretch}{1.2}
\begin{array}{|c|| c| c| c| c| c| c| c| c| c| c|}
\hline
\nu  &  4 & 6 & 8 & 13 & 15 & 35 & 77  \\
\hline
\sigma &  8 & 8 & 8 & 7 & 7 & 6 & 5 \\
\hline
[a\sb 1\dots a\sb k]  &  [411] & [321] & [222] & [3111] & [2211]& [21111]& [111111] \\
\hline
\end{array}
}
$$
\caption{The pairs of $\CCL{\nu}$ and $[a\sb 1\dots a\sb k]$}\label{table:generic}
\end{table}
\begin{example}\label{example:pascal}
Let $G$ be a polynomial of $\Uts$,
and let $A\sb 1, \dots, A\sb 6$ and $B$ be  distinct words of $\CCG$.
We say that $(A\sb 1, \dots, A\sb 6, B)$ is a {\it Pascal configuration}
if the following hold:
\begin{itemize}
\item The words $A\sb 1$, \dots, $A\sb 6$ are  of weight $5$.
\item The word $B$ is of weight $8$ and irreducible in $\CCG$.
\item Let $P\sb{ij}$ be the point of $A\sb{i}\cap A\sb{j}$ for $i\ne j$.
Then the six points $P\sb{12}$, $P\sb{23}$, $P\sb{34}$, $P\sb{45}$, $P\sb{56}$ and  $P\sb{61}$
are distinct and contained in $B$.
\end{itemize}
The code $\CCG$ contains a Pascal configuration if and only if
there exists a hexagon $L\sb 1 L\sb 2 L\sb 3 L\sb 4 L\sb 5 L\sb 6$
formed by lines splitting in $\XG$
that is inscribed in a smooth conic $Q$.
(See Figure~\ref{fig:pascal}.)
Note that the conic $Q$ is also splitting in $\XG$ by Proposition~\ref{prop:convnum}.
\begin{figure}
\includegraphics[height=2.5cm]{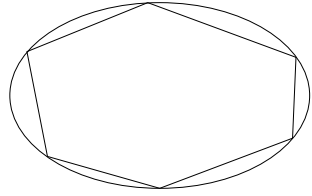}
\caption{The Pascal configuration}\label{fig:pascal}
\end{figure}
If $\CCG$ is in the class $\CCL{136}$,
then $\CCG$ contains a Pascal configuration.
If $\CCG$ contains a Pascal configuration,
then $\CCL{136}\leqq [\CCG]$ holds.
Because the moduli of pairs 
of a smooth conic $Q$  and a hexagon inscribed in $Q$
is irreducible,
we conclude that the locus ${\,\UUU}\sb{2, 6, \CCL{136}}$
is irreducible.
\par
We fix a smooth conic $Q\sb 1\st\Pt$,
and let $P\sb 1, \dots, P\sb 6$
be general points on $Q\sb 1$.
We put
$$
L\sb i:= \overline{P\sb i P\sb {i+1}}\quad(i=1, \dots, 5),\qquad
L\sb 6:= \overline{P\sb 6 P\sb 1}.
$$
Let $G\sb{L\sb i}=0$ be a defining equation of the line $L\sb i$.
Then
$$
G:= G\sb{L\sb 1} G\sb{L\sb 2} G\sb{L\sb 3} G\sb {L\sb 4} G\sb {L\sb 5} G\sb {L\sb 6}
$$
is a point of ${\,\UUU}\sb{2, 6, \CCL{136}}$.
The points
$L\sb 1\wedge L\sb 4$, 
$L\sb 2\wedge L\sb 5$,
and 
$L\sb 3\wedge L\sb 6$
are distinct,
because
$P\sb 1, \dots, P\sb 6$
are general on $Q\sb 1$.
By Pascal's theorem,
these three points are  on a line $M$.
By the converse to Pascal's theorem,
the hexagons
$$
L\sb{1}L\sb{5}L\sb{3} L\sb{4}L\sb{2}L\sb{6},
\quad
L\sb{1}L\sb{2}L\sb{6} L\sb{4}L\sb{5}L\sb{3},
\quand
L\sb{1}L\sb{5}L\sb{6} L\sb{4}L\sb{2}L\sb{3},
$$
are also inscribed in smooth conics.
Let $Q\sb 2$, $Q\sb 3$ and $Q\sb 4$ be those conics.
Then the lines $L\sb 1, \dots, L\sb 6, M$
and the smooth conics $Q\sb 1, \dots, Q\sb 4$ are splitting in $\XG$.
\end{example}
\begin{figure}
\includegraphics[height=3cm]{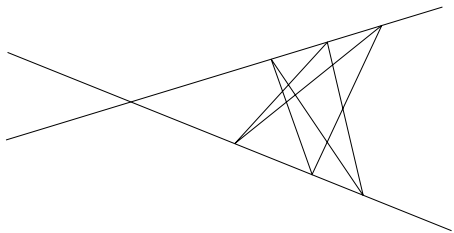}
\caption{The Pappos configuration}\label{fig:pappos}
\end{figure}
\begin{example}\label{example:pappos}
The class $\CCL{177}$ corresponds to the {\it Pappos configuration}
(Figure~\ref{fig:pappos})
in the same way as $\CCL{136}$ corresponds to the Pascal configuration.
Hence ${\,\UUU}\sb{2, 6, \CCL{177}}$
is irreducible.
\quad\end{example}
\bibliographystyle{amsplain}

\providecommand{\bysame}{\leavevmode\hbox to3em{\hrulefill}\thinspace}

\end{document}